\documentclass[preprint,11pt]{elsarticle}
\usepackage[bookmarks=true,colorlinks=true,linkcolor=blue]{hyperref}

\pdfoutput=1

\biboptions{sort&compress}

\usepackage[usenames]{color}

\newcommand{\new}[1]{{\color{blue!50!black}#1}}

\usepackage{soul}

\usepackage{amssymb,amsmath}
\usepackage{amsthm}

\usepackage{calc}
\usepackage[margin=1.in]{geometry}

\renewcommand{\url}[1]{}
\newcommand{\citeCount}[1]{}
\newtheorem{theorem}{Theorem}

\newtheorem{condition}{Condition}

\newtheorem*{algorithmVP}{Velocity Projection Algorithm}

\newcommand{\etah}{\hat{\eta}}

\newcommand{\Hf}{H}

\newcommand{\dt}{\Delta t}

\newcommand{\bogus}[1]{{}}

\newcommand{\fv}{\mathbf{ f}}
\newcommand{\gv}{\mathbf{ g}}

\newcommand{\iv}{\mathbf{ i}}
\newcommand{\jv}{\mathbf{ j}}
\newcommand{\kv}{\mathbf{ k}}

\newcommand{\nv}{\mathbf{ n}}

\newcommand{\rv}{\mathbf{ r}}
\newcommand{\sv}{\mathbf{ s}}
\newcommand{\tv}{\mathbf{ t}}
\newcommand{\uv}{\mathbf{ u}}
\newcommand{\vv}{\mathbf{ v}}
\newcommand{\wv}{\mathbf{ w}}
\newcommand{\xv}{\mathbf{ x}}

\newcommand{\zv}{\mathbf{ z}}

\newcommand{\Fv}{\mathbf{ F}}

\newcommand{\Iv}{\mathbf{ I}}

\newcommand{\Lv}{\mathbf{ L}}

\newcommand{\half}{{1\over2}}

\newcommand{\Real}{{\mathbb R}}

\newcommand{\Gs}{{\mathcal G}}

\newcommand{\Gc}{{\mathcal G}}

\newcommand{\Pc}{{\mathcal P}}

\newcommand{\tauv}{\boldsymbol{\tau}}

\newcommand{\sigmav}{\boldsymbol{\sigma}}
\newcommand{\etav}{\boldsymbol{\eta}}

\newcommand{\grad}{\nabla}

\newcommand{\tableFont}{\footnotesize}
\newcommand{\tableFontSize}{\tableFont}
\newcommand{\num}[2]{#1e#2} 
\newcommand{\errFormat}[1]{$\text{E}_j^{(#1)}$} 

\newcommand{\rateLabel}{rate}

\clearpage

\newcommand{\As}{\bar{A}}
\newcommand{\Pss}{\bar{\mathcal{P}}}
\newcommand{\Ks}{\bar{K}}

\newcommand{\rhos}{\bar{\rho}}

\newcommand{\xs}{\bar{x}}
\newcommand{\xsv}{\bar{\xv}}
\newcommand{\xsbv}{{\xsv_b}}

\newcommand{\usv}{\bar{\uv}}
\newcommand{\vsv}{\bar{\vv}}

\newcommand{\wsv}{\bar{\wv}}

\newcommand{\vvh}{\hat{\vv}}
\newcommand{\vh}{\hat{v}}
\newcommand{\ph}{\hat{p}}

\newcommand{\nsv}{\bar{\nv}}
\newcommand{\vvs}{\bar{\vv}}
\newcommand{\uvs}{\bar{\uv}}

\newcommand{\fsv}{\bar{\fv}}

\newcommand{\amp}{A}

\newcommand{\bfss}{\sffamily\bfseries}

\newcommand{\normalss}{\sffamily}

\newcommand{\hs}{\bar{h}}
\newcommand{\hsHalf}{\frac{\hs}{2}}

\newcommand{\Af}{A_f}

\newcommand{\OmegaF}{\Omega}
\newcommand{\OmegaS}{\bar\Omega}
\newcommand{\OmegaSbar}{\bar\Omega}

\newcommand{\OmegaFh}{\Omega_h}
\newcommand{\OmegaSh}{\bar\Omega_h}

\newcommand{\GammaIh}{\Gamma_h}

\newcommand{\Lt}{\tilde{L}}

\newcommand{\kx}{k_x} 

\newcommand{\tanv}{\tv}

\newcommand{\strutt}{\rule{0pt}{10pt}}

\newcommand{\nd}{{n_d}} 

\newcommand{\Fig}{Figure} 

\newcommand\Ts{{\bar T}}

\newcommand{\xp}{x'}

\newcommand{\Gbcf}{\Gc_{bcf}} 

\newcommand\tnv{{\mathbf{t}}}

\renewcommand\tv{\tnv}

\newcommand{\xb}{x_0}
\newcommand{\xvb}{{\bar\xv_0}}
\newcommand{\xvz}{{\bar\zv}}

\newcommand{\partialt}{\frac{\partial}{\partial t}}

\newcommand{\Es}{\bar{E}}
\newcommand{\Is}{\bar{I}}
\newcommand{\Ls}{\bar{L}}

\newcommand{\xbv}{\xsv_0}

\renewcommand{\ss}{s}
\newcommand{\ls}{\bar{l}}
\newcommand{\Lsv}{\bar{\Lv}}
\newcommand{\Lvs}{\bar{\Lv}}

\newcommand{\ww}{\eta}

\newcommand{\wwv}{\etav}
\newcommand{\Nn}{N_n}

\newcommand{\Gcfc}{\Gc_{\text{et}}}
\newcommand{\pMax}{p_{\text{max}}}
\newcommand{\tMax}{t_{\text{max}}}

\usepackage{tikz}

\newlength{\tfwidth}
\newlength{\tfheight}
\newlength{\tfxa}
\newlength{\tfxb}
\newlength{\tfya}
\newlength{\tfyb}
\newcommand{\trimFigWithBox}[6]{%
\setlength\fboxsep{0pt}%
\setlength\fboxrule{1.0pt}
\fbox{\includegraphics[width=#2, clip, trim=#3 #4 #5 #6]{#1}}%
}
\newcommand{\trimFigNoBox}[6]{%
\setlength\fboxsep{1pt}
\setlength\fboxrule{0.0pt}
\fbox{\includegraphics[width=#2, clip, trim=#3 #4 #5 #6]{#1}}%
}
\newcommand{\trimFigb}[6]{%
\setlength{\tfwidth}{(#2+#2*\real{#3})+#2*\real{#4}}
\setlength{\tfheight}{(#2+#2*\real{#5})+#2*\real{#6}}%
\setlength{\tfxa}{\tfwidth*\real{#3}}%
\setlength{\tfxb}{\tfwidth*\real{#4}}%
\setlength{\tfya}{\tfheight*\real{#5}}%
\setlength{\tfyb}{\tfheight*\real{#6}}%
\trimFigWithBox{#1}{#2}{\tfxa}{\tfya}{\tfxb}{\tfyb}%
}
\newcommand{\trimFig}[6]{%
\setlength{\tfwidth}{(#2+#2*\real{#3})+#2*\real{#4}}
\setlength{\tfheight}{(#2+#2*\real{#5})+#2*\real{#6}}%
\setlength{\tfxa}{\tfwidth*\real{#3}}%
\setlength{\tfxb}{\tfwidth*\real{#4}}%
\setlength{\tfya}{\tfheight*\real{#5}}%
\setlength{\tfyb}{\tfheight*\real{#6}}%
\trimFigNoBox{#1}{#2}{\tfxa}{\tfya}{\tfxb}{\tfyb}%
}
\usepackage{algorithm} 
\usepackage{algpseudocode} 

\begin{document}

\small

\begin{frontmatter}
\title{A stable partitioned FSI algorithm for incompressible flow\\ and deforming beams}

\author[rpi]{L.~Li\fnref{LongfeiThanks}}
\ead{lil19@rpi.edu}

\author[rpi]{W.~D.~Henshaw\corref{cor1}\fnref{DOEThanks,NSFgrantNew}}
\ead{henshw@rpi.edu}

\author[rpi]{J.~W.~Banks\fnref{DOEThanks,PECASEThanks}}
\ead{banksj3@rpi.edu}

\author[rpi]{D.~W.~Schwendeman\fnref{DOEThanks,NSFgrantNew}}
\ead{schwed@rpi.edu}

\author[duke]{G.A.~Main\fnref{CGSF}}
\ead{amain8511@gmail.com}

\address[rpi]{Department of Mathematical Sciences, Rensselaer Polytechnic Institute, Troy, NY 12180, USA.}

\address[duke]{Department of Civil and Environmental Engineering, Duke University, Durham, NC, 27708, USA.}

\cortext[cor1]{Department of Mathematical Sciences, Rensselaer Polytechnic Institute, 110 8th Street, Troy, NY 12180, USA.}


\fntext[LongfeiThanks]{Research supported by the Margaret A. Darrin Postdoctoral Fellowship.}

\fntext[DOEThanks]{This work was supported by contracts from the U.S. Department of Energy ASCR Applied Math Program.}

\fntext[NSFgrantNew]{Research supported by the National Science Foundation under grant DMS-1519934.}


\fntext[PECASEThanks]{Research supported by a U.S. Presidential Early Career Award for Scientists and Engineers.}

\fntext[CGSF]{Research supported by the U.S. Department of Energy NNSA Stewardship Science Graduate Fellowship.}



\begin{abstract}

An added-mass partitioned (AMP) algorithm is described for solving fluid-structure interaction (FSI) problems coupling
incompressible flows with thin elastic structures undergoing finite deformations.  The new AMP scheme
is fully second-order accurate and stable, without sub-time-step iterations, even for very light
structures when added-mass effects are strong.
The fluid, governed by the incompressible Navier-Stokes equations, is solved in velocity-pressure form using 
a fractional-step method; large deformations are treated with a mixed Eulerian-Lagrangian approach on deforming 
composite grids. The motion of the thin structure is governed by a generalized
Euler-Bernoulli beam model, and these equations are solved in a Lagrangian frame using two
approaches, one based on finite differences and the other on finite elements.  The key AMP interface
condition is a generalized Robin (mixed) condition on the fluid pressure.  This condition, which is derived
at a continuous level, has no adjustable parameters and 
is applied at the discrete level to couple the partitioned domain solvers.
Special treatment of the AMP condition is required to couple the finite-element beam solver with the
finite-difference-based fluid solver, and two coupling approaches are described. 
A normal-mode stability analysis is performed for a linearized model problem involving a beam separating two fluid
domains, and it is shown that the AMP scheme is stable independent of the ratio of the mass of the
fluid to that of the structure.  A traditional partitioned (TP) scheme using a Dirichlet-Neumann
coupling for the same model problem is shown to be unconditionally unstable if the added mass of the
fluid is too large.  A series of benchmark problems of increasing complexity are considered to
illustrate the behavior of the AMP algorithm, and to compare the behavior with that of the TP
scheme. The results of all these benchmark problems verify the stability and accuracy of the AMP scheme. 
Results for one benchmark problem modeling blood flow in a deforming artery are also compared
with corresponding results available in the literature.

\end{abstract}


\end{frontmatter}

\clearpage
\tableofcontents

\clearpage
\section{Introduction}


Fluid-structure interaction (FSI) problems that describe the motion of an incompressible fluid coupled to a thin walled structure
(beam or shell) arise in many applications such as those in structural engineering and biomedicine.  Such FSI problems are often
modeled mathematically by suitable partial differential equations for the fluids and structures in their respective domains,
together with matching conditions at the boundaries of the domains where the solutions of the equations interact.  Numerical algorithms used
to solve these FSI problems can be classified into two main categories.  Algorithms belonging to one category, called monolithic
schemes, treat the equations for the fluids and structures along with interface and boundary conditions as a large system of
evolution equations, and then advance the solutions together.  The other category of algorithms are partitioned schemes (also known
as modular or sequential schemes), and these algorithms employ separate solvers for the fluids and structures which are coupled at
the interface.  Sub-iterations are often performed at each time step of partitioned
algorithms for stability.
Even though many existing partitioned schemes suffer from moderate to serious stability issues in certain problem regimes, they are
often preferred since they can make use of existing solvers and can be more efficient than monolithic schemes.

The traditional partitioned algorithm for beams (or shells) uses the velocity and/or acceleration of the
solid as a boundary condition on the fluid.  The force of the fluid is accounted for through a body forcing on the beam.  It has
been found that partitioned schemes may be unstable, or require multiple sub-iterations per time step, when the density of the
structure is similar to or lighter than that of the fluid~\cite{CausinGerbeauNobile2005,vanBrummelen2009}.  These instabilities are
attributed to the {\em added-mass effect} whereby the force required to accelerate a structure immersed in a fluid must also account
for accelerating the surrounding fluid.  The added-mass effect has been found to be especially problematic in many biological flows
such as haemodynamics since the density of the fluid (blood) is similar to that of the adjacent structure (arterial
walls)~\cite{YuBaekKarniadakis2013}.

In our previous companion papers~\cite{fib2014,fis2014}, we developed stable added-mass partitioned (AMP) algorithms
for linearized FSI problems involving incompressible Stokes fluids coupled to elastic bulk solids and to beams. 
The key ingredient of our approach is the use of non-traditional Robin interface conditions.
These condition are derived at a continuous level, have no adjustable parameters, and are amenable for
incorporation in high-order accurate schemes.
Using mode analysis, the AMP approach was shown to be stable without sub-iterations per time step, and the numerical results 
demonstrated second-order accuracy in the max-norm. 
In the present paper, we focus on the case of thin structures and our primary purpose is to extend the AMP scheme in~\cite{fis2014}
to the case of finite amplitude motions in more general geometries. 
In particular, the AMP interface conditions are derived for beams of finite thickness,
for beams with fluid on two sides, and for beams with a free end immersed in the fluid.
Finite amplitude motions of the beam results in finite deformations of the surrounding fluid domain.  Our numerical
approach for the solution of the fluid equations in an evolving fluid domain is based on the use of deforming composite
grids (DCG).  The DCG approach for FSI problems was described first in~\cite{fsi2012} inviscid
compressible flow coupled to a linearly elastic solid, and later 
in~\cite{flunsi2014r} for compressible flow coupled to nonlinear hyperelastic solids.  The approach is extended
here for incompressible flow coupled to beams.
Both finite-difference and finite-element approximations to the equations governing the motion of the beam are considered, and
issues concerning the interface coupling of a finite-element based beam solver to a finite-difference based fluid solver are discussed.
The stability analysis in~\cite{fis2014} is also extended to treat beams with fluid on two sides.
Finally, several benchmark FSI problems are presented to illustrate the stability and accuracy of the present AMP algorithm.









The investigation of FSI problems and the development of numerical approximations for their solution
are very active areas of research, see for example~\cite{FernandezReview2011,Degroote2013,BukacCanicMuha2015} and the
references therein.
For the coupling of incompressible flows and bulk solids, the work in~\cite{NobilePozzoliVergara2013,YuBaekKarniadakis2013,LiuJaimanGurugubelli2014,fib2014,BukacCanicMuha2015} describes recent developments.
The development of partitioned schemes for the case of incompressible fluids coupled to thin structures, the focus of this paper, 
is also an active area. 
%
The first stable partitioned scheme for incompressible flow coupled to thin structures
was the ``kinematically coupled scheme'' of Guidoboni et~al.~\cite{GuidoboniGlowinskiCavalliniCanic2009}
which uses an operator splitting of the kinematic interface condition (matching the fluid and beam velocities). 
This scheme was further advanced in~\cite{CanicMuhaBukac2012,BukacCanicGlowinskiTambacaQuaini2013}.
Luk\'a\v{c}ov\'a-Medvid'ov\'a et~al.~\cite{LukacovaMedvidovaRusnakovaHundertmarkZauskova2013}  developed
a partitioned scheme based on the kinematically
coupled scheme that uses Strang splitting.
Fernandez and collaborators have also developed  stable 
partitioned schemes 
for incompressible flow coupled to
thin structures that are based on 
a time-splitting of the kinematic boundary condition~\cite{Fernandez2011,Fernandez2012,FernandezMullaertVidrascu2013,FernandezLandajuela2013Note,FernandezLandajuela2014,FernandezLanajuelaVidrascu2015}. 
However, it appears to be difficult to achieve higher than first-order accuracy with the time-splitting approach. 
The AMP approximation, in contrast, is not based on a time splitting and is amenable to second- or even higher-order accuracy. Second-order accuracy in the max-norm was demonstrated in~\cite{fis2014} for linearized problems, while in this paper second-order accuracy is shown for the full nonlinear problem with deforming domains.

The remainder of the paper is organized as follows. In Section~\ref{sec:governingEquations} we describe the governing equations.
The AMP interface conditions are derived in Section~\ref{sec:AMPalgorithm}.
A second-order accurate predictor-corrector
algorithm based on these conditions is described in Section~\ref{sec:timeStepping}.
A specific choice for the beam model is given in Section~\ref{sec:beamModels}.
In Section~\ref{sec:analysis}, the stability of the AMP scheme is shown for 
a linearized FSI problem involving a beam with fluid on two sides.  
The numerical approach used for moving domains based on deforming composite grids, our numerical
approaches for the beam solver, and the treatment of the AMP interface conditions
are described in Section~\ref{sec:numericalApproach}. 
Numerical results presented in Section~\ref{sec:numericalResults}
carefully demonstrate the stability and accuracy of the AMP algorithm.
Conclusions are provided in Section~\ref{sec:conclusions}.


\renewcommand{\sv}{{s}}
\renewcommand{\OmegaSbar}{{\bar{\cal R}}}

\section{Governing equations}\label{sec:governingEquations}

We consider the fluid-structure coupling of an incompressible fluid and a thin deformable structure.  The fluid occupies the domain $\xv\in\OmegaF(t)$ while the structure lies in the domain $\xv\in\OmegaS(t)$, where $\xv$ is position and $t$ is time.  The coupling of the fluid and structure occurs along the interface $\Gamma(t)$, see Figure~\ref{fig:threeDimensionalBeam}.  It is assumed that the fluid is governed by the incompressible Navier-Stokes equations, which in an Eulerian frame are given by 
\begin{alignat}{3}
  &  \frac{\partial\vv}{\partial t} + (\vv\cdot\grad)\vv 
                 =  \frac{1}{\rho} \grad\cdot\sigmav  , \qquad&& \xv\in\OmegaF(t) ,  \label{eq:fluidMomentum}  \\
  & \grad\cdot\vv =0,  \quad&& \xv\in\OmegaF\new{(t)} ,  \label{eq:fluidDiv3d}
\end{alignat}
where $\rho$ is the (constant) fluid density and $\vv=\vv(\xv,t)$ is the fluid velocity.  The fluid stress tensor, $\sigmav=\sigmav(\xv,t)$, is given by 
\begin{equation}
  \sigmav = -p \Iv + \tauv,\qquad \tauv = \mu \left[ \grad\vv + (\grad\vv)^T \right],
  \label{eq:fluidStress}
\end{equation}
where $p=p(\xv,t)$ is the pressure, $\Iv$ is the identity tensor, $\tauv$ is the viscous stress tensor, and $\mu$ is the (constant) fluid viscosity. For future reference, the components of a vector such as $\vv$ will be denoted by $v_m$, $m=1,2,3$ (i.e.~$\vv=[v_1, v_2, v_3]^T$), while components of a tensor such as $\sigmav$ will be denoted by $\sigma_{mn}$, $m,n=1,2,3$.  The velocity-divergence form of the equations given by~\eqref{eq:fluidMomentum} and~\eqref{eq:fluidDiv3d} require appropriate initial and boundary conditions, as well as conditions on $\Gamma(t)$ where the behavior of the fluid is coupled to that of the solid (as discussed below).

%
%
%

{
\def\xa{-.5}
\def\xb{12.5}
\def\xs{.3}
\def\ys{.3}
\def\lb{8}
\def\a{.075}%
\def\b{45}
\def\hsh{.5}
\def\extra{1.}
\def\yTop{2.5}
\def\yBot{-2.75}
\begin{figure}[hbt]
\newcommand{\textFont}{\normalss}
\begin{center}
\begin{tikzpicture}[scale=1]
\useasboundingbox (\xa,.25) rectangle (12.5,5.25);  
\begin{scope}[xshift=0cm,yshift=\yTop cm]
 \fill[fill=blue!20,draw=blue!20,line width=0pt,yshift=4pt] 
        (\xa,\yTop) -- (\xb,\yTop) -- (\xb,\yBot) -- (\xa,\yBot) -- cycle;
  \begin{scope}[xshift=1.5cm,rotate=-5]
    %
    %
    \fill[fill=red!20,draw=red,line width=2pt,xshift=.3cm,yshift=.3cm] 
           plot[samples=100, domain=0.:\lb] (\x, {(\hsh+\a*cos(\b*\x))}) --
           plot[samples=100, domain=\lb:0] (\x, {(-\hsh+\a*cos(\b*\x))}) -- cycle;
    %
    %
    \fill[fill=red!20,draw=red,line width=2pt,xshift=.3cm,yshift=.3cm] 
           plot[samples=100, domain=0.:\lb] (\x, {(-\hsh+\a*cos(\b*\x))}) --
           plot[samples=100, domain=\lb:0] (\x-2*\xs, {(-2*\ys-\hsh+\a*cos(\b*\x))}) -- cycle;
    \foreach \hShift in {-\hsh,\hsh}
    { 
      \draw[draw=red,line width=2pt,xshift=.3cm,yshift=.3cm,dotted] 
               plot[samples=100, domain=-\extra:\lb+\extra] (\x, {(\hShift+\a*cos(\b*\x))});
    }
    \foreach \hShift in {-\hsh,\hsh}
    { 
      \draw[draw=red,line width=2pt,xshift=-.3cm,yshift=-.3cm,dotted] 
               plot[samples=100, domain=-\extra:\lb+\extra] (\x, {(\hShift+\a*cos(\b*\x))});
    }
    %
    %
    \foreach \xx in {0,4,8}
    {
      \draw[-,red,fill=red!20,line width=2pt,fill opacity=0.5] (\xx-\xs, {( \hsh+\a*cos(\b*\xx)-\ys)}) -- (\xx+\xs, {( \hsh+\a*cos(\b*\xx)+\ys)}) --
                                  (\xx+\xs, {(-\hsh+\a*cos(\b*\xx)+\ys)}) -- (\xx-\xs, {(-\hsh+\a*cos(\b*\xx)-\ys)}) -- cycle;
    }
    \draw[-,black,line width=2pt] plot[samples=100, domain=-\hsh-\hsh:\lb+\hsh+\hsh] (\x, {\a*cos(\b*\x)}); 
    %
    \fill[fill=red!20,draw=red,line width=2pt,xshift=-.3cm,yshift=-.3cm,fill opacity=0.4] 
           plot[samples=100, domain=0.:\lb] (\x, {(\hsh+\a*cos(\b*\x))}) --
           plot[samples=100, domain=\lb:0] (\x, {(-\hsh+\a*cos(\b*\x))}) -- cycle;
    %
    %
    %
    \foreach \xx in {0,4,8}
    {
      \draw[-,red,line width=2pt] (\xx-\xs, {( \hsh+\a*cos(\b*\xx)-\ys)}) -- (\xx-\xs, {(-\hsh+\a*cos(\b*\xx)-\ys)});
    }
    %
    \draw[black] (5.25,1.5) node[anchor=south,xshift=9pt,yshift=-2pt] {$\xbv(\ss,t)$};
    \draw[->,black,line width=1.pt] (5.25,1.5) -- (4,-.02);    
    \newcommand{\xx}{3.75}
    \draw[->,black,line width=1pt] (\xx+\xs, {( \hsh+\a*cos(\b*\xx)+\ys)}) -- (\xx-\xs, {( \hsh+\a*cos(\b*\xx)-\ys)}) --
                                  (\xx-\xs, {(-\hsh+\a*cos(\b*\xx)+\ys)}) node[anchor=north] {$\theta$};
    %
    \draw[black] (4.75,-1.25) node[anchor=north,xshift=5pt] {\textFont $\xsv_b(\theta,s,t)$};
    \draw[->,black,line width=1.pt] (4.75,-1.25) -- (3.7,-.3); 
    %
    \draw[black] (1.,-1.3) node[anchor=north,xshift=-3pt,yshift=0pt] {\Large$\OmegaS(t)$}; 
    \draw[->,black,line width=1.pt] (1,-1.3) -- (1.8,-.5); 
    \draw[black] (2,1.5) node[anchor=south,xshift=0pt,yshift=-3pt] {\Large$\Gamma(t)$}; 
    \draw[->,black,line width=1.pt] (2,1.5) -- (2.5,.5); 
  \end{scope}
  \draw[black] (11,2) node[anchor=east,xshift=0pt,yshift=-20pt] {\Large$\OmegaF(t)$}; 
\end{scope}
%
%
\end{tikzpicture}
\end{center}
\caption{Geometry of a three-dimensional beam immersed in a fluid. The solid domain is $\OmegaS(t)$. The fluid domain
is $\OmegaF(t)$. The interface between the fluid and solid is $\Gamma(t)$. 
The beam reference-line is $\xvb(s,t)$. 
A point on the surface of the beam is denoted by $\xsv_b(\theta,s,t)$ where $\theta\in\Pc(s)$ is the parameterization
of the cross-sectional surface curve.
The fluid tractions on the surface points $\xsv_b(\theta,s,t)$, for $\theta\in\Pc(s)$,  contribute to the force
on the beam at $\xvb(s,t)$.
 }
\label{fig:threeDimensionalBeam}
\end{figure}
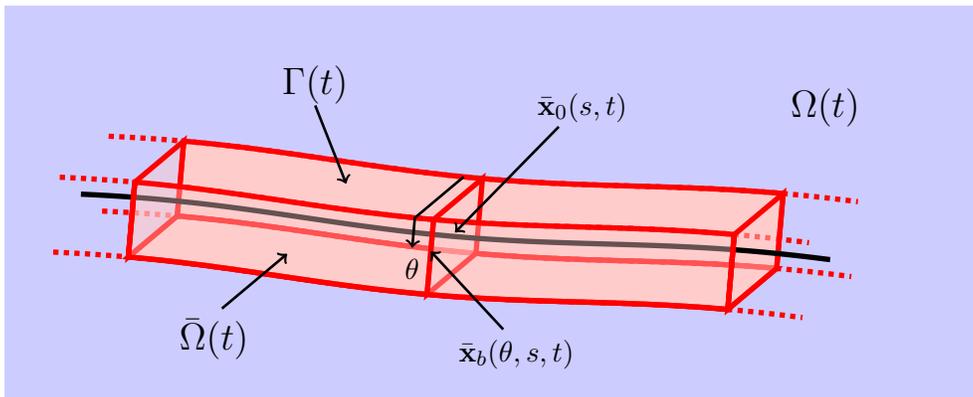
}

An elliptic equation for the fluid pressure can be derived from~\eqref{eq:fluidMomentum}--\eqref{eq:fluidStress} as
\begin{alignat}{3}
  &  \Delta p &= - \rho \grad\vv:(\grad\vv)^T ,  \qquad&& \xv\in\OmegaF(t) . \label{eq:fluidPressure}
\end{alignat}
In the velocity-pressure form of the incompressible Navier-Stokes equations, the Poisson equation~\eqref{eq:fluidPressure} is used in place of~\eqref{eq:fluidDiv3d}. This form requires an additional boundary condition, and a natural choice is $\grad\cdot\vv=0$ for $\xv\in\partial\OmegaF(t)$, see~\cite{ICNS} for example.

The structure is assumed to be thin and governed by a beam model, which describes the evolution of the displacement, $\uvs=\uvs(\sv,t)$, of a reference centerline curve of the beam as a function of a Lagrangian coordinate $\sv\in\OmegaSbar$ and time $t$.  (Overbars are used throughout to denote quantities belonging to the structure.)  The beam model has the form
\begin{alignat}{3}
  \rhos \As \frac{\partial^2 \uvs}{\partial t^2} =& \Lvs(\uvs,\vvs) +\fsv(\sv,t), \qquad && \sv\in\OmegaSbar,
     \label{eq:beamEquation}
\end{alignat}
where $\rhos$ is the (constant) beam density, $\As=\As(\sv)$ is the cross-sectional area of the beam, $\Lvs$ is the internal
restoring force per unit length which depends on $\uvs$ and on the velocity, $\vvs=\partial \uvs/\partial t$, of the reference curve
(in general), and $\fsv$ is the external force per unit length on the beam from the fluid (defined below).  The effects of torsion
are neglected in the present beam model.  In physical space, the position $\xvb$ of the reference curve
is given by
\begin{equation}
\xvb(\sv,t)=\xvb(\sv,0)+\uvs(\sv,t),\qquad \sv\in\OmegaSbar,
\label{eq:xcenterline}
\end{equation}
where $\xvb(\sv,0)$ is the initial position of the reference curve of the beam.  The position $\xsbv$ of the surface of the beam is given by
\begin{equation}
\xsbv(\theta,\sv,t)=\xvb(\sv,t)+\xvz(\theta,\sv,t),\qquad \theta\in\bar{\mathcal{P}},\quad\sv\in\OmegaSbar, 
\label{eq:surfaces}
\end{equation}
where $\xvz$ is distance from the reference curve to the perimeter $\bar{\mathcal{P}}=\bar{\mathcal{P}}(\sv)$ of the cross section of the beam at $\sv$, parameterized by arclength $\theta$ along the perimeter. 
Note that $\xvz$ is defined in terms of the beam variables; 
its form for an Euler-Bernoulli beam is
given in Section~\ref{sec:beamModels}.
 The surface of the beam given by~\eqref{eq:surfaces} defines the position of the fluid-structure interface $\Gamma(t)$ (see Figure~\ref{fig:threeDimensionalBeam}).

At the fluid-structure interface, the velocity of the fluid at a point on the surface of the beam matches the corresponding velocity of the beam surface,
\begin{align}
    \vv(\xsbv(\theta,\sv,t),t) &= \vsv_b(\theta,\sv,t), \qquad \theta\in\bar{\mathcal{P}},\quad\sv\in\OmegaSbar, \label{eq:kinematicII}
\end{align}
where 
\begin{align}
   \vsv_b(\theta,\sv,t)=  \frac{\partial}{\partial t}\xsbv(\theta,\sv,t) = \vsv(\sv,t) + \wsv(\theta,\sv,t) ,\qquad \wsv(\theta,\sv,t)=\frac{\partial}{\partial t}\xvz(\theta,\sv,t).\label{eq:beamSurfaceVelocity}
\end{align}
Formulas for the finite-thickness corrections, $\xvz$ and $\wsv$, are given in Section~\ref{sec:beamModels} for the case of an Euler-Bernoulli beam model.  The dynamic matching condition (balance of forces) defines the force per unit length on the beam, $\fsv$, in (\ref{eq:beamEquation}) in terms of the fluid force on the surface of the beam,
\begin{align}
    \fsv(\sv,t)  &= -\int_{{\bar{\mathcal{P}}}}(\sigmav\nv)(\hat\theta,\sv,t)\,d\hat\theta, \qquad \sv\in\OmegaSbar , \label{eq:dynamic} 
\end{align} 
where
\begin{align}
    (\sigmav\nv)(\theta,\sv,t) &\equiv \sigmav(\xsbv(\theta,\sv,t),t)\,\nv(\xsbv(\theta,\sv,t),t).
\end{align} 
Here, $\nv$ denotes the unit normal to the beam surface (outward pointing from the fluid domain). 

%
%
%

{
\def\xa{0}
\def\xb{13}
\def\lb{8}
\def\a{.075}%
\def\b{45}
\def\hsh{.5}
\def\yTop{2.5}
\def\yBot{-2.5}
\begin{figure}[hbt]
\newcommand{\textFont}{\normalss}
\begin{center}
\begin{tikzpicture}[scale=1]
\useasboundingbox (0,.5) rectangle (13,5.25);  
\begin{scope}[xshift=0cm,yshift=\yTop cm]
 \fill[fill=blue!20,draw=blue,line width=2pt,yshift=4pt] 
        (\xa,\yTop) -- (\xb,\yTop) -- (\xb,\yBot) -- (\xa,\yBot) -- cycle;
\fill[fill=red!20,draw=red,line width=2pt] 
       plot[samples=100, domain=0.:\lb] (\x, {(\hsh+\a*cos(\b*\x))}) --
       plot[samples=100, domain=90:-90] ({\lb + cos(\x)},{\a+\hsh*sin(\x)}) --  
       plot[samples=100, domain=\lb:0] (\x, {(-\hsh+\a*cos(\b*\x))}) -- cycle;
%
\draw[-,black,line width=2pt] plot[samples=100, domain=0.:\lb+\hsh+\hsh] (\x, {\a*cos(\b*\x)}); 
\def\xA{1.5} \pgfmathparse{\hsh+\a*cos(\b*\xA)} \let\yB\pgfmathresult
 \draw[black] (\xA,\yB) node[anchor=south] {\textFont $\xv_{b,+}(\sv,t)$};
\def\xA{1.5} \pgfmathparse{-\hsh+\a*cos(\b*\xA)} \let\yA\pgfmathresult
 \draw[black] (\xA,\yA) node[anchor=north] {\textFont $\xv_{b,-}(\sv,t)$};
%
\draw[black] (3.25,-1) node[anchor=north] {\textFont $\xbv(\sv,t)$};
\draw[->,black,line width=1.pt] (3.25,-1) -- (3.5,-.1); 
\draw[black] (11,2) node[anchor=east,xshift=0pt,yshift=-20pt] {\Large$\OmegaF(t)$}; 
\draw[black] (3,1) node[anchor=south,xshift=10pt,yshift=-4pt] {\Large$\OmegaS(t)$}; 
\draw[->,black,line width=1.pt] (3,1) -- (2.75,.25);
\draw[black] (8,-1.) node[anchor=north west,xshift=-2pt,yshift=2pt] {\large$\Gamma(t)$}; 
\draw[->,black,line width=1.pt] (8,-1) -- (7.5,-.5); 
%
\def\xA{5.5} \pgfmathparse{-\hsh+\a*cos(\b*\xA)} \let\yA\pgfmathresult
\def\xB{5.42}  \pgfmathparse{\hsh+\a*cos(\b*\xB)} \let\yB\pgfmathresult
 \draw[-,black] (\xA,\yA) -- (\xB,\yB);
 \filldraw[fill=white,draw=black] (\xA,\yA) circle (3pt);
 \filldraw[fill=white,draw=black] (\xB,\yB) circle (3pt);
 \draw[black] (\xA,\yB) node[anchor=south] {\textFont $(\sigmav\nv)_+$};
 \draw[black] (\xA,\yA) node[anchor=north] {\textFont $(\sigmav\nv)_-$};
%
%
%
\end{scope}
%
%
\end{tikzpicture}
\end{center}
\caption{Geometry of a two-dimensional beam immersed in a fluid. 
The beam reference surface is $\xvb(\sv,t)$. The top and bottom surfaces of the beam are 
the $\xsv_+(\sv,t)$ and $\xsv_-(\sv,t)$, respectively. 
 }
\label{fig:genericBeamCartoon}
\end{figure}
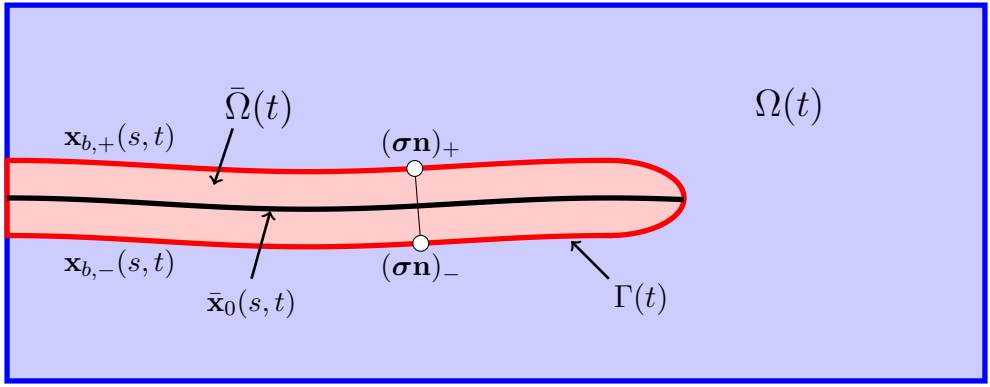
}
\newcommand{\sbf}{\mathbf{s}}

We note that for the case of a beam in two dimensions as shown in Figure~\ref{fig:genericBeamCartoon},
the surface position in~\eqref{eq:surfaces} reduces to
\begin{equation}
\xsv_{b,\pm}(\sv,t)=\xvb(\sv,t)+\xvz_\pm(\sv,t),\qquad \sv\in\OmegaSbar,  
\label{eq:surfaces2d}
\end{equation}
where the plus and minus signs denote the upper and lower surfaces of the beam, respectively.  The velocity on the surfaces of the beam are now given by
\begin{equation}
   \vsv_{b,\pm}(\sv,t)=  \frac{\partial}{\partial t}\xsv_{b,\pm}(\sv,t) = \vsv(\sv,t) + \wsv_\pm(\sv,t) ,\qquad \wsv_\pm(\sv,t)=\frac{\partial}{\partial t}\xvz_\pm(\sv,t),
\label{eq:beamSurfaceVelocity2d}
\end{equation}
and if fluid regions exist on both surfaces of the beam, then the the external force is given by
\begin{equation}
\fsv(\sv,t)  = \Big[ (\sigmav\nv)_+(\sv,t)+ (\sigmav\nv)_-(\sv,t) \Big], \qquad \sv\in\OmegaSbar , \label{eq:dynamic2d} 
\end{equation} 
where
\begin{align}
    (\sigmav\nv)_\pm(\sv,t) &\equiv \sigmav(\xsv_{b,\pm}(\sv,t),t)\,\nv(\xsv_{b,\pm}(\sv,t),t). \label{eq:beam2D_sigma}
\end{align} 
The formulas~\eqref{eq:surfaces2d}--\eqref{eq:beam2D_sigma} also hold for the case of a thin plate
or shell provided the dependent variables, $\xsv_0(\sbf,t)$, $\vsv(\sbf,t)$, etc., are taken to depend on two
independent parameters $(s_1,s_2)=\sbf\in\Real^2$ that parameterize the reference surface of the shell.

%

\section{Added-Mass Partitioning} \label{sec:AMPalgorithm}

In the traditional approach to partitioned fluid-structure interaction, the kinematic condition~\eqref{eq:kinematicII}
is satisfied by imposing the motion of the solid as an interface condition on the fluid. This inherently
one-sided approach is one root cause leading to added-mass instabilities for partitioned fluid-structure
interaction. Various avenues to
more symmetric imposition of~\eqref{eq:kinematicII} have been taken with various degrees of
success. Here we pursue the Added Mass Partitioned (AMP) approach to coupling. 
The core AMP algorithm for coupling a linear incompressible fluid with a thin shell or beam
was described in~\cite{fis2014}. In that work a simplified model was considered where 
linearized equations of motion were used to describe infinitesimal amplitude perturbations
from the specific geometry of fluid in a box with a beam on one of the faces of the box.
This manuscript considers the more general case of the nonlinear Navier-Stokes equations
with finite amplitude deformations applied to more general geometries. 
%
%

\subsection{Continuous AMP Interface Condition}
The AMP interface condition can be derived at a continuous level by first 
taking the time derivative of the kinematic interface condition~\eqref{eq:kinematicII} and using
the definition of the velocity on the surface of the beam~\eqref{eq:beamSurfaceVelocity} to give
\begin{align}
  \frac{D}{Dt}\vv(\xsbv(\theta,\sv,t),t)  & = \partialt \vsv(\sv,t) + \partialt\wsv(\theta,\sv,t)  \nonumber\\
     & = \frac{\partial^2 \uvs}{\partial t^2}(\sv,t) + \partialt\wsv(\theta,\sv,t), \label{eq:interfaceAcceleration}
\end{align}
where $D/Dt=\partial_t+\vv\cdot\grad$ is the material derivative.
The fluid acceleration, $D\vv/Dt$, is given by the fluid momentum equation~\eqref{eq:fluidMomentum},
while the acceleration of the structure, $\partial^2 \uvs / \partial t^2$, is given by the beam equation~\eqref{eq:beamEquation}. Substituting
these definitions into~\eqref{eq:interfaceAcceleration} gives 
\begin{equation}
  \frac{\rhos\As}{\rho}\grad\cdot\sigmav(\xsbv(\theta,\sv,t),t)  = \Lsv\bigl(\uvs,\vsv\bigr) + \rhos\As \partialt \wsv(\theta,\sv,t)  + \fsv(\sv,t),
          \qquad \theta\in\Pss,\quad \sv\in\OmegaSbar.
    \label{eq:AMPconditionSigma}
\end{equation}
%
Finally, using~\eqref{eq:dynamic} for the force $\fsv$ yields the continuous 
compatibility interface condition 
\begin{equation}
     \int_{{\Pss}}(\sigmav\nv)(\hat\theta,\sv,t)\,d\hat\theta + 
\frac{\rhos\As}{\rho}\grad\cdot\sigmav(\xsbv(\theta,\sv,t),t)  
         = \Lsv\bigl(\uvs,\vsv\bigr) + \rhos\As \partialt \wsv(\theta,\sv,t) , \qquad\theta\in\Pss ,\quad \sv\in\OmegaSbar,
    \label{eq:AMP_ContinuousConditions}
\end{equation}
which can be interpreted as a mixed Robin-type condition along the surface of the beam coupling the behaviors of the fluid and structure.

\subsection{Partitioned Schemes Using the AMP Condition}
%
Following the approach used in~\cite{fis2014}, we employ an AMP interface condition
based on the compatibility condition~\eqref{eq:AMP_ContinuousConditions} 
as a stable condition for partitioned schemes coupling incompressible fluids and thin structural beams.
Equation~\eqref{eq:AMP_ContinuousConditions} appears to be a condition that fully couples
the fluid and solid values on the interface. However, a critical observation made in~\cite{fis2014} is
that~\eqref{eq:AMP_ContinuousConditions}
can be used as a stable condition for partitioned schemes by using {\em predicted} values
for the solid variables $\usv^{(p)}$, $\vsv^{(p)}$ and $\wsv^{(p)}$ when evaluating the right-hand-side 
of~\eqref{eq:AMP_ContinuousConditions}.  
These predicted values for the solid
are obtained in a first step of the partitioned scheme and 
the interface condition then becomes a mixed Robin-type condition for the fluid in a subsequent step of the scheme.  
This leads to the primary AMP interface condition:

\begin{condition}
The AMP interface condition for the fluid is 
\begin{equation}
\int_{{\Pss}}(\sigmav\nv)(\hat\theta,\sv,t)\,d\hat\theta + 
\frac{\rhos\As}{\rho}\grad\cdot\sigmav(\xsbv(\theta,\sv,t),t)  
         = \Lsv^{(p)} + \rhos\As \partialt \wsv^{(p)}(\theta,\sv,t), \qquad \theta\in\Pss,\quad \sv\in\OmegaSbar,
    \label{eq:AMP_CombinedConditions}
\end{equation}
where $\usv^{(p)}$, $\vsv^{(p)}$ and $\wsv^{(p)}$ are predicted solid variables, and $\Lsv^{(p)}\equiv \Lsv\bigl(\uvs^{(p)},\vsv^{(p)}\bigr)$.
\end{condition}

The essential philosophy in this approach is to embed the local behavior of solutions of the FSI problem
to define the interface dynamics.
Therefore, in a partitioned scheme the solid would be advanced first to a new time step, providing the predicted
values used in the interface conditions~\eqref{eq:AMP_CombinedConditions}
for the subsequent fluid solve. This process can be viewed as 
computing the solution of a generalized fluid-solid Riemann problem at the interface, thus extending the approach used in the AMP
scheme for elastic solids coupled to inviscid compressible fluids discussed in~\cite{BanksSjogreen2011,fsi2012,flunsi2014r}.
In these previous studies, the FSI problems were purely hyperbolic and the fluid-solid Riemann problem involved local, piecewise constant states of the fluid and solid on either side of the interface.  The solution was then computed using the method of characteristics, giving interface values that are impedance-weighted averages of the fluid and solid states. 
In the present case, involving an incompressible fluid, the initial states 
for the solid domain are $\usv^{(p)}$, $\vsv^{(p)}$ and $\wsv^{(p)}$, while the solution of the generalized
fluid-solid Riemann problem requires the solution of an elliptic problem in the entire fluid domain with~\eqref{eq:AMP_CombinedConditions}  providing the appropriate condition along the interface.



If a fractional-step scheme is used to integrate the fluid equations, 
with the velocity and pressure being updated in separate stages,
then suitable boundary conditions can be derived by
decomposing~\eqref{eq:AMP_CombinedConditions} into normal and tangential components.  The normal component gives a boundary condition for the pressure, while the tangential components, along with the continuity equation~\eqref{eq:fluidDiv3d}, provide boundary conditions for the velocity. 
\begin{condition}
The AMP interface condition for the pressure equation~\eqref{eq:fluidPressure} is given by
\begin{equation}
\nv\sp{T}\int_{{\Pss}}(p\nv)(\hat\theta,\sv,t)\,d\hat\theta + \frac{\rhos\As}{\rho} \frac{\partial p}{\partial n}
=\nv\sp{T}\left[-\Lsv^{(p)} - \rhos\As \partialt \wsv^{(p)}
          +   \frac{\mu\rhos\As}{\rho}\Delta \vv + \int_{{\Pss}}(\tauv\nv)(\hat\theta,\sv,t)\,d\hat\theta\right],
\label{eq:AMPpressureBCI}
\end{equation}
for $\theta\in\Pss$ and $\sv\in\OmegaSbar$, where
\[
(p\nv)(\theta,\sv,t) \equiv p(\xsbv(\theta,\sv,t),t)\,\nv(\xsbv(\theta,\sv,t),t),\qquad
(\tauv\nv)(\theta,\sv,t) \equiv \tauv(\xsbv(\theta,\sv,t),t)\,\nv(\xsbv(\theta,\sv,t),t).
\]
\end{condition}

\bigskip
\begin{condition}
The AMP interface conditions for the velocity when integrating the momentum equation~\eqref{eq:fluidMomentum}
are given by
\begin{equation}
\tanv_{m}^T\left[\int_{{\Pss}}(\tauv\nv)(\hat\theta,\sv,t)\,d\hat\theta+\frac{\mu\rhos\As}{\rho}\Delta \vv\right]=
\tanv_{m}^T\left[\frac{\rhos\As}{\rho}\grad p  + \Lsv^{(p)}+  \rhos\As \partialt \wsv^{(p)}
            + \int_{{\Pss}}(p\nv)(\hat\theta,\sv,t)\,d\hat\theta \right] ,
\label{eq:AMPtangentialVelocityBCI}
\end{equation}
for $\theta\in\Pss$ and $\sv\in\OmegaSbar$, where $\tanv_{m}=\tanv_{m}(\theta,\sv,t)$, $m=1,2$,  denote mutually orthonormal tangent vectors. 
\end{condition}


We note that if the beam motion only depends on one degree of freedom (e.g.~as with 
the Euler-Bernoulli beam model described in Section~\ref{sec:beamModels})
then only~\eqref{eq:AMPpressureBCI} is 
used and~\eqref{eq:AMPtangentialVelocityBCI} should be replaced with a
standard condition on the tangential velocity. For example, \eqref{eq:AMPtangentialVelocityBCI}
can be replaced by setting the tangential components of the fluid velocity equal to the
tangential components of the velocity on the beam surface, i.e.~a {\em no-slip} type condition on the tangential velocity,
\begin{equation}
  \tv_m^T\vv(\xsbv(\theta,\sv,t),t) = \tv_m^T\vsv(\sv,t)+\tv_m^T\wsv(\theta,\sv,t), \qquad m=1,2,
    \label{eq:TangentialVelocityNoSlip}
\end{equation}
for $\theta\in\Pss$ and $\sv\in\OmegaSbar$.  Another choice for the boundary condition could be a {\em slip wall} condition which sets the tangential components of the traction to zero,
\begin{equation}
  \tv_m^T(\tauv\nv)(\theta,\sv,t) = 0 ,\qquad m=1,2,
    \label{eq:TangentialVelocitySlip}
\end{equation}
for $\theta\in\Pss$ and $\sv\in\OmegaSbar$.

\subsection{Velocity projection} 
The velocity at the interface in the AMP scheme is defined to be a weighted average of the predicted velocities from
the fluid and solid. 
This averaging can be motivated by noting that 
for heavy beams the motion of the interface is primarily determined from the beam and thus the interface velocity should
be nearly equal to that of the structure.  For very light beams, on the other hand, 
the fluid interface
approaches a free surface where the motion of the interface is primarily determined by the fluid itself.
The weighted average given in~\cite{fis2014} correctly accounts for these two limits, and this average is extended here
to account for beams of finite thickness.

\begin{algorithmVP}\small
~

\begin{enumerate}
\item
Compute the projected surface velocity $\vv^I(\theta,\sv,t)$ on the surface of the beam as a weighted average of the predicted surface velocities of the fluid $\vv_b\sp{(p)}(\theta,\sv,t)$ and the beam $\vsv_b\sp{(p)}(\theta,\sv,t)$,
\begin{align}
    \vv^I(\theta,\sv,t) = \gamma \vv_b\sp{(p)}(\theta,\sv,t) + (1-\gamma) \vsv_b\sp{(p)}(\theta,\sv,t) ,  \qquad \theta\in\Pss,\quad \sv\in\OmegaSbar,  \label{eq:velocityProjection}
\end{align}
where
\begin{equation}
\gamma=\frac{1}{1 + (\rhos\As)/(\rho \Af)}.
\label{eq:vaverage}
\end{equation}
Here, $\Af=\Af(\sv,t)$ is the area of a characteristic annular region in the fluid surrounding the beam at a position $\xvb(\sv,t)$ on the reference curve.  It has been found that algorithm is insensitive to the choice of this characteristic area~\cite{fis2014}.
\item
Transfer the interface velocities $\vv^I$ (minus the surface correction $\wsv$) onto the beam reference curve using an average of the values on the
cross-section perimeter curve~$\Pss$ 
\begin{equation}
\vsv_0(s,t)={1\over\vert\Pss\vert}\int_{\Pss} \left( \vv^I(\hat\theta,\sv,t) -\wsv(\hat\theta,\sv,t) \right) \,d\hat\theta,
    \qquad \sv\in\OmegaSbar ,
\label{eq:vref}
\end{equation}
where $\vert\Pss\vert$ denotes the length of the perimeter curve $\Pss$.
We note that if the fluid density is not constant, or if the beam separates two fluid regions in two dimensions, then the average in~\eqref{eq:vref} could be replaced by a density-weighted average.  Redefine the reference-curve velocity $\vsv(s,t)$ to equal $\vsv_0(s,t)$.
\item
Set the fluid velocity on the surface of the beam to be
\begin{equation}
\vv(\xsbv(\theta,\sv,t),t)=\vsv(\sv,t) + \wsv(\theta,\sv,t),\qquad \theta\in\Pss,\quad\sv\in\OmegaSbar,
\label{eq:vfref}
\end{equation}
following the kinematic matching conditions in~\eqref{eq:kinematicII} and~\eqref{eq:beamSurfaceVelocity} so that the beam reference-curve velocity provides a consistent definition of the fluid surface velocity.

\end{enumerate}
\end{algorithmVP}

From the perspective of numerical stability, this velocity projection appears to be needed only when the 
beam is quite {\em light}. Otherwise the fluid velocity on the interface could be taken as the velocity of the beam.
This observation is supported by the linear stability analysis in~\cite{fis2014} as well as numerical experiments. However, in practice we
find that the velocity projection is never harmful and we therefore advocate its use. In addition, for very light beams the
interface acts as a free surface with no surface tension (assuming $\Lsv$ is also very small). 
In this case, the fluid interface velocity is smoothed using a few steps of a fourth-order filter
to remove oscillations that can be excited on the interface, see Section~\ref{sec:velocityFilter} for
more details. Finally note that for cases where the beam supports motion only in the normal direction to the beam, only the 
normal component of the fluid velocity is projected, i.e.
\begin{align}
   (\nv^T\vv)^I(\theta,\sv,t) = \gamma (\nv^T\vv_b^{(p)})(\theta,\sv,t) + (1-\gamma) (\nv^T\vsv_b^{(p)})(\theta,\sv,t)  ,  \qquad \theta\in\Pss,\quad \sv\in\OmegaSbar, \label{eq:velocityProjectionNormalComponent}
\end{align}
where $\gamma$ is given by~\eqref{eq:vaverage} as before.
\newcommand{\bc}[1]{\mbox{\bfss#1}}   
\newcommand{\cc}[1]{\mbox{$//$  #1}}  
\newcommand{\ia}{\qquad\qquad}        
\newcommand{\ib}{\ia\quad}     
\newcommand{\ic}{\ib\quad}     
\newcommand{\id}{\ic\quad}     
\newcommand{\ie}{\id\quad}     

\newcommand{\FUNC}[1]{{\color{blue}#1}}
\newcommand{\RETURN}{{\color{blue}Return}}
\newcommand{\IF}{{\color{blue}if}}
\newcommand{\ELSE}{{\color{blue}else}}
\newcommand{\ELSEIF}{{\color{blue}else if}}
\newcommand{\FOR}{{\color{blue}for}}
\newcommand{\COM}[1]{{\color{purple}\em #1}}

\newcommand{\indentI}{\hangindent\parindent\hangafter=0\noindent}

\renewcommand{\OmegaSh}{{\OmegaSbar_h}}

\section{The AMP FSI time-stepping algorithm} \label{sec:timeStepping}


In this section we outline the FSI time-stepping algorithm, including
deforming geometry, for an incompressible fluid governed by the Navier-Stokes equations and a thin structure determined by the
general beam model in~\eqref{eq:beamEquation}. The intent is to provide a concrete
description of the AMP algorithm with sufficient details to indicate the primary extensions required from the original scheme 
described in~\cite{fis2014}. 
The fluid equations are solved in velocity-pressure form using a 
second-order accurate fractional-step algorithm~\cite{fis2014,splitStep2003,ICNS}.
The beam equations are advanced with a second-order accurate predictor-corrector scheme.
In practice we have also used a predictor-corrector version of the Newmark-beta scheme
for the beam equations. 

Given the discrete solution at the current time $t^n$, and one previous time level $t^{n-1}$, 
the goal of the FSI algorithm is to determine the solution at time $t^{n+1}$. 
The algorithm is written for a fixed time-step, $\dt$, so that $t^n=n\dt$, but is easily extended 
to a variable $\dt$. 
A predictor-corrector scheme that implements the FSI algorithm is as follows:


\newcommand{\SumPerim}{\sum_{\displaystyle\kv\in\Pss_{h,\jv}}}

\medskip
\noindent{\bf Begin predictor.}

\medskip
\indentI {\bf Stage I - structure}: Predicted values for the displacement and velocity of the shell are obtained using a second-order accurate {\em leap-frog scheme}
given by
\begin{alignat*}{3}
  \frac{\usv_\jv^{(p)} - \usv_\jv^{n-1}}{2\dt} &= \vsv_\jv^{n}, \qquad&& \jv\in\OmegaSh ,\\
   \rhos\As_\jv\, \frac{\vsv_\jv^{(p)} - \vsv_\jv^{n-1}}{2\dt} &= \Lsv_{h}\bigl(\usv_\jv^n,\vsv_\jv^{n}\bigr) + \fsv_\jv^{n},
               \qquad&& \jv\in\OmegaSh ,
\end{alignat*}
where the external force at $t=t^n$ is given in terms of the fluid stress on the surface of the beam
by a discrete approximation to~\eqref{eq:dynamic},
\begin{align*}
   \fsv_\jv^{n} = \SumPerim (\sigmav\nv)_\kv^n \, \Delta\theta .
\end{align*}
Here $\jv\in\OmegaSh$ denotes the indices of the discrete points along the beam reference curve and $\iv\in\GammaIh$ denotes
the corresponding indices of the discrete fluid points on the beam surface.  There is a cross section of the beam identified for each point $\jv$ along the reference curve, and the points $\kv\in\Pss_{h,\jv}$ denote the discrete representation of the curve on the perimeter of the cross section.

\medskip
\indentI {\bf Stage II - grid generation}: Predicted values for the deformed interface positions are obtained based on the
predicted reference-curve deformation using discrete versions of the formulas in~\eqref{eq:xcenterline} and~\eqref{eq:surfaces}. The predicted fluid grid $\Gs^{(p)}$ is then defined in terms of the interface positions using
an appropriate mapping.  In the context of overlapping grids (see Section~\ref{sec:dcg}) an interface-fitted fluid
grid is generated with a hyperbolic marching algorithm and then a new overlapping grid is generated. 


\newcommand{\pressureSum}{\SumPerim p_{\kv}\sp{(p)}\nv_{\kv}\,\Delta\theta}
\medskip
\indentI {\bf Stage III - fluid velocity}: Predicted values for the fluid velocity are obtained using a second-order accurate
 Adams-Bashforth scheme
\begin{alignat*}{3}
  & \rho \frac{\vv_\iv^{(p)} - \vv_\iv^{n}}{\dt} = \frac{3}{2} \Fv_\iv^n - \frac{1}{2} \Fv_\iv^{n-1} , \qquad&& \iv\in\OmegaFh,
\end{alignat*}
where
\begin{alignat*}{3}
  &  \Fv_\iv^n \equiv  - \big((\uv_\iv^n-\wv_\iv^n)\cdot\grad\big) \uv_\iv^n -\grad_h p_\iv^n + \mu\Delta_h\vv_\iv^n . 
\end{alignat*}
On moving grids the equations are solved in a moving coordinate system (see Section~\ref{sec:dcg}) and $\wv_\iv^n$ is the grid velocity. 
From~\eqref{eq:AMPtangentialVelocityBCI}, the boundary conditions for $\vv_\iv^{(p)}$ on $\GammaIh$ are
\begin{equation}
\left.
\begin{array}{c}
\displaystyle{
\tanv_{m,\iv}^T\Biggl[\;\SumPerim\tauv_{\kv}\sp{(p)}\nv_{\kv}\,\Delta\theta+\frac{\mu\rhos\As_\jv}{\rho}\Delta_h \vv_\iv\sp{(p)}\Biggr]=\tanv_{m,\iv}^T\Biggl[\frac{\rhos\As_\jv}{\rho}\grad_h p_\iv^{(p)}+\bar H_{\jv}^{(p)} + \pressureSum \Biggr] ,
}\medskip\\
\displaystyle{
\grad_h\cdot\vv_\iv^{(p)} =0,
}
\end{array}
\right\}  
\label{eq:velocityBCs}
\end{equation}
%
where $\iv\in\GammaIh$, $\jv\in\OmegaSh$, $m=1,2$, and
\newcommand{\wsvDot}{\dot{\wsv}}
\begin{equation}
\tauv_\iv\sp{(p)} = \mu \left[ \grad_h\vv_\iv\sp{(p)} + \left(\grad_h\vv_\iv\sp{(p)}\right)^T \right],\qquad
\bar H_{\jv}\sp{(p)}=\Lsv_{h}\bigl(\usv_\jv^{(p)},\vsv_\jv^{(p)}\bigr)+ 
      \rhos\As_\jv \wsvDot_\jv^{(p)} .
\label{eq:defH}
\end{equation}
The term $\wsvDot_\jv^{(p)}$ in~\eqref{eq:defH} is a predicted value for $\partial\wsv/\partial t$,
which for an Euler-Bernoulli beam can be computed from the predicted beam position,
velocity and acceleration.
At this stage of the time-stepping algorithm, a predicted value for the pressure, $p_\iv\sp{(p)}$, on the right-hand side of~\eqref{eq:velocityBCs} is not yet available.  However, a sufficiently accurate value for the purposes of the boundary conditions can be obtained using a third-order extrapolation in time, $p_\iv^{(p)} = 3 p_\iv^n - 3p_\iv^{n-1} + p_\iv^{n-2}$.  This extrapolated pressure is replaced in the next stage of the algorithm by solving a discrete Poisson problem.  Appropriate \new{boundary} conditions for $\vv_\iv^{(p)}$ should also be applied on other boundaries.

\medskip
\indentI{\bf Stage IV - fluid pressure}: Predicted values for the pressure are determined by solving the pressure equation
\begin{align}
   \Delta_h p_\iv^{(p)} & = - \rho\, \grad \vv_\iv^{(p)}:\left(\grad\vv_\iv^{(p)}\right)\sp{T}  , \qquad \iv\in\OmegaFh,  \label{eq:pDiscrete}
\end{align}
subject to the boundary conditions, 
%
\[
\nv_\iv\sp{T}\Biggl[\;\pressureSum + \frac{\rhos\As_\jv}{\rho}\grad_h p_\iv^{(p)}\Biggr]
=\nv_\iv\sp{T}\Biggl[-\bar H_{\jv}\sp{(p)}
          +   \frac{\mu\rhos\As_\jv}{\rho}\Delta_h \vv_\iv\sp{(p)} + \SumPerim \tauv_{\kv}\sp{(p)}\nv_{\kv}\,\Delta\theta\Biggr],\qquad\iv\in\GammaIh ,
\]
with $\jv\in\OmegaSh$ and appropriate conditions for $p_\iv^{(p)}$ applied on the other boundaries.  
Note that in practice it is useful to add a divergence damping term to the right-hand side of the the pressure 
equation~\eqref{eq:pDiscrete} following~\cite{ICNS,splitStep2003}.

\medskip
\noindent{\bf End predictor.}

\vskip\baselineskip
\noindent{\bf Begin corrector.}

\medskip
\indentI {\bf Stage V - structure}:
Corrected values for the displacement and velocity of the shell are obtained using a second-order accurate Adams-Moulton scheme, 
\begin{alignat*}{3}
\frac{\usv_\jv^{n+1} - \usv_\jv^{n}}{\dt} &= \vsv_\jv^{n+\half}, \qquad&& \jv\in\OmegaSh  ,\\
\rhos\As_{\jv} \frac{\vsv_\jv^{n+1} - \vsv_\jv^{n}}{\dt} &= \Lsv_{h}\left(\usv_\jv^{n+\half},\vsv_\jv^{n+\half}\right) + \fsv_\jv^{n+\half} , \qquad&&\jv\in\OmegaSh , 
\end{alignat*}
where
\[
\usv_\jv^{n+\half} \equiv {1\over2}\left(\usv_\jv^{(p)} + \usv_\jv^{n}\right) ,\qquad \vsv_\jv^{n+\half} \equiv {1\over2}\left(\vsv_\jv^{(p)} + \vsv_\jv^{n}\right) ,\qquad 
    \fsv_\jv^{n+\half} \equiv {1\over2}\left(\fsv_\jv^{(p)} + \fsv_\jv^{n}\right).  
\]

\medskip
\indentI {\bf Stage VI - grid generation}: The fluid grid could be corrected and regenerated at this stage of the algorithm. In practice,
     however, this stage is generally skipped to avoid the cost of an extra grid generation step. 
   Note that the predicted grid is already second-order accurate and there
   appears to be no noticeable adverse effect on the stability of the algorithm by skipping this stage.

\medskip
\indentI {\bf Stage VII - fluid velocity}: Corrected values for the fluid velocity are obtained using a second-order accurate Adams-Moulton scheme,
\begin{align*}
   \rho \frac{\vv_\iv^{n+1} - \vv_\iv^{n}}{\dt} &= \frac{1}{2} ( \Fv_\iv^{(p)} + \Fv_\iv^{n} ), \qquad \iv\in\OmegaFh.
\end{align*}
The boundary conditions have the same form as those used in Stage III with $\usv_\jv^{(p)}$, $\vsv_\jv^{(p)}$ and $\wsv_\jv^{(p)}$ replaced by $\usv_\iv^{n+1}$, $\vsv_\iv^{n+1}$ and $\wsv_\jv^{(n)}$.

\medskip
\indentI {\bf Stage VIII - fluid pressure}: Corrected values for the pressure are determined by solving the discrete equations in Stage IV with the predicted values replaced by values at $t^{n+1}$.

\medskip
\indentI {\bf Stage IX - project interface velocity}:  Define an interface velocity as a weighted average
of the fluid and solid interface velocities, using
\begin{align*}
   \vv_\iv^I = \gamma \vv_\iv^{n+1} + (1-\gamma) \left( \vsv_\jv^{n+1} +\wsv_\jv^{n+1} \right) ,
            \qquad \iv\in\GammaIh, \quad\jv\in\OmegaSh,
\end{align*}
where $\gamma$ is defined in~\eqref{eq:vaverage}. 
Redefine the beam reference-curve velocity as an average of the surface values (adjusted for $\wsv_\kv^{n+1}$),
\begin{align*}
   \vsv_\jv^{n+1} = \frac{1}{|\Pss_{h,\jv}|} \SumPerim \left( \vv_\kv^I - \wsv_\kv^{n+1}\right)\,\Delta\theta, 
         \qquad\jv\in\OmegaSh,
\end{align*}
and then redefine the fluid surface velocity to be consistent with the new beam velocity, 
\begin{align*}
    \vv_\iv^{n+1} = \vsv_\jv^{n+1} + \wsv_\jv^{n+1} ,\qquad \iv\in\GammaIh,  \quad\jv\in\OmegaSh.
\end{align*}

\medskip
\noindent{\bf End corrector.}

\medskip

We emphasize that the AMP algorithm is stable with no corrector step, although if the predictor step is used alone, then Stage IX should be performed
following the predictor to project the interface velocity. 
We typically use the corrector step, since for the fluid in isolation the 
scheme has a larger stability region than the predictor step alone, and
the stability region includes the imaginary axis so that the scheme can be used for inviscid problems ($\mu=0$).

\newcommand{\Kt}{\Ks_1}
\newcommand{\Kxxt}{\Ts_1}
\newcommand{\stressJump}{\big[\, \nv^T\sigmav\nv \,\big]_\Gamma}

\section{An Euler-Bernoulli beam model} \label{sec:beamModels}

For the purposes of this paper, we consider a two-dimensional beam governed by a generalized Euler-Bernoulli (EB) beam model.  The model considers the behavior of the displacement, $\ww(s,t)$, of the reference curve of the beam in the direction normal to its initially flat position, where $s$ measures distance along the reference curve.  The model assumes small deviations from its initial position (i.e.~small slopes) so that displacement in the direction tangent to the reference curve is neglected.  The generalized EB model includes damping terms and takes the form
\begin{align}
  \rhos \hs \frac{\partial^2\ww}{\partial t^2} = 
                    & -\Ks_0\ww +
                   \frac{\partial}{\partial s}\left( \Ts \frac{\partial\ww}{\partial s}\right)
                   - \frac{\partial^2}{\partial s^2}\left( \Es \Is \frac{\partial^2\ww}{\partial s^2}\right) 
                  - \Kt \frac{\partial\ww}{\partial t}
                   + \Kxxt \frac{\partial^2}{\partial s^2}\left(\frac{\partial\ww}{\partial t}\right)
                   + f(s,t), 
                  \qquad s\in\OmegaSbar,  \label{eq:BeamModel}
\end{align}
where $\rhos$ is the density of the beam, $\hs=\hs(s)$ is its thickness, $\Ks_0$ is the linear stiffness coefficient, $\Ts$ is the tension coefficient, $\Es$ is Young's modulus, and $\Is$ is the area-moment of inertial of the beam.  The term with coefficient $\Kt$ is a linear damping term, while the term with 
coefficient $\Kxxt$ damps higher spatial frequencies in $\ww(s,t)$ faster than lower frequencies and thus
serves to smooth the beam.  We assume that all of the coefficients in~\eqref{eq:BeamModel} are constants, except for the thickness of the beam which may vary with $s$ such as near the rounded end of a cantilevered beam (see Figure~\ref{fig:beamGrids}, for example).

%
%
%

{
\newcommand{\yH}{3}
\newcommand{\xL}{8}
\newcommand{\yUpper}{3}
\newcommand{\yLower}{0}
\def\xa{0}
\def\xb{9}
\def\a{.2}%
\def\b{35}
\def\hsh{1.}
\begin{figure}[hbt]
\newcommand{\textFont}{\normalss}
\begin{center}
\begin{tikzpicture}[scale=1]
\useasboundingbox (0,.1) rectangle (10,5.75);  
\begin{scope}[xshift=0cm,yshift=2.5cm]
\fill[fill=red!20,draw=red,line width=2pt] 
       plot[samples=100, domain=0.:\xb] (\x, {\hsh+\a*sin(\b*\x)}) --
       plot[samples=100, domain=\xb:0] (\x, {-\hsh+\a*sin(\b*\x)}) -- cycle;
\draw[-,red,line width=2pt] plot[samples=100, domain=0.:\xb] (\x, {\a*sin(\b*\x)}); 
\def\xA{\xb} \pgfmathparse{\hsh+\a*sin(\b*\xA)} \let\yA\pgfmathresult
 \draw[black] (\xA,\yA) node[anchor=west] {\textFont $\xsv_+(\ss,t)$};
\def\xA{\xb} \pgfmathparse{\a*sin(\b*\xA)} \let\yA\pgfmathresult
 \draw[black] (\xA,\yA) node[anchor=west] {\textFont $\xbv(\ss,t)$};
\def\xA{\xb} \pgfmathparse{-\hsh+\a*sin(\b*\xA)} \let\yA\pgfmathresult
 \draw[black] (\xA,\yA) node[anchor=west] {\textFont $\xsv_-(\ss,t)$};
%
\def\xr{4.} \pgfmathparse{\a*sin(\b*\xr)} \let\yr\pgfmathresult
\fill[red] (\xr,\yr) circle (3 pt);
\def\xp{4.08} \pgfmathparse{\hsh+\a*sin(\b*\xp)} \let\yp\pgfmathresult
\fill[red] (\xp,\yp) circle (3 pt);
\def\xm{3.92} \pgfmathparse{-\hsh+\a*sin(\b*\xm)} \let\ym\pgfmathresult
\fill[red] (\xm,\ym) circle (3 pt);
  \draw[<->,line width=1.pt,black] (\xm,\ym) -- (\xr,\yr) node[anchor=south west,yshift=0pt] {$\hs$} -- (\xp,\yp); 
\def\xA{5.95} \pgfmathparse{\a*sin(\b*\xA)} \let\yA\pgfmathresult
\def\xB{6.05} \pgfmathparse{.8*\hsh+\a*sin(\b*\xA)} \let\yB\pgfmathresult
 \draw[->,line width=1.pt,black] (\xA,\yA) -- (\xB,\yB) node[anchor=east,yshift=-10pt] {$\nsv$};
%
\def\yTop{2.5}
\def\yBot{-2.5}
\fill[fill=blue!20,draw=blue,line width=2pt,yshift=4pt] 
       (\xa,\yTop) --
       plot[samples=100, domain=0.:\xb] (\x, {\hsh+\a*sin(\b*\x)}) -- (\xb,\yTop)  -- cycle;
\fill[fill=blue!20,draw=blue,line width=2pt,yshift=-4pt] 
       (\xa,\yBot) --
       plot[samples=100, domain=0.:\xb] (\x, {-\hsh+\a*sin(\b*\x)}) -- (\xb,\yBot)  -- cycle;
%
\def\xA{1.97} \pgfmathparse{\hsh+\a*sin(\b*\xA)} \let\yA\pgfmathresult
\def\xB{2.03} \pgfmathparse{.2+\a*sin(\b*\xA)} \let\yB\pgfmathresult
\draw[->,line width=1.pt,black,yshift=3pt] (\xA,\yA) -- (\xB,\yB) node[anchor=east,yshift=+10pt] {$\nv_+$}; 
\def\xA{2.02} \pgfmathparse{-\hsh+\a*sin(\b*\xA)} \let\yA\pgfmathresult
\def\xB{1.98} \pgfmathparse{-.2+\a*sin(\b*\xA)} \let\yB\pgfmathresult
\draw[->,line width=1.pt,black,yshift=-3pt] (\xA,\yA) -- (\xB,\yB) node[anchor=east,xshift=2pt,yshift=-10pt] {$\nv_-$}; 
\end{scope}
%
%
\end{tikzpicture}
\end{center}
\caption{Geometry of an Euler-Bernoulli beam. The beam reference-line is $\xvb=\xvb(s,t)$ and the beam normal is $\nsv$.
 The beam has thickness $\hs$ and meets the fluid on the curves $\xsv_-$ and $\xsv_+$. Outward normals to the fluid domain are $\nv_-$ and $\nv_+$. 
 }
\label{fig:beamCartoon}
\end{figure}
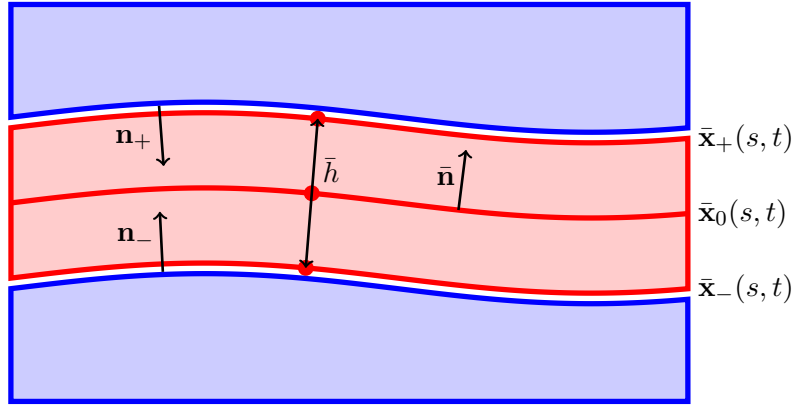
}

%

The geometry of an FSI problem involving a section of an EB beam is illustrated in Figure~\ref{fig:beamCartoon}.  For this configuration, the displacement of the reference curve is given by $\usv(s,t)=(0,\eta(\ss,t))$ so that its position is $\xbv(\ss,t)=(\ss,\eta(\ss,t))$.  The unit normal to the reference curve, $\nsv(\ss,t)$, can be computed in terms of the slope $\partial\ww/\partial s$, and this determines the position of the upper and lower surfaces of the beam as
\begin{align}
  \xsv_{\pm}(\ss,t)=\xbv(\ss,t)+\xvz_\pm(\sv,t),\qquad \xvz_\pm(\sv,t)=\pm \frac{\hs}{2} \nsv(\ss,t), \qquad s\in\OmegaSbar,\label{eq:EBsurface}
\end{align}
which is a specific case of the formula given previously in~\eqref{eq:surfaces2d} with the $b$ subscript dropped for notational convenience.  From~\eqref{eq:EBsurface}, the velocity of the upper and lower surfaces of the beam are given by
\begin{align}
   \vsv_\pm(s,t) =  \frac{\partial}{\partial t}\xsv_{\pm}(\ss,t) = \vsv(\ss,t) + \wsv_\pm(s,t), \qquad \wsv_\pm(s,t) = \pm \frac{\hs}{2} \partialt \nsv(\ss,t).
      \label{eq:beamSurfaceVelocityEB}
\end{align}
Finally, a modeling choice is required to determine the force $f(s,t)$ on the reference curve of the beam, and we use
\begin{align}
    f(s,t) &= -\Big[ (\nv^T\sigmav\nv)_+ + (\nv^T\sigmav\nv)_- \Big]. \label{eq:EB_BeamForce}
\end{align}
where $\nv$ is the unit fluid normal on the beam surface so that $(\nv^T\sigmav\nv)_\pm$ is the component of the fluid traction in the
direction normal to the upper and lower surfaces of the beam.  This modeling approximation is valid for small beam deflections, and
it is consistent with the use of an EB beam model and (a reduced form of) the AMP pressure condition in~\eqref{eq:AMPpressureBCI}
where the right-hand side of~\eqref{eq:EB_BeamForce} becomes the left-hand side of~\eqref{eq:AMPpressureBCI}.

For a cantilevered beam, the free end of the beam is usually rounded as shown in Figure~\ref{fig:beamGrids}. The flow can
be complex near a sharp corner and by rounding the corner a high resolution grid can be
generated that resolves the flow in this region. In addition, the force on the beam surface 
is computed using the beam-surface normal $\nv$ and the approximation $\nv^T\sigmav\nv$. Near the rounded end of a beam 
this is a poor approximation since the fluid normal may be orthogonal to the beam reference curve. Consequently, 
in practice, the fluid forces on the beam tip are ignored in computing force on the beam (in particular, any points on 
the surface of the beam that extend beyond the length of the beam reference line are ignored). 
This should be considered a modeling assumption. Similarly, contributions from the end points are also ignored 
when projecting the interface velocity onto the beam.

\newcommand{\pJumpN}{\big[\, \ph^{n} \,\big]_\Gamma}
\newcommand{\pJumpNp}{\big[\, \ph^{n+1} \,\big]_\Gamma}
\newcommand{\pJump}{\big[\, p \,\big]_\Gamma}
\newcommand{\phJump}{\big[\, \ph \,\big]_\Gamma}
\section{Analysis of a model problem for a two-dimensional beam} \label{sec:analysis}

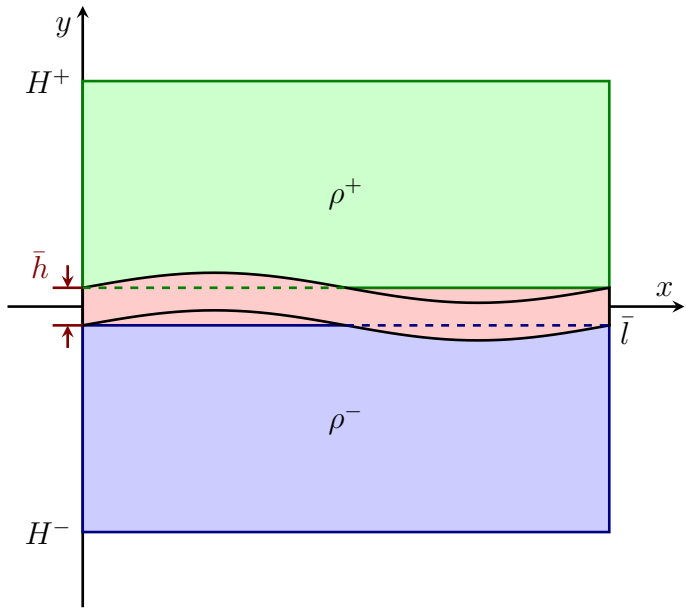
\begin{figure}[hbt]
\newcommand{\textFont}{\normalss}
\begin{center}
\begin{tikzpicture}[line width= 1pt,>=stealth]
\useasboundingbox (0,1.5) rectangle (9,9);

\draw[->] (1,1) -- (1,9);
\filldraw[fill=red!20!white, draw=red!50!black] (1,4.75) rectangle (8,5.25) coordinate;
\filldraw[fill=blue!20!white, draw=blue!50!black] (1,2) rectangle (8,4.75) coordinate;
\filldraw[fill=green!20!white, draw=green!50!black] (1,5.25) rectangle (8,8);

\draw[->] (0,5) -- (9,5);
\filldraw[fill=red!20!white] plot[domain=1:8,smooth]({\x},{0.2*sin(deg((\x-1)*2*pi/7))+5+0.25})-- plot[domain=8:1,smooth]({\x},{0.2*sin(deg((\x-1)*2*pi/7))+5-0.25}) -- cycle;

\draw[green!50!black,dashed] (1,5.25) -- (4.5,5.25);
\draw[blue!50!black,dashed] (4.5,4.75)--(8,4.75);

\draw[red!50!black] (1,4.75) -- +(-0.4,0);
\draw[red!50!black] (1,5.25) -- +(-0.4,0);
\draw[<-,red!50!black] (0.8, 4.75) -- +(0,-0.3);
\draw[<-,red!50!black] (0.8, 5.25) -- +(0,0.3);
\draw[red!50!black]  (0.7,5.55) node[anchor=east]{\large$\bar{h}$};
\draw (4.5,6.5) node {\large$\rho^{+}$};
\draw (4.5,3.5) node {\large$\rho^{-}$};
\draw (1,8) node[anchor=east] {\large$H^{+}$};
\draw (1,2) node[anchor=east] {\large$H^{-}$};
\draw (8,5) node[anchor=north west] {\large$\bar{l}$};
\draw (1,9) node[anchor=north east] {\large$y$};
\draw (9,5) node[anchor=south east] {\large$x$};
\end{tikzpicture}
\end{center}
\caption{The geometry of the linearized model  for  a beam with fluid on two sides.}
\label{fig:linearizedTwoSidedBeam}
\end{figure}

In this section, the stability of the AMP algorithm is investigated for a linearized model problem consisting of an EB beam coupled to inviscid, incompressible fluids with different constant fluid densities in the regions above and below the beam as shown in
Figure~\ref{fig:linearizedTwoSidedBeam}.  For comparative purposes, we also consider the stability of a traditional partitioned (TP) scheme for the same model problem.  We begin by obtaining the exact solution of the model problem to explicitly identify the contribution of the
added mass of the fluid for two different choices of the boundary conditions on the fluid. The stability of the discrete AMP and TP schemes is then investigated using a mode analysis that
extends the theory for the case of a beam with fluid on one side only presented in~\cite{fis2014}.
 For both the AMP and TP schemes, the fluid
equations are advanced using a backward-Euler time integrator as opposed to the predictor-corrector integrator
described in Section~\ref{sec:AMPalgorithm}.  While this simplifies the analysis somewhat, the essential results are unaffected.  It is found,
in agreement with~\cite{fis2014}, that the AMP scheme is stable for arbitrarily light structures 
when the usual time-step restriction is satisfied.  The TP scheme without sub-time step iteration, on the other hand, becomes unconditionally unstable for light structures.  Note that the TP scheme has been analyzed previously by Causin, Gerbeau
and Nobile~\cite{CausinGerbeauNobile2005}, among others, and the results presented here extend those works.

In the following analysis the two fluids
are coupled to an
EB beam with constant thickness $\hs$ and without damping terms. The equations describing the evolution of the fluid-structure system are then linearized 
about flat undeflected beam surfaces at $y=-\hs/2$ and $\hs/2$ as shown in Figure~\ref{fig:linearizedTwoSidedBeam}.
Let $\rho^+$, $\vv^+(x,y,t)$ and $p^+(x,y,t)$, denote the fluid density, velocity and pressure, respectively, in the upper domain, $\OmegaF_{+}=[0,\ls]\times[\hs/2,\Hf^{+}]$, 
and $\rho^-$, $\vv^-(x,y,t)$ and $p^-(x,y,t)$, the corresponding quantities in the lower domain $\OmegaF_{-}=[0,\ls]\times[\Hf^{-},-\hs/2]$, where $\ls$ is the length of the beam.
The governing equations and boundary conditions for this model are 
\begin{equation}
\begin{array}{rl}
\text{Fluid:} &
\displaystyle{
 \left\{ 
   \begin{alignedat}{3}
  &  \rho^{\pm}\frac{\partial\vv^{\pm}}{\partial t} 
                 + \grad p^{\pm} =0, \qquad&& \xv\in\OmegaF_{\pm} ,  \\
  & \grad\cdot\vv^{\pm} =0,  \qquad&& \xv\in\OmegaF_{\pm} ,   \\
  & v_2^\pm=0\text{ or } p^{\pm}=0  , \qquad&& x\in(0,\ls), \quad y=\Hf^{\pm}
   \end{alignedat}  \right. 
} 
\medskip\\
\text{Beam:} &
\displaystyle{
  \left\{
   \begin{alignedat}{3}
  &  \rhos\hs\frac{\partial^2\eta}{\partial t^2}  = -\Ks_0 \eta + \Ts\frac{\partial^2\eta}{\partial x^2} - \Es\Is\frac{\partial^4\eta}{\partial x^4} -\pJump, \qquad&&  x\in(0,\ls),
   \end{alignedat}
   \right.
}
\medskip\\
\text{Interface:} &
\displaystyle{
  \left\{
   \begin{alignedat}{3}
  & v_2^{\pm}=\frac{\partial\eta}{\partial t}, \qquad && x\in(0,\ls), \quad y=\pm\hsHalf, 
   \end{alignedat}  \right.
}
\end{array}
  \label{eq:MP1}
\end{equation}
where 
\begin{align*}
 \pJump \equiv p^{+}(x,\hs/2,t) -p^{-}(x,-\hs/2,t),
\end{align*}
denotes the jump of the pressure across the beam, and where we take
periodic boundary conditions in the $x$-direction.
Initial conditions for $\vv^\pm$, $\eta$ and $\partial\eta/\partial t$ are needed to complete the specification of the problem, but these are not required in the subsequent analysis. 
%

\subsection{Exact solution of the continuous model problem} \label{subsec:slipWallBCs}

It is instructive to first examine the behavior of the exact solution of the FSI problem in~\eqref{eq:MP1}.  Since this model problem is a linear constant-coefficient system of equations with periodic boundary conditions in the $x$-direction, the solution can be expressed in terms of a Fourier series in~$x$, e.g.
\[
\vv^{\pm}(x,y,t)=\sum_{k=-\infty}\sp\infty \vvh^{\pm}(k,y,t) \, e^{i\kx x},\qquad \kx = 2\pi k/\ls,
\]
and similar expressions for $p\sp\pm$ and $\eta\sp\pm$.  Here, $\kx$ identifies the wave number of the Fourier mode.  Transforming the equations for the fluid to Fourier space gives
\begin{equation}
\left. \begin{array}{r}
\displaystyle{
  \rho^{\pm} \frac{\partial \vh_1^{\pm}}{\partial t} + i\kx \ph^{\pm} =0
} \smallskip\\
\displaystyle{
  \rho^{\pm} \frac{\partial \vh_2^{\pm}}{\partial t} + \frac{\partial \ph^{\pm}}{\partial y} =0
} \smallskip\\
\displaystyle{
  i\kx\vh_1^{\pm}+{\partial\vh_2^{\pm}\over\partial y} =0
}
\end{array}
\right\},\qquad
\begin{array}{l}
  \hbox{$y \in (\Hf^{-},-\hs/2)$ for the lower fluid},\medskip\\
  \hbox{$y \in (\hs/2,\Hf^{+})$ for the upper fluid}.\\
\end{array}
\label{eq:modelFluidEquationFT}
\end{equation}
The transformed equation for the beam is
\begin{equation}
  \rhos\hs \frac{\partial ^2\etah}{\partial t^2} = - \Lt \etah -\phJump, 
\label{eq:modelShellEquationFT}
\end{equation}
where $\Lt=\Ks_0+\Ts\kx^2+\Es\Is\kx^4$, 
and for notational brevity we have suppressed the functional dependence on the wave number.
The transformed equations for the fluid in~\eqref{eq:modelFluidEquationFT} can be manipulated to derive equations for the pressures of the form
\begin{equation}
{\partial\sp2\ph\sp\pm\over\partial y\sp2}-\kx\sp2\ph\sp\pm=0.
\label{eq:modelPressureEquationFT}
\end{equation}
It is assumed that $\kx\ne0$ in the present analysis while the special case when $\kx=0$ is considered later.  For the case of the rigid-wall boundary conditions, $\vh_2^{\pm}=0$ at $y=\Hf^{\pm}$, we use the $y$-momentum equations in~\eqref{eq:modelFluidEquationFT} to derive the Neumann conditions for the pressures,
\begin{equation}
  \frac{\partial \ph^{\pm}}{\partial y}(\Hf^{\pm},t)=0.
  \label{eq:modelBCFT}
\end{equation}
General solutions of the pressure equations in~\eqref{eq:modelPressureEquationFT} satisfying the homogeneous Neumann boundary conditions in~\eqref{eq:modelBCFT}~are
\begin{equation}
  \ph^{\pm} = a^{\pm}\cosh\left(\kx(y-H^{\pm})\right), \label{eq:pSolution}
\end{equation}
where $a^{\pm}=a^{\pm}(t)$ are a time-dependent coefficients.  Similarly, the kinematic conditions at the interface,
\[
\vh_2^{\pm}(\pm\hs/2,t)=\frac{\partial\etah}{\partial t},
\]
together with the $y$-momentum equations in~\eqref{eq:modelFluidEquationFT} give the interface conditions
\begin{equation}
  \frac{\partial \ph^{\pm}}{\partial y}(\pm\hs/2,t)=  -\rho^{\pm}\frac{\partial^2\etah}{\partial t^2}.
  \label{eq:modelInterfaceFTpy}
\end{equation}
The general solutions in~\eqref{eq:pSolution} and the interface conditions in~\eqref{eq:modelInterfaceFTpy} imply
\begin{equation}
  a^{\pm} = \pm\frac{\rho^{\pm}}{\kx}\left(\frac{\partial^2\etah}{\partial t^2}\right){\rm csch}\left(\kx\mathcal{D}\sp\pm\right),
\label{eq:acoeffs}
\end{equation}
where $\mathcal{D}^{\pm}=\pm\Hf^{\pm}-\hs/2$ measure the vertical extent of the upper and lower fluid domains.  Using~\eqref{eq:pSolution} and~\eqref{eq:acoeffs},
the applied force appearing in the beam equation~\eqref{eq:modelShellEquationFT} takes the form
\begin{equation}
   \phJump = \hat{M}_a^{-}\frac{\partial^2\etah}{\partial t^2} + \hat{M}_{a}^{+}\frac{\partial^2\etah}{\partial t^2},
 \label{eq:appliedForceFT}
\end{equation}
where
\begin{equation}
  \hat{M}_a^\pm =\rho^\pm\mathcal{D}^\pm \, \left[{\coth(\kx\mathcal{D}^\pm)\over\kx\mathcal{D}^\pm}\right] ,\qquad\hbox{(rigid-wall BCs),}
  \label{eq:addedMass}
\end{equation}
are the added masses from the upper and lower fluids for the case of rigid-wall boundary conditions at $y=\Hf^{\pm}$.  
Note that $\rho^\pm\mathcal{D}^\pm$ represents the total mass of fluid in the upper and lower domains.
For the alternate case of the pressure boundary conditions, $\ph^{\pm}=0$ at $y=\Hf^{\pm}$, a similar analysis leads to the added masses given by
\begin{equation}
  \hat{M}_a^\pm =\rho^\pm\mathcal{D}^\pm \, \left[{\tanh(\kx\mathcal{D}^\pm)\over\kx\mathcal{D}^\pm}\right] ,\qquad\hbox{(pressure BCs).}
  \label{eq:addedMassPressure}
\end{equation}
Finally, using the applied force (\ref{eq:appliedForceFT}) in the equation of motion
for the beam (\ref{eq:modelShellEquationFT}) gives
\begin{equation}
 \Big( \rhos\hs+\hat{M}_a^{-}+\hat{M}_a^{+}\Big) \frac{\partial ^2\etah}{\partial t^2} = - \Lt \etah,
 \label{eq:beamMotionFT}
\end{equation}
which explicitly identifies the contribution of the added mass of the fluid in the beam equation. A solution
to~\eqref{eq:beamMotionFT} satisfying initial conditions for $\etah$ and $\partial\etah/\partial t$ 
can be found easily to complete the solution 
for $\etah$. Given $\etah$, the pressure is then fully determined from~\eqref{eq:pSolution} and~\eqref{eq:acoeffs}, which in
turn determines $\vh_1^{\pm}$ and $\vh_2^{\pm}$ from~\eqref{eq:modelFluidEquationFT}
together with the initial conditions and boundary conditions.
This completes the solution of the model problem in Fourier space.

We now return to the special case when $\kx=0$.
For rigid-wall boundary conditions, the full solution for $\kx=0$ is found to be
\begin{align}
  \vh_1^\pm = \vh_1^\pm(y,0), \qquad \vh_2^\pm =0, \qquad \etah=\etah(0), \qquad \ph^\pm=\ph^\pm_0(t),  \label{eq:modelSolutionModeZero}
\end{align}
where $\ph^\pm_0(t)$ are spatially uniform fluid pressures whose difference $\Delta\ph_0=\ph^+_0-\ph^-_0$ is a constant
related to the beam displacement by
\[
\Ks_0\etah(0)=-\Delta\ph_0.
\]
Note that the $\kx=0$ mode of the beam is stationary, consistent
with the fact that the volumes of the lower and upper domains must remain constant for incompressible fluids.
The added mass in this case is interpreted to be infinite. This is fully consistent with the observation that
the added mass in~\eqref{eq:addedMass}
tends to infinity as $\kx\mathcal{D}^\pm$ tends to zero.
 For the case of pressure boundary conditions with $\kx=0$, solutions for the fluid pressures become
\[
\ph\sp\pm=-\rho\sp\pm\left({\partial\sp2\etah\over\partial t\sp2}\right)\bigl(y-H\sp\pm\bigr),
\]
so that the added masses in~\eqref{eq:appliedForceFT} are
\[
\hat{M}_a^\pm=\rho\sp\pm\mathcal{D}^\pm.
\]
This is in agreement with the limit $\kx\mathcal{D}^\pm\rightarrow0$ for the formula for $\hat{M}_a^\pm$ in~\eqref{eq:addedMassPressure}.

Returning to the equation for $\etah(t,\kx)$ in~\eqref{eq:beamMotionFT}, we note that added mass associated with the fluids on
either side of the beam tends to reduce the acceleration of the beam.  Since $\hat{M}_a^{\pm}\rightarrow0$ as
$\kx\mathcal{D}^\pm\rightarrow\infty$ for either~\eqref{eq:addedMass} or~\eqref{eq:addedMassPressure}, the effect of the added mass
is small for high-frequency components of the solution.  For the case of rigid-wall boundary conditions,
$\hat{M}_a^{\pm}\rightarrow\infty$ as $\vert\kx\mathcal{D}\vert$ becomes small so that the effect of the added mass is large in this
limit.  This suggests that the rigid-wall case could be difficult numerically if the algorithm is sensitive to the effects of
added mass.  For the case of the pressure boundary condition, $\hat{M}_a^{\pm}$ becomes constant when $\kx\mathcal{D}=0$ and so this
case would likely be less difficult numerically than the case of rigid walls.

\subsection{Stability analysis of the AMP and TP schemes}
The stability of a partitioned scheme for an FSI problem depends crucially on the treatment of the interface and the resulting boundary conditions applied to the fluid and structure on either side.  AMP interface conditions have been developed and used in the AMP FSI time-stepping algorithm, and the purpose of this section is to analyze the stability of a form of this algorithm applied to the model problem in~\eqref{eq:MP1}.  A starting point for this analysis is the system of equations in~\eqref{eq:modelFluidEquationFT} and~\eqref{eq:modelShellEquationFT} already written in Fourier space.  For this system, we consider a discretization in time, but let the dependence on $y$ remain continuous.  This approach reveals the stability of the partitioned algorithm, while keeping the algebraic details as simple as possible.

We begin with an analysis of the AMP algorithm, and then compare the results with a traditional partitioned (TP) in the end.  For both cases, let
\[
\vvh^{n\pm}(y)\approx\vvh^{\pm}(y,t^n), \qquad \ph^{n\pm}(y)\approx\ph^{\pm}(y,t^n), \qquad \etah^{n}\approx\etah(t^n),
\]
where $t^n=n\dt$. Furthermore, let $D_{+t}$ and $D_{-t}$ denote the forward and backward divided different operators in time 
(e.g. $D_{+t}\etah^n=(\etah^{n+1}-\etah^n)/\dt$ and $D_{-t}\etah^n=(\etah^{n}-\etah^{n-1})/\dt$).

The following AMP scheme is an abbreviated version of the full scheme described in Section~\ref{sec:timeStepping}:

\begin{algorithm}[H] {\bf AMP scheme:} 

\noindent Stage I: Advance the beam displacement using a leap-frog scheme,
\begin{equation}
  \rhos\hs D_{+t}D_{-t}\etah^n  = -\Lt \etah^n - \pJumpN.
\label{eq:etaLeapFrogAM}
\end{equation}
Stage II: Advance the fluid velocity and pressure using a backward-Euler scheme applied to the equations in pressure-velocity form,
\begin{equation}
\left. \begin{array}{r}
\displaystyle{
\rho^{\pm} D_{+t} \vh_1^{n\pm}  + i\kx\, \ph^{(n+1)\pm} =0
} \smallskip\\
\displaystyle{
\rho^{\pm} D_{+t} \vh_2^{n\pm} + \frac{\partial \ph}{\partial y}^{(n+1)\pm} =0
} \smallskip\\
\displaystyle{
\frac{\partial^2 \ph}{\partial y^2}^{(n+1)\pm} - \kx^2\, \ph^{(n+1)\pm} =0
}
\end{array} \right\},
\begin{array}{l}
  \qquad y \in (\Hf^{-},-\hs/2) \hbox{ for the lower fluid},\smallskip\\
  \qquad y \in (\hs/2,\Hf^{+}) \hbox{ for the upper fluid},\\
\end{array}
\label{eq:fluidAMP}
\end{equation}
with the AMP Robin interface conditions derived from \eqref{eq:AMPpressureBCI},
\begin{equation}
-\pJumpNp + \frac{\rhos\hs}{\rho^{\pm}} \frac{\partial \ph}{\partial y}^{(n+1)\pm} = \Lt \etah^{n+1},\qquad y=\pm\hsHalf,
\label{eq:AMPinterfaceBC}
\end{equation}
and boundary conditions,
\begin{equation}
\frac{\partial \ph}{\partial y}^{(n+1)\pm}=0, \qquad y=\Hf^{\pm}.
\label{eq:rbbcs}
\end{equation}
\label{alg:AMP}
\end{algorithm}

A key ingredient of the AMP scheme is the 
Robin condition in~\eqref{eq:AMPinterfaceBC} for the fluid pressure. This condition, along with~\eqref{eq:etaLeapFrogAM}, implies
\begin{equation}
\frac{\partial \ph}{\partial y}^{(n+1)\pm} = -\rho^{\pm}  D_{+t}D_{-t}\etah^{n+1},\qquad y=\pm\hsHalf,
\label{eq:balance}
\end{equation}
which is a discrete version of the interface condition in~\eqref{eq:modelInterfaceFTpy}.
We observe that the acceleration term on the right-hand side of~\eqref{eq:balance} is evaluated at the {\em same time level} as the pressure gradient term on the left-hand side, and this plays an important role in the stability of the AMP scheme as is noted below.
After solving the pressure equation for the fluid subject to the boundary condition in~\eqref{eq:rbbcs} and interface condition in~\eqref{eq:balance}, 
the applied force is found to be
\begin{equation}
 -\pJumpNp =-\hat{M}_a^{-}D_{+t}D_{-t}\etah^{n+1}-\hat{M}_{a}^{+}D_{+t}D_{-t}\etah^{n+1},
 \label{eq:appliedForceAMP}
\end{equation}
where the added masses, $\hat{M}_a^\pm$, are given in~\eqref{eq:addedMass}.
Substituting the expression for the applied force~\eqref{eq:appliedForceAMP}, with $n+1$ replaced by $n$, into the leap-frog scheme~\eqref{eq:etaLeapFrogAM} gives the difference equation for $\etah^n$, 
\begin{align}
  \left(\rhos\hs+\hat{M}_a^{-}+\hat{M}_a^{+}\right)D_{+t}D_{-t}\etah^{n}& = -\Lt \etah^n,
  \label{eq:etaAMP}
\end{align}
which is analogous to the evolution equation for $\etah(t)$ in~\eqref{eq:beamMotionFT}.  A solution of the form $\etah^n=A^n\etah^0$ is sought which leads to the quadratic equation
\begin{align}
  A^2-2\left(1-\frac{\dt^2\Lt}{2(\rhos\hs+\hat{M}_a^{-}+\hat{M}_a^{+})}\right)A+1=0,
\label{eq:AMP_A} 
\end{align}
for the amplitude $A$.  Stability requires $\vert\amp\vert \le 1$ and that the roots be simple, if they
have modulus one, which leads to the following theorem:
\begin{theorem}
   The AMP scheme applied to the model problem in~\eqref{eq:MP1} is stable if and only if 
\begin{align}
   \dt < 2 \sqrt{\frac{ \rhos\hs+\hat{M}_a^{-}+\hat{M}_a^{+} }{\Lt}}, \label{eq:AMPstability}
\end{align}
where $\hat{M}_a^{\pm}$ is given by~\eqref{eq:addedMass} for case of the rigid-wall BCs and by~\eqref{eq:addedMassPressure} for the case of pressure BCs.  Moreover, if~\eqref{eq:AMPstability} holds, then the roots of~\eqref{eq:AMP_A} satisfy $\vert\amp\vert=1$ and the AMP scheme for the model problem is non-dissipative.
\end{theorem}

We note that if a stable time step is chosen for the leap-frog scheme in~\eqref{eq:etaLeapFrogAM} with no external forcing due to the fluid, namely $\dt<2 \sqrt{\rhos\hs/\Lt}$, then this value of $\dt$ would be stable for the AMP scheme as well since $\hat{M}_a^{\pm}\ge0$.
We also note that a velocity projection has been omitted in the prototype AMP scheme considered in the stability analysis, and thus the fluid and structure velocities 
at the interface would only match to first-order accuracy.  An additional projection step to match the velocities on the interface exactly could be included, but this is not essential for the present analysis.

Having considered the stability of the AMP scheme, we now turn our attention to a TP scheme defined as follows:

\begin{algorithm}[H] {\bf Traditional partitioned (TP) scheme:}

\noindent Stage I: Advance the beam displacement using~\eqref{eq:etaLeapFrogAM}.

\medskip\noindent
Stage II: Advance the fluid velocity and pressure using~\eqref{eq:fluidAMP} with interface condition,
\begin{equation} 
\frac{\partial \ph}{\partial y}^{(n+1)^{\pm}} = -\rho^{\pm}  D_{+t}D_{-t}\etah^n,  \qquad y=\pm\hsHalf, \label{eq:interfaceTraditional}
\end{equation}
and the boundary condition in~\eqref{eq:rbbcs}.
\label{alg:TP}
\end{algorithm}

The only difference between the AMP scheme and the TP scheme is the coupling at the interface.  For the TP and AMP schemes, the stress from the fluid at time level $t^n$ provides a forcing to the leap-frog scheme in~\eqref{eq:etaLeapFrogAM} to advance the displacement of the beam to time level $t^{n+1}$.  The velocity and pressure in the fluid are then advanced to time level $t^{n+1}$ using the discrete equations in~\eqref{eq:fluidAMP}, the boundary condition in~\eqref{eq:rbbcs} and a coupling condition at the interface.  The AMP scheme uses the coupling condition in~\eqref{eq:AMPinterfaceBC}, whereas the TP scheme uses~\eqref{eq:interfaceTraditional}.  The latter is derived from the time derivative of the kinematic interface condition matching the vertical velocity of the beam to that of the fluid on interface.  We note, however, that a similar interface condition was derived from the leap-frog scheme in~\eqref{eq:etaLeapFrogAM} and the AMP interface coupling condition in~\eqref{eq:AMPinterfaceBC}, which appears in~\eqref{eq:balance}.  The {\em key} difference is that the time level for the vertical acceleration of the beam in~\eqref{eq:interfaceTraditional} lags that of fluid in contrast to~\eqref{eq:balance}, and this is the essential issue leading to the instability of the TP scheme for light structures as the subsequent analysis shows.

Following the steps used in the analysis of the AMP scheme, the pressure equation in \eqref{eq:fluidAMP} is solved with interface conditions in~\eqref{eq:interfaceTraditional} and
boundary conditions in~\eqref{eq:rbbcs} to give the applied force
\begin{equation}
 -\pJumpNp =-\hat{M}_a^{-}D_{+t}D_{-t}\etah^n-\hat{M}_{a}^{+}D_{+t}D_{-t}\etah^n.
 \label{eq:appliedForceTP}
 \end{equation}
Using the beam equation in~\eqref{eq:etaLeapFrogAM} together with~\eqref{eq:appliedForceTP} to eliminate the applied force yields a difference equation for $\etah\sp{n}$ of the form
\begin{align}
  \rhos\hs D_{+t}D_{-t}\etah^{n} +\Big(\hat{M}_a^{-}+\hat{M}_a^{+}\Big)D_{+t}D_{-t}\etah^{n-1}& = -\Lt \etah^n.
\label{eq:etaTraditional} 
\end{align}
Again we observe a time lag in the added-mass terms in~\eqref{eq:etaTraditional} as compared to the previous evolution equation for $\etah^n$ in~\eqref{eq:etaAMP} for the AMP scheme. This time-lag is well known and is the primary source of the added-mass instability for the TP scheme.  Using $\etah^n=A^n\etah^0$ in~\eqref{eq:etaTraditional}  leads to the cubic equation for $A$, 
\begin{align}
  \rhos\hs A(A-1)^2+\dt^2\Lt A^2+(\hat{M}_a^{-}+\hat{M}_a^{+})(A-1)^2=0.
\label{eq:TP_A} 
\end{align}
The TP scheme is weakly stable if all roots of \eqref{eq:TP_A} satisfy $|A|\le 1$, and with no repeated roots of modulus one. Using the theory
of von Neumann polynomials~\cite{Miller1971,Strikwerda89} leads to the following theorem.
\begin{theorem} \label{th:traditional}
  For $\rhos\hs>0$ and $\Lt>0$, the traditional partitioned scheme applied to the model problem in~\eqref{eq:MP1} is weakly stable if and only if
\begin{align}
  \hat{M}_a^{-}+\hat{M}_a^{+} < \rhos\hs \label{eq:TP_addedMassLimit}
\end{align}
  and
\begin{align}
  &  \dt <  2\sqrt{ \frac{\rhos\hs-\hat{M}_a^{-}-\hat{M}_a^{+}}{\Lt}}, \label{eq:traditionalStabLimitII}
\end{align}
where $\hat{M}_a^{\pm}$ is given by~\eqref{eq:addedMass} for case of the rigid-wall BCs and by~\eqref{eq:addedMassPressure} for the case of pressure BCs.
\end{theorem}

The condition in \eqref{eq:TP_addedMassLimit} implies that the TP scheme is unconditionally unstable if 
the contribution from the added mass is too 
large, i.e.~if $\hat{M}_a^{-}+\hat{M}_a^{+} >\rhos \hs$, regardless of the choice for the time-step, $\dt$.
In general, the results of the analysis suggest that TP-type schemes would suffer from stability issues when
$\hat{M}_a^\pm/\rhos\hs$ is large.  We note from~\eqref{eq:addedMass} or~\eqref{eq:addedMassPressure} that $\hat{M}_a^\pm/\rhos\hs$ is large if the mass of the fluid domain is large compared to that of the structure, i.e.~if $\rho\sp\pm \mathcal{D}^{\pm}/\rhos\hs$
is large.  However, the ratio is also large for the case of rigid-wall boundary conditions if $\vert\kx\mathcal{D}\sp\pm\vert$ is small in
view of the bracketed term in~\eqref{eq:addedMass}.  This latter
condition can be re-interpreted as $\mathcal{D}\sp\pm/\ls$ is small, i.e. the fluid domain is long and slender, since $\kx=2\pi k/\ls$.
These observations are in  agreement with the results of Causin et~al.~\cite{CausinGerbeauNobile2005}, who 
studied a simple beam with fluid on one side.  The observations are also in agreement with the numerical results presented here and in~\cite{fis2014}.
If the condition $\hat{M}_a^{-}+\hat{M}_a^{+}<\rhos\hs$ is satisfied, then \eqref{eq:traditionalStabLimitII} indicates how the time-step 
must be reduced as the total added-mass increases.




\subsection{Pressure regularization for rigid-wall boundary conditions}

Note that for both the TP and AMP schemes the limit $\kx=0$ with rigid-wall boundary conditions must be treated 
with some care since the continuous fluid pressure is only determined up to an arbitrary additive constant.
In the AMP scheme $p^+$ and $p^-$ are coupled so that only the difference in pressure is arbitrary, whereas
both pressures are arbitrary for the TP scheme.
In a numerical implementation these singularities need to be removed since, due to round-off or truncation
errors, the singular systems may have no solution. The specific choice of regularization can affect the overall stability
of the approach. For example, one common regularization would follow the approach in~\cite{splitStep2003},
which for $\kx=0$ would solve the augmented pressure equations 
\[
   \frac{\partial^2 \ph^{(n+1)\pm}}{\partial y^2} +\alpha^{(n+1)\pm} =0,
\]
with the constraints
\[
   \int_{\Hf^{-}}^{-\hsHalf} \ph^{(n+1)-}dy=\int^{\Hf^{+}}_{\hsHalf} \ph^{(n+1)+}dy  = 0.
\]
Here $\alpha^{(n+1)\pm}$ are two additional unknowns which serve to regularize the solution and ensure that the augmented 
system has a unique solution. The effect of regularization can be significant for the $\kx=0$ mode.  Using this
regularization the traditional partitioned scheme is found to be stable for the $\kx=0$ mode under the mollified time-step restriction
\[
  \dt <  2\sqrt{ \frac{\rhos\hs-\frac{1}{3} \rho^{-}\mathcal{D}^{-}-\frac{1}{3} \rho^{+}\mathcal{D}^{+}}{\Lt}}.
\]
This is significant difference from \eqref{eq:traditionalStabLimitII} which required $\dt\to 0$ when
$\kx\to 0$. This less severe time-step restriction is a result of the pressure regularization
and applies only to the $\kx=0$ mode. All other
small $\kx$ modes in the system must satisfy \eqref{eq:traditionalStabLimitII}. 
An analogous regularization for the AMP scheme involves the pressure equations
\[
   \frac{\partial^2 \ph^{(n+1)\pm}}{\partial y^2} +\alpha^{(n+1)} =0,
\]
and the single constraint
\[
   \int_{\Hf^{-}}^{-\hsHalf} \ph^{(n+1)-}dy+\int^{\Hf^{+}}_{\hsHalf} \ph^{(n+1)+}dy  = 0,
\]
which determines the one unknown $\alpha^{(n+1)}$.
In this case it can be shown that the AMP scheme is stable with no time-step restriction for $k_x=0$, and the  solution is found to be
\[
  \vh_1^{(n+1)\pm} = \vh_1^{0\pm}, \qquad \vh_2^{(n+1)\pm} =0, \qquad \etah^{n+1}=\etah^0, \qquad \ph^{(n+1)\pm}=\mp \Ks_0\etah^0\frac{\mathcal{D}^{\mp}}{\mathcal{D}^{+}+\mathcal{D}^{-}}.
\]
Note that, by adding one pressure regularization, the AMP scheme is able to determine the solution uniquely,
and the solution is compatible with the exact solution for $k_x=0$ given in the continuous case by~\eqref{eq:modelSolutionModeZero}.

\section{Numerical approach using deforming composite grids (DCG)} \label{sec:numericalApproach}

Our numerical approach for the solution of the equations governing an
FSI initial-boundary-value problem is based on the use of deforming
composite grids (DCG).  This FSI-DCG approach was first described
in~\cite{fsi2012} for the case of an inviscid
compressible flow coupled to a linearly elastic solid, and later 
in~\cite{flunsi2014r} for the case of compressible flow coupled to nonlinear
hyperelastic solids.  Here, we extend the approach to FSI problems involving an incompressible flow
coupled to deforming beams.



\subsection{Deforming composite grids and the fluid domain solver} \label{sec:dcg}

Deforming composite grids (DCGs) are used to discretize the evolving fluid domains in physical space.
An overlapping grid, $\Gs$, consists of a set of structured component grids, $\{G_g\}$,
$g=1,\ldots,{\mathcal N}$, that cover each fluid domain, $\OmegaF_k(t)$, and overlap where
the component grids meet. Typically, boundary-fitted curvilinear grids are used near the
boundaries while one or more
background Cartesian grids are used to handle the bulk of the fluid domain. 
Each component grid is a logically rectangular, curvilinear grid in $\nd$~space
dimensions, and is defined by a smooth
mapping from parameter space~$\rv$ (the unit square or cube) to physical
space~$\xv$,
\[
  \xv = \gv(\rv,t),\qquad \rv\in[0,1]^\nd,\qquad \xv\in\Real^\nd.
\]
Typically, the background grids are static, while the boundary-fitted grids evolve in time to
match the motion of the boundary.

{
\newcommand{\figWidth}{7.5cm}
\newcommand{\trimfig}[2]{\trimFig{#1}{#2}{.1}{.4}{.35}{.375}}
\begin{figure}[htb]
\begin{center}
\begin{tikzpicture}[scale=1]
  \useasboundingbox (0.0,.75) rectangle (16.,5.5);  
  \draw(0.0,0) node[anchor=south west,xshift=-4pt,yshift=+0pt] {\trimfig{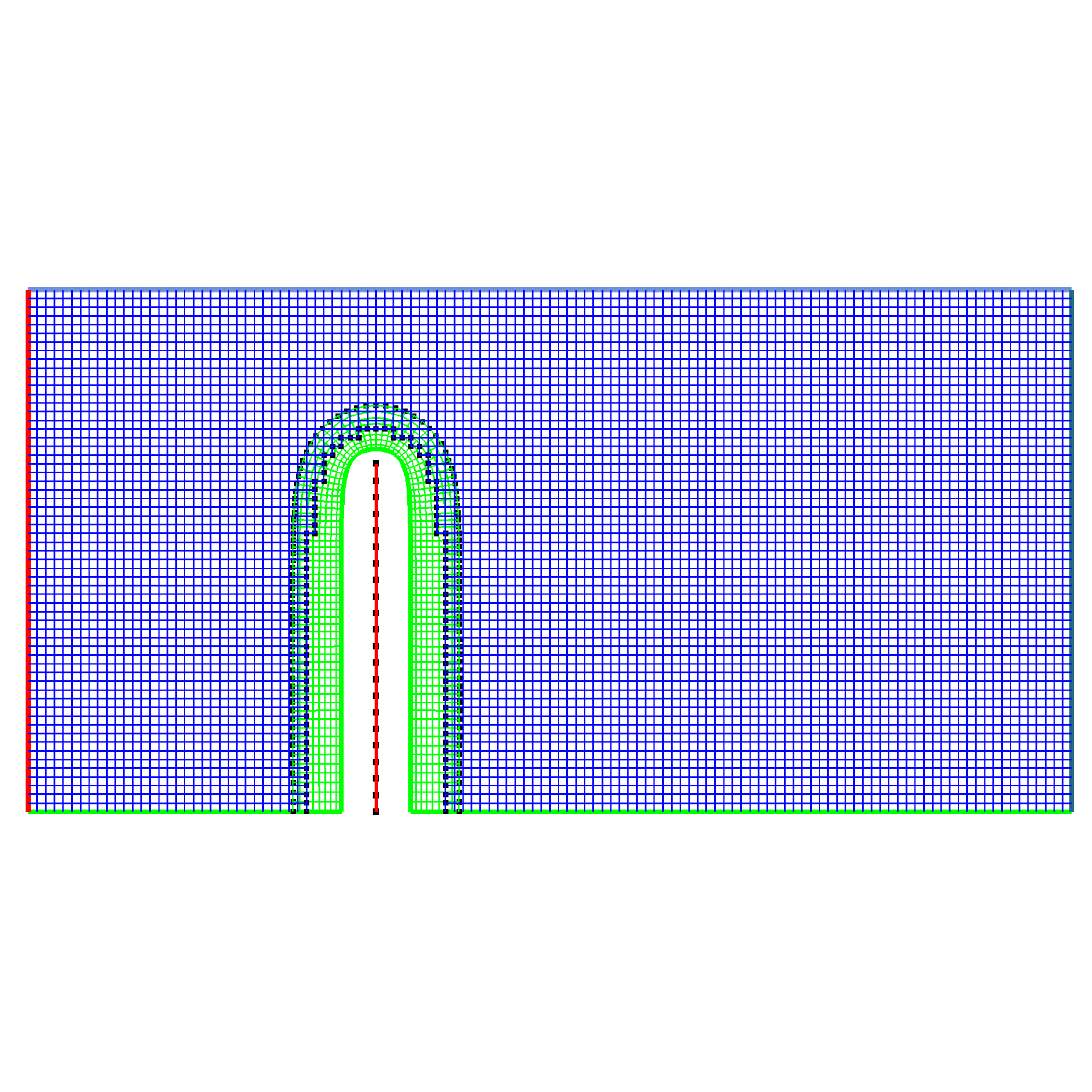}{\figWidth}};
  \draw(8.0,0) node[anchor=south west,xshift=-4pt,yshift=+0pt] {\trimfig{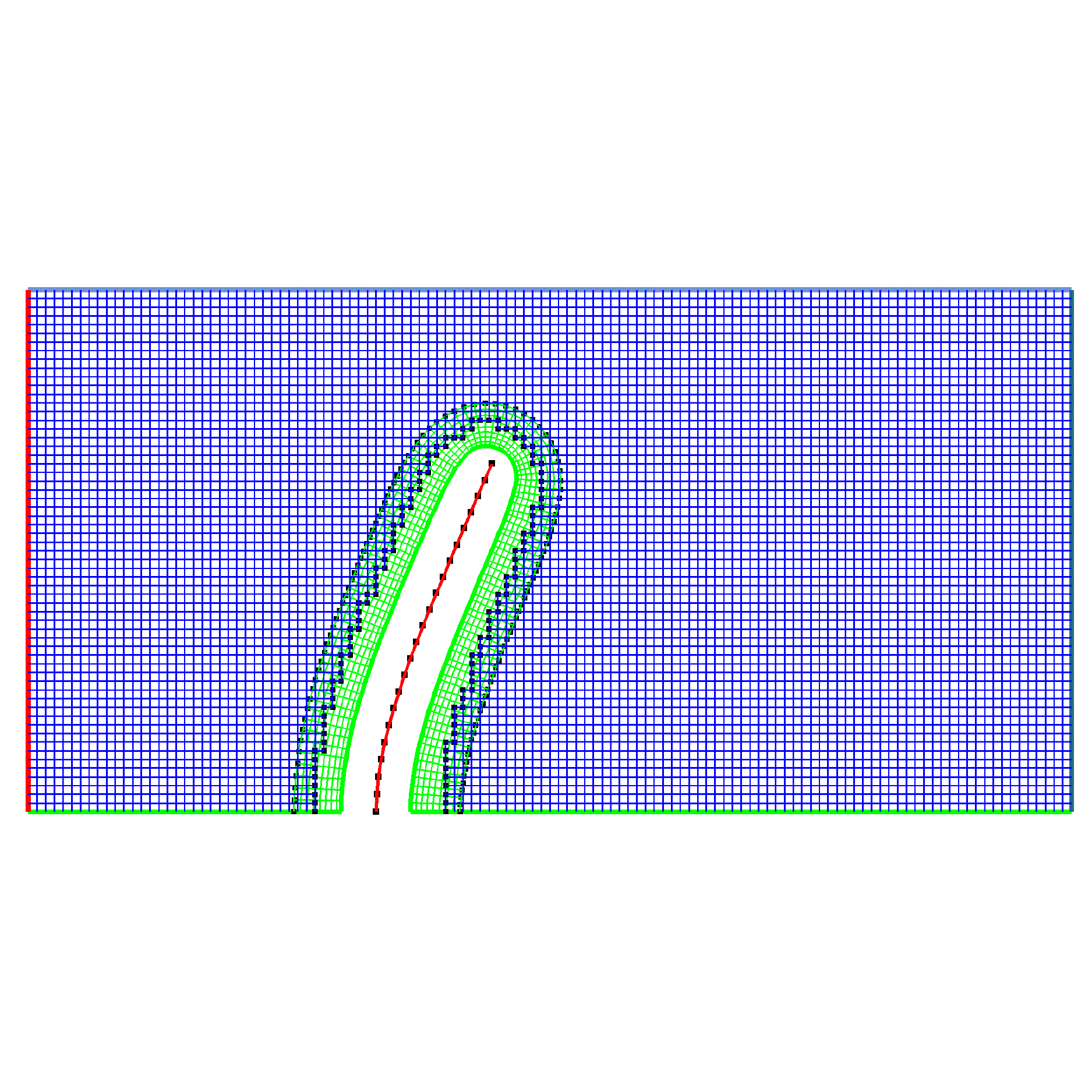}{\figWidth}};
\end{tikzpicture}
\end{center}
  \caption{Composite grids for a beam in a channel at times  $t=0.0$, and $1.0$. The green fluid grid deforms over
time to match the evolving beam (shown in white) and overlaps with the background Cartesian fluid grid (shown in blue).  The reference curve
of the beam is shown in red and the thickness of the beam describes the position of the surface of the beam.}
  \label{fig:beamGrids}
\end{figure}
}

In the FSI-DCG approach, component grids next to a fluid-structure interface deform over time
to match the beam motion. This is illustrated in \Fig~\ref{fig:beamGrids} for the case of a beam
in a fluid channel.
A green fluid grid is fitted to the boundary of the beam, which is defined in terms of displacement of the reference curve of the beam (shown in red) and the thickness of the beam (c.f.~\eqref{eq:EBsurface}).
The beam thickness is a function of the position $s$ along the reference curve, and remains fixed as the beam deforms.
%
After each time step, the points on the fluid interface corresponding to the surface of the beam are recomputed according to the current displacement of the beam and the normal to the reference curve using~\eqref{eq:EBsurface},
and a hyperbolic grid generator~\cite{HyperbolicGuide} is used to regenerate the local interface grids.
After all deforming component grids have been regenerated, the {\tt Ogen} grid generator~\cite{OGEN}
is called to regenerate the overlapping grid connectivity (e.g., cut holes, determine interpolation points) between the evolving interface-fitted fluid grids and the fixed background fluid grids, typically Cartesian grids, shown in blue in \Fig~\ref{fig:beamGrids}.

The incompressible Navier-Stokes equations are solved in velocity-pressure form 
on overlapping grids using a 
second-order accurate fractional-step algorithm~\cite{fis2014,splitStep2003,ICNS} similar to that
described in Section~\ref{sec:timeStepping}. The treatment of moving overlapping grids, including the 
assignment of {\em exposed points}, follows the approach
discussed in~\cite{mog2006} for compressible flow.
The governing equations for the fluid on a given component grid are solved in a coordinate frame moving with the grid. 
Consider the momentum equations in~\eqref{eq:fluidMomentum}
for the fluid in an Eulerian frame,
\begin{align}
  & \frac{\partial v_i }{\partial t} + v_j\frac{\partial v_i}{\partial x_j} = 
      {1\over\rho}\,\frac{\partial \sigma_{ji}}{\partial x_j}, \qquad i=1,\ldots,\nd,  \label{eq:fluidEulerian}
\end{align}
with summation convention. Under a general moving coordinate transformation, $\xv=\gv(\rv,t)$, 
the equations in~\eqref{eq:fluidEulerian} are transformed using the chain rule to become 
\begin{align}
  \frac{\partial v_i }{\partial t} 
     + (v_j-w_j)\frac{\partial r_k}{\partial x_j}\frac{\partial v_i}{\partial r_k} =
     {1\over\rho}\,\frac{\partial r_k}{\partial x_j}\frac{\partial \sigma_{ji}}{\partial r_k}, \qquad i=1,\ldots,\nd,\label{eq:GMmomentum}
\end{align}
where the component of velocity $v_i$ and stress $\sigma_{ji}$ are now considered to be functions of $\rv$ and $t$, and where
$w_j=\partial g_j/\partial t$ are the components of the {\em grid velocity}, $\wv$.
The pressure equation in~\eqref{eq:fluidPressure} is transformed to $(\rv,t)$ in a similar way.

\subsection{Structural domain solvers and the AMP condition}  \label{sec:beamSolverNumerical}

Two approaches are implemented to advance the governing equations for the beam.  The first is a second-order accurate predictor-corrector scheme described previously in Section~\ref{sec:timeStepping}, which also uses standard second-order accurate finite differences to approximate the spatial derivatives in~\eqref{eq:BeamModel} assuming an Euler-Bernoulli beam.  The basic AMP condition in~\eqref{eq:AMP_CombinedConditions} for this beam model becomes
\begin{equation}
    (\sigmav\nv)_+ + (\sigmav\nv)_- + 
\frac{\rhos\hs}{\rho}\grad\cdot\sigmav(\xsv_\pm(\ss,t),t)  
         = \Lsv\bigl(\uvs^{(p)},\vsv^{(p)}\bigr) + \rhos\hs \partialt \wsv^{(p)}_\pm(\ss,t) , \qquad \ss\in\OmegaSbar,
    \label{eq:AMP_CombinedConditionsII}
\end{equation}
where the predicted quantities on the right-hand side of~\eqref{eq:AMP_CombinedConditionsII} involving the displacement and velocity of the beam can be computed in terms of the discrete values for $\eta$, $\partial\eta/\partial t$, and $\partial^2\eta/\partial t^2$, and their spatial derivatives, using standard finite differences.  The normal and tangential components of~\eqref{eq:AMP_CombinedConditionsII} then become boundary conditions for the fractional-step fluid solver.

\newcommand{\fh}{f_h}
\newcommand{\fvh}{\fv_h}
\newcommand{\wvh}{\wwv_h}
\newcommand{\dwvh}{\dot\wwv_h}
\newcommand{\ddwvh}{\ddot\wwv_h}

The second approach implemented to advance the equations in~\eqref{eq:BeamModel} for an Euler-Bernoulli beam uses a standard finite-element approximation in space and the Newmark-beta time-stepping scheme (with parameters chosen for second-order accuracy in time).  The spatial representation of the beam displacement in the finite-element approximation is taken as 
\begin{align}
 \ww_h(s,t) = \sum_{j=1}^{\Nn}\Big\{\, \ww_j(t) \phi_j(s) +\ww_j'(t) \psi_j(s)  \,\Big\}, \label{eq:beamFEM}
\end{align}
where the degrees of freedom are the
nodal displacement $\ww_j(t)$ and slope $\ww\sp\prime_j(t)$, and where $\phi_j(s)$ and $\psi_j(s)$ are 
cubic Hermite polynomials 
with support on the interval $[s_{j-1},s_{j+1}]$ and satisfying
\[
\begin{array}{ll}
  \phi_i(s_j)=\delta_{ij}, \quad &\phi_i'(s_j)=0, \medskip\\
  \psi_i(s_j)=0,                  &\psi_i'(s_j)=\delta_{ij},
\end{array}
\]
where $\delta_{ij}$ is the Kronecker delta.  Following the usual Galerkin approximation in space, the beam equation become
\begin{align}
   M_h \frac{d^2\wvh}{dt^2}  = - K_h \wvh - B_h \frac{d\wvh}{dt}  + \fvh, \label{eq:discreteBeam}
\end{align}
where $M_h$ is a mass matrix, $K_h$ is a stiffness matrix, and $B_h$ is a damping matrix.  The vector $\wvh$ contains the nodal displacements and slopes, $\ww_j(t)$ and $\ww\sp\prime_j(t)$, $j=1,\ldots,\Nn$, and the components of the forcing vector $\fvh$ are given by integrals of $f(\sv,t)$ in~\eqref{eq:BeamModel} with the Hermite basis functions\footnote{Appropriate adjustments are made to the matrices and right-hand-side to account for boundary conditions.}.
The ODEs in~\eqref{eq:discreteBeam} are advanced in time using a second-order accurate Newmark-beta scheme.

\newcommand{\beamInnerProd}[1]{\left(#1\right)_{\OmegaSbar}} \newcommand{\partialx}{\partial_x}

Given predicted values for $\wvh^{(p)}$ and $d\wvh^{(p)}/dt$ from the first stage of the AMP algorithm, pointwise values
of the beam operator, $\Lsv\sp{(p)}$, are required to evaluate the right-hand side of the AMP condition
in~\eqref{eq:AMP_CombinedConditionsII}.  
This is generally a straightforward computation in the context of finite difference or
finite volume approximations to the beam operator, whereas the calculation may be more difficult within the finite-element
approximation.  As an example, consider a beam operator of the form
\[
   \Ls = - (\Es\Is\ww_{ss})_{ss},
\]
where the subscripts denote partial derivatives with respect to $s$.  It is sufficient to only consider the highest derivative term in~\eqref{eq:BeamModel} for the purposes of the present discussion.
Let $\Ls$ have the Hermite finite-element representation
\begin{equation}
 \Ls_h(s,t) = \sum_{j=1}^{\Nn} \Big\{\, \gamma_j(t) \phi_j(s) + \gamma_j'(t) \psi_j(s) \,\Big\}  ,\label{eq:beamOperatorAMP_APPROX}
\end{equation}
where $\gamma_j(t)$ and $\gamma_j'(t)$ are coefficient functions of time, and $\phi_j(s)$ and $\psi_j(s)$ are the cubic polynomial basis functions in~\eqref{eq:beamFEM}.  
The $2\Nn$ unknowns, $\gamma_j(t)$ and $\gamma_j'(t)$, are determined from 
the $2\Nn$ Galerkin equations
\begin{equation}
    \beamInnerProd{ \chi_j,\Ls_h}  = (\chi_j, 
        - \beamInnerProd{\Es\Is\ww_{h,ss})_{ss} }, \qquad j=1,2,\ldots,\Nn, \label{eq:galerkin}
\end{equation}
where $\chi_j=\phi_j$ or $\psi_j$, and where $\beamInnerProd{f,g}\equiv \int_0^{\ls} f(s) g(s)\, ds$. 
Integration by parts (twice) gives 
\begin{equation}
    \beamInnerProd{\chi_j,\Ls_h} =  
       - \beamInnerProd{\chi_{j,ss},\Es\Is\ww_{h,ss}} + \left.\Big\{ 
            - \chi_j (\Es\Is\ww_{h,ss})_{s} + \chi_{j,s} \Es\Is\ww_{h,ss})\Big\}\right\vert_0^{\ls}, \qquad j=1,2,\ldots,\Nn .
\label{eq:integrationByParts}
\end{equation}
The problem now arises in computing second-order accurate approximations to the boundary terms $(\Es\Is\ww_{h,ss})_{s}$ and $\Es\Is\ww_{h,ss}$
appearing in~\eqref{eq:integrationByParts}. The cubic Hermite representation is not of
sufficient accuracy for this purpose since, for example, $\ww_{h,sss}$ is constant over an element and at best first-order accurate. 

There are various possible solutions to the difficulty of evaluating the right-hand side of~\eqref{eq:AMP_CombinedConditionsII}
in the context of a finite-element approximation, 
and we have considered and tested two different approaches that require relatively minor changes. 
The properties of the two approaches are evaluated in Section~\ref{sec:numericalResults} where some
conclusions are drawn.

\subsubsection{AMP-PBA : Predict Beam Acceleration} \label{sec:ampadj}

 In the first approach, which we refer to as AMP-PBA (Predict Beam Acceleration), the vector form of~\eqref{eq:BeamModel} is used in~\eqref{eq:AMP_CombinedConditionsII} to eliminate $\Lsv\sp{(p)}$, which gives an AMP condition of the form,
\begin{equation}
    (\sigmav\nv)_+ + (\sigmav\nv)_- + 
\frac{\rhos\hs}{\rho}\grad\cdot\sigmav(\xsv_\pm(\ss,t),t)  
         =  \rhos \hs \frac{\partial^2}{\partial t^2}\usv^{(p)} - \fsv\sp{(p)}(s,t) + \rhos\hs \partialt \wsv^{(p)}_\pm(\ss,t) , \qquad \ss\in\OmegaSbar.
    \label{eq:AMP_CombinedConditionsIIadjusted}
\end{equation}
The predicted beam acceleration 
$\partial\sp2\usv\sp{(p)}/\partial t\sp2$ in~\eqref{eq:AMP_CombinedConditionsIIadjusted} is readily computed from the acceleration term in the weak form~\eqref{eq:discreteBeam}, while the predicted fluid forcing, $\fsv\sp{(p)}(s,t)=-(\sigmav\nv)_+^{(p)}-(\sigmav\nv)_-^{(p)}$, can be obtained from the discrete fluid stress.  This gives the following AMP-PBA algorithm:
\begin{enumerate}
   \item Obtain predicted values for $\fsv\sp{(p)}=-(\sigmav\nv)_+^{(p)}-(\sigmav\nv)_-^{(p)}$ 
         using extrapolation in time during the predictor stage or the current guess during the corrector stage (c.f.~the algorithm in Section~\ref{sec:timeStepping}).
   \item Solve~\eqref{eq:discreteBeam} for $d^2\wvh^{(p)}/dt^2$ using the predicted force $\fsv\sp{(p)}$, 
        and evaluate pointwise values for 
        $\rhos \hs (\partial\sp2\usv\sp{(p)}/\partial t\sp2)$
          from the Hermite representation for this quantity. 
   \item Evaluate the right-hand side to~\eqref{eq:AMP_CombinedConditionsIIadjusted} and use this when applying the AMP conditions to the fluid.
\end{enumerate}

\subsubsection{AMP-PBF : Predict Beam internal Force}

In the second approach, referred to as AMP-PBF (Predict Beam internal Force), a mixed finite-element/finite-difference approach is used to compute the beam internal force
$\Lsv\sp{(p)}$ directly from equation~\eqref{eq:integrationByParts} (extended to include all terms in the full EB model). 
In particular, using the solution nodal values and derivatives on the boundary and two adjacent points, a
fourth-order accurate approximation to $\Es\Is\ww_{ss}$ and a second-order accurate approximation to $(\Es\Is\ww_{ss})_{s}$ on the boundary can be
derived. This allows the coefficients $\gamma_j(t)$ and $\gamma_j'(t)$
to be determined and then pointwise values of $\Lsv\sp{(p)}$ can be computed from~\eqref{eq:beamOperatorAMP_APPROX}.  
This gives the following AMP-PBF algorithm:
\begin{enumerate}
   \item Given predicted values for $\ww_h^{(p)}$, determine second-order accurate approximations to the boundary terms $(\Es\Is\ww_{h,ss})_{s}$ 
       and $\Es\Is\ww_{h,ss}$ in~\eqref{eq:integrationByParts} using a finite-difference approximation involving the boundary   
       points and two adjacent points. 
   \item Compute values for $\Lsv^{(p)}$ by solving the Galerkin equations in~\eqref{eq:integrationByParts} 
         (extended to include all terms in the full EB model) for
         $\gamma_j$ and $\gamma_j'$. 
   \item Evaluate the right-hand side to~\eqref{eq:AMP_CombinedConditionsII} using~\eqref{eq:beamOperatorAMP_APPROX},
         and use this when applying the AMP conditions to the fluid.
\end{enumerate}

\subsection{Filtering the projected interface velocity} \label{sec:velocityFilter}

The projected interface velocity~\eqref{eq:velocityProjection}
is optionally smoothed using a few iterations of a fourth-order filter as described in~\cite{smog2012}.
The beam solution can also be smoothed to prevent high-frequency numerical oscillations, 
for example by choosing the coefficient $\Kxxt$ in~\eqref{eq:BeamModel} to be proportional to $\Delta\sv^2$. However, for light
beams the interface velocity is determined primarily from the fluid and thus the 
filter on the interface velocity is needed for light beams to prevent numerical oscillations on the interface.
Note that as the beam becomes very light,
with beam parameters $\rhos\hs$, $\Ts$, $\Es\Is$, etc. going to zero, the AMP interface
condition approaches that of a free surface boundary condition for the fluid.
Such a free surface condition with no surface tension is susceptible to physical and numerical instabilities, and the fourth-order filter acts to suppress these instabilities and thus serves as a {\em numerical surface tension}. 
This filter does not affect the overall second-order accuracy of the AMP FSI time-stepping scheme.

%
%
%
%
%
%
%
\subsection{Time step determination} \label{sec:timeStep}

The global time step, $\dt$, for the AMP FSI time-stepping scheme is chosen to be the 
minimum of the stable time steps for the fluid domain solver and beam solver. 
The $\dt$  for the fluid solver is chosen to be a factor 
$\lambda_{{\rm SF}}$ (``SF'' for stability factor)
times an estimate of the 
largest stable time step based on a von Neumann analysis of the linearized Navier-Stokes system.
Typically we choose $\lambda_{{\rm SF}}=0.9$.
In this work, the FEM beam solver is implicit in time and stable for any time step. However, for 
accuracy reasons it is convenient to choose $\dt$ based on a condition relating to explicit time-stepping. 
Thus, the time-step for the beam is chosen as a factor $\bar\lambda_{{\rm SF}}$ 
times an estimate of the largest time step such that an explicit time-stepping procedure for the beam would 
satisfy a CFL stability condition.  Here $\bar\lambda_{{\rm SF}}$  is chosen in the range~1 to~10.

%
%
%
%
%
%
%
%
%

\section{Numerical Results} \label{sec:numericalResults}

We now present the results for a series of simulations chosen to demonstrate the properties of our
numerical approach for FSI problems based on the use of the AMP interface conditions. In order to illustrate the basic properties
of the approach in a simple setup, we begin with two cases that use a single
deforming fluid grid coupled to an Euler-Bernoulli (EB) beam with negligible thickness and finite mass. This 
simplification eliminates the finite-thickness corrections
given in~\eqref{eq:EBsurface} and~\eqref{eq:beamSurfaceVelocityEB} by setting $\hs=0$. The first test involves a time-dependent problem in which an exact solution is known following the method of manufactured solutions, while the second test considers a steady state FSI problem with known analytical solution.  The solutions for both tests are obtained using the beam solver based on finite differences, and they illustrate the accuracy and stability of our numerical approach for a relatively simple geometric configuration.  The next set of tests consider finite-thickness EB beams in more complex geometries.  Solutions of these FSI problems are obtained using deforming composite grids and the beam solver based on finite elements as discussed in Section~\ref{sec:numericalApproach}.  These problems include the dynamics of a flexible beam separating two fluid chambers, fluid flow in a channel
with flexible walls, and finally the deformation of a flexible cantilevered beam in a cross fluid flow.

\newcommand{\etatz}{\tilde{\eta}}
\newcommand{\vtz}{\tilde{v}}
\newcommand{\ustz}{\tilde{u}}
\newcommand{\ptz}{\tilde{p}}
\newcommand{\fx}{f_x}
\newcommand{\ft}{f_t}
\newcommand{\pd}[2]{\frac{\partial #1}{\partial #2}} 
\newcommand{\pdn}[3]{\frac{\partial^{#3} #1}{\partial #2^{#3}}} 
\newcommand{\tableTZAmp}{%
  \begin{figure}[hbt]\tableFont %
    \begin{center}
    
\begin{tabular}{|c|c|c|c|c|c|c|c|c|c|c|} \hline 
\multicolumn{11}{|c|}{\strutt  light beam, $\rhos\hs=10\sp{-3}$} \\ \hline 
\strutt~~$h_j$~~& $E_j^{(p)}$ & $r$ & $E_j^{(v_1)}$ & $r$ & $E_j^{(v_2)}$ & $r$ & $E_j^{(\eta)}$ & $r$   & $E_j^{(\eta_t)}$  & $r$  \\ \hline 
1/10 & 6.55e-02 &    & 2.64e-02  &     & 2.23e-02 &     & 9.06e-04  &     & 2.19e-02  &      \\ \hline
1/20 & 1.70e-02 &  3.84  & 7.26e-03  &  3.64  & 5.41e-03 &  4.12  & 2.23e-04  &  4.07  & 5.15e-03  &  4.25   \\ \hline
1/40 & 4.45e-03 &  3.83  & 1.86e-03  &  3.91  & 1.39e-03 &  3.90  & 6.03e-05  &  3.69  & 1.31e-03  &  3.94   \\ \hline
1/80 & 1.12e-03 &  3.98  & 4.70e-04  &  3.95  & 3.55e-04 &  3.91  & 1.54e-05  &  3.90  & 3.29e-04  &  3.97   \\ \hline
rate &  1.96 &   &  1.94  &   &  1.99  &   &  1.95  &   &  2.01  &    \\ \hline
\end{tabular}

\medskip

\begin{tabular}{|c|c|c|c|c|c|c|c|c|c|c|} \hline 
\multicolumn{11}{|c|}{\strutt  medium beam, $\rhos\hs=1$} \\ \hline 
\strutt~~$h_j$~~& $E_j^{(p)}$ & $r$ & $E_j^{(v_1)}$ & $r$ & $E_j^{(v_2)}$ & $r$ & $E_j^{(\eta)}$ & $r$   & $E_j^{(\eta_t)}$  & $r$  \\ \hline 
1/10 & 1.18e-01 &     & 2.07e-02  &     & 2.40e-02 &     & 3.56e-04  &     & 9.09e-03  &      \\ \hline
1/20 & 2.53e-02 &  4.69  & 5.86e-03  &  3.53  & 5.80e-03 &  4.14  & 7.70e-05  &  4.63  & 1.89e-03  &  4.80   \\ \hline
1/40 & 6.49e-03 &  3.90  & 1.54e-03  &  3.80  & 1.52e-03 &  3.82  & 2.02e-05  &  3.81  & 4.61e-04  &  4.11   \\ \hline
1/80 & 1.67e-03 &  3.89  & 3.95e-04  &  3.90  & 3.87e-04 &  3.92  & 5.10e-06  &  3.96  & 1.15e-04  &  4.02   \\ \hline
rate &  2.04 &   &  1.91  &   &  1.98  &   &  2.03  &   &  2.10  &    \\ \hline
\end{tabular}

\medskip

\begin{tabular}{|c|c|c|c|c|c|c|c|c|c|c|} \hline 
\multicolumn{11}{|c|}{\strutt  heavy beam, $\rhos\hs=10\sp{3}$} \\ \hline 
\strutt~~$h_j$~~& $E_j^{(p)}$ & $r$ & $E_j^{(v_1)}$ & $r$ & $E_j^{(v_2)}$ & $r$ & $E_j^{(\eta)}$ & $r$   & $E_j^{(\eta_t)}$  & $r$  \\ \hline 
1/10 & 1.73e-01 &     & 2.06e-02  &     & 2.80e-02 &     & 8.45e-05  &     & 2.91e-03  &      \\ \hline
1/20 & 3.33e-02 &  5.20  & 5.63e-03  &  3.66  & 6.68e-03 &  4.20  & 2.25e-05  &  3.76  & 7.72e-04  &  3.77   \\ \hline
1/40 & 8.83e-03 &  3.77  & 1.51e-03  &  3.74  & 1.77e-03 &  3.77  & 6.37e-06  &  3.53  & 1.91e-04  &  4.05   \\ \hline
1/80 & 2.31e-03 &  3.82  & 3.90e-04  &  3.86  & 4.54e-04 &  3.90  & 1.63e-06  &  3.92  & 4.75e-05  &  4.01   \\ \hline
rate &  2.06 &   &  1.91  &   &  1.98  &   &  1.89  &   &  1.98  &    \\ \hline
\end{tabular}

\caption{Maximum-norm errors and estimated convergence rates using manufactured solutions. The column labeled ``r'' provides the
  ratio of the errors at the current grid spacing to that on the next coarser grid. }
\label{tab:tableTZAmp}
\end{center}
\end{figure}
}

\subsection{Verification results for a beam with negligible thickness}

\newcommand{\xx}{x}
\newcommand{\yy}{y}

We first consider solutions of two FSI problems involving the interaction of an incompressible fluid and a beam of negligible thickness and finite mass.  The geometric configuration of both problems is illustrated in Figure~\ref{fig:lightBeamError}.  The fluid domain is given by $\OmegaF(t) =[0,1]\times[0,1+\eta(\xx,t)]$, where $\eta(\xx,t)$ is the vertical displacement of the Euler-Bernoulli (EB) beam which forms the upper surface of the fluid domain.  The fluid domain is represented by a single transfinite interpolation (TFI) grid which evolves in time according to the computed solution for $\eta$.  The TFI mapping is given by
\[
(\xx,\yy) = \mathbf{g}(\rv,t)\equiv(r_1,(1+\eta(r_1,t))r_2),\qquad \rv\in [0,1]\times[0,1].
\]
The equations for the fluid are solved in the mapped domain on the unit square with uniform grid spacings $\Delta r_1=\Delta r_2$.  The equations for the EB beam are solved for $\sv=r_1\in[0,1]$ using the same grid spacing as the fluid grid in the $r_1$~direction.

\subsubsection{Manufactured Solutions} \label{sec:manufacturedSolutions}

The method of manufactured solutions~\cite{CGNS,Roache1998} can be used to construct exact
solutions of FSI problems by adding forcing functions to the governing
equations. The forcing is specified so that a chosen function becomes an exact solution to the forced equations. Here the approach is used to verify the stability and accuracy of the AMP scheme for the nonlinear 
FSI problem described above consisting of an initially rectangular fluid domain bounded by an EB beam along its top surface.
 
The parameters for the fluid are chosen to be $\rho=1$ and $\mu=0.05$.  The thickness of the beam is taken to be zero while the finite mass of the beam given by $\rhos\hs$ is varied to consider the cases of light, medium and heavy beams.  The remaining parameters of the beam are assumed to be $\Ks=1$, $\Ts=\rhos\hs$ and $\Es\Is=\Kt=\Kxxt=0$.  The exact displacement of the beam is chosen to be
\begin{equation}
\eta_e(\xx,t) = \frac{a}{\pi f_t}\sin(f_x\pi \xx)\sin(f_t\pi t),
\label{eq:tzsolid}
\end{equation}
where $a=0.5$ is an amplitude and $f_x=f_t=2$ are frequencies.  (Note that $s=\xx=r_1$ for this problem.)  The exact components of velocity and the pressure in the fluid are taken to be
\begin{equation}
\left.
\begin{array}{l}
\displaystyle{
v_{1,e}(\xv,t) = -a\cos(f_x\pi \xx)\sin(f_x\pi\tilde \yy)\cos(f_t\pi t)},\medskip\\
\displaystyle{
v_{2,e}(\xv,t) = a\sin(f_x\pi \xx)\cos(f_x\pi\tilde \yy)\cos(f_t\pi t)-a\pd{\eta_e}{\xx}\cos(f_x\pi \xx)\sin(f_x\pi\tilde \yy)\cos(f_t\pi t)},\medskip\\
\displaystyle{
p_e(\xv,t) = \cos(f_x\pi \xx)\cos(f_x\pi \yy)\cos(f_t\pi t)},
\end{array}
\right\}
\label{eq:tzfluid}
\end{equation}
where $\tilde \yy=\yy-(1+\eta_e)$.  
The components of velocity have been chosen so that the continuity equation is satisfied in the exact solution, and the vertical velocities of the fluid and the beam match along the interface at $\yy=\eta_e$.  Forcing functions are added to the equations for the fluid and structure so that the functions in~\eqref{eq:tzsolid} and~\eqref{eq:tzfluid} are exact solutions.

The stability of the AMP scheme for light beams is illustrated in figure~\ref{fig:lightBeamError} for the case $\rhos\hs=10^{-3}$.  The image on the left shows the TFI grid for the fluid domain at $t=0.1$.  The beam lies on the top surface of the domain, while the boundary conditions on the other three sides of the domain are taken to be no-slip conditions on the left and right sides ($\vv=\vv_e$) and Dirichlet conditions on the bottom ($\vv=\vv_e$ and $p=p_e$).  The plot in the center of the figure shows the numerical solution for the vertical component of the fluid velocity, $v_2$, while the plot on the right shows the difference between the computed solution for $v_2$ and the exact solution given by $v_{2,e}$.  Here we observe that the error in $v_2$ is well behaved in that the magnitude is small and it is smooth throughout the domain including the boundaries.  The behavior of the error in the other components of the solution are similar.
 
%
%
%
{
\newcommand{\figWidth}{5.5cm}
\newcommand{\cbWidth}{.2cm}
\newcommand{\cbHeight}{4.1cm}
\newcommand{\trimfigcb}[3]{\includegraphics[width=#2, height=#3, clip, trim=17cm 2.35cm 1.65cm 2.35cm]{#1}}
\newcommand{\trimfig}[2]{\trimFig{#1}{#2}{.2}{0.12}{.5}{.5}}
\begin{figure}[htb]
\begin{center}
\begin{tikzpicture}[scale=1]
\useasboundingbox (0.0,.75) rectangle (18.,6.);  
\draw(-0.5,0.0) node[anchor=south west,xshift=-4pt,yshift=+0pt] {\trimfig{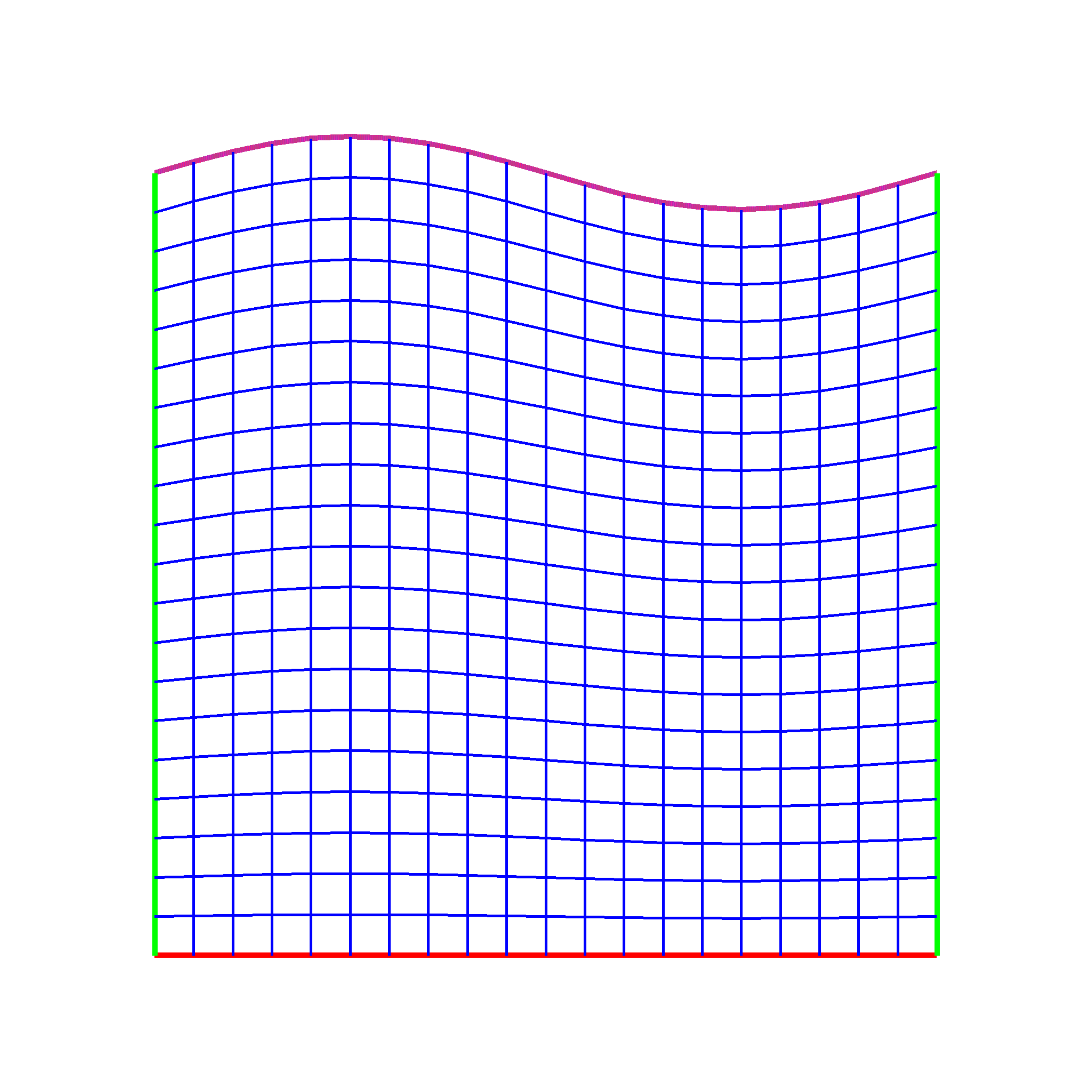}{\figWidth}};
\draw(0.1,6.) node[draw,fill=white,anchor=west,xshift=2pt,yshift=-16pt] {\scriptsize Grid at $t=0.1$};
\draw(4.5,0.0) node[anchor=south west,xshift=-4pt,yshift=+0pt] {\trimfig{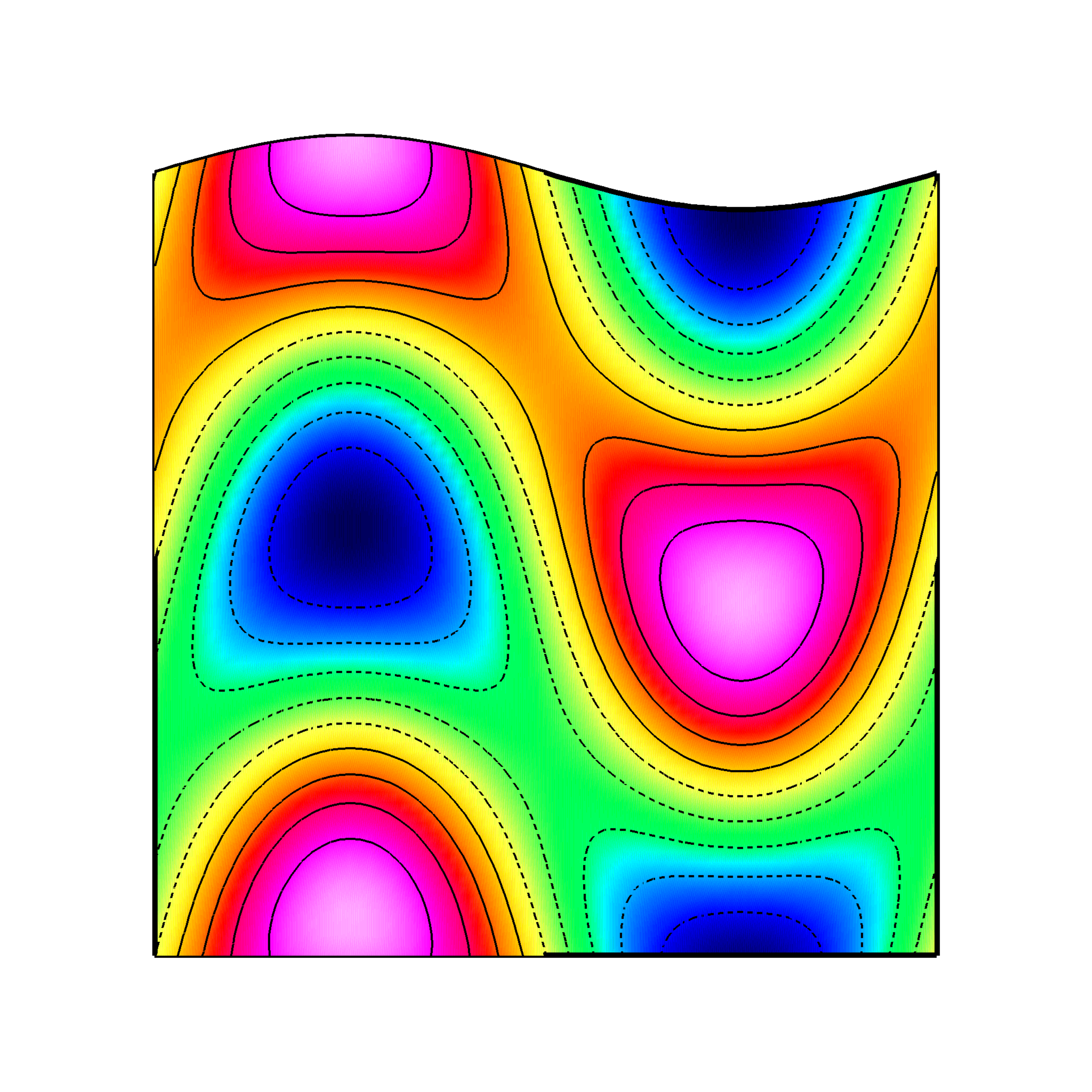}{\figWidth}};
\draw (9.25,0.3) node[anchor=south west,xshift= +0pt,yshift=+0pt] {\trimfigcb{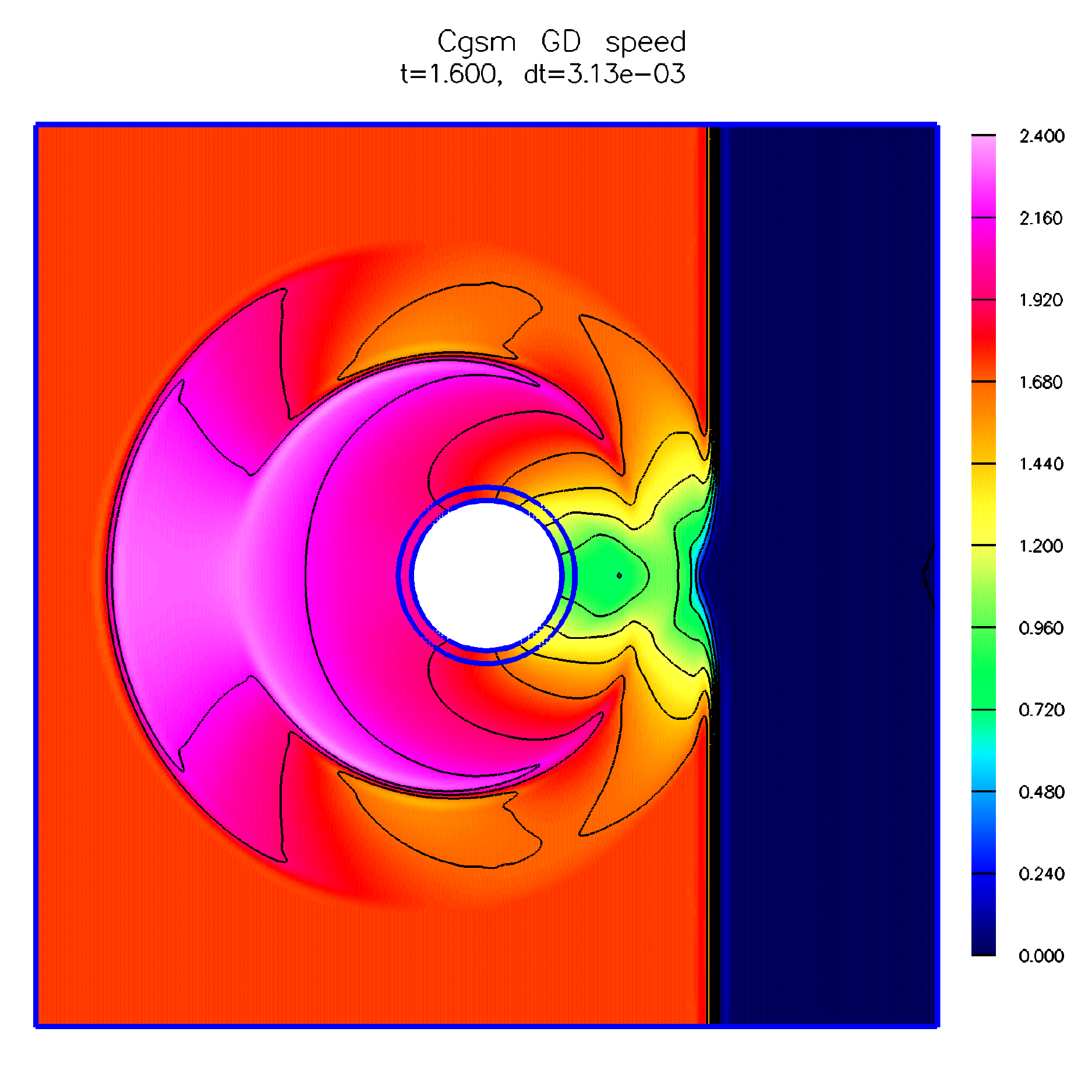}{\cbWidth}{\cbHeight}};
\draw (9.25,0.3) node[anchor=south west,xshift= +8pt,yshift=+1pt] {\scriptsize $-0.4$};
\draw (9.25,4.5) node[anchor=south west,xshift= +8pt,yshift=-6pt] {\scriptsize $0.4$};
\draw(5.1,6.) node[draw,fill=white,anchor=west,xshift=2pt,yshift=-16pt] {\scriptsize $v_2$ at $t=0.1$};
\draw(10.25,0.0) node[anchor=south west,xshift=-4pt,yshift=+0pt] {\trimfig{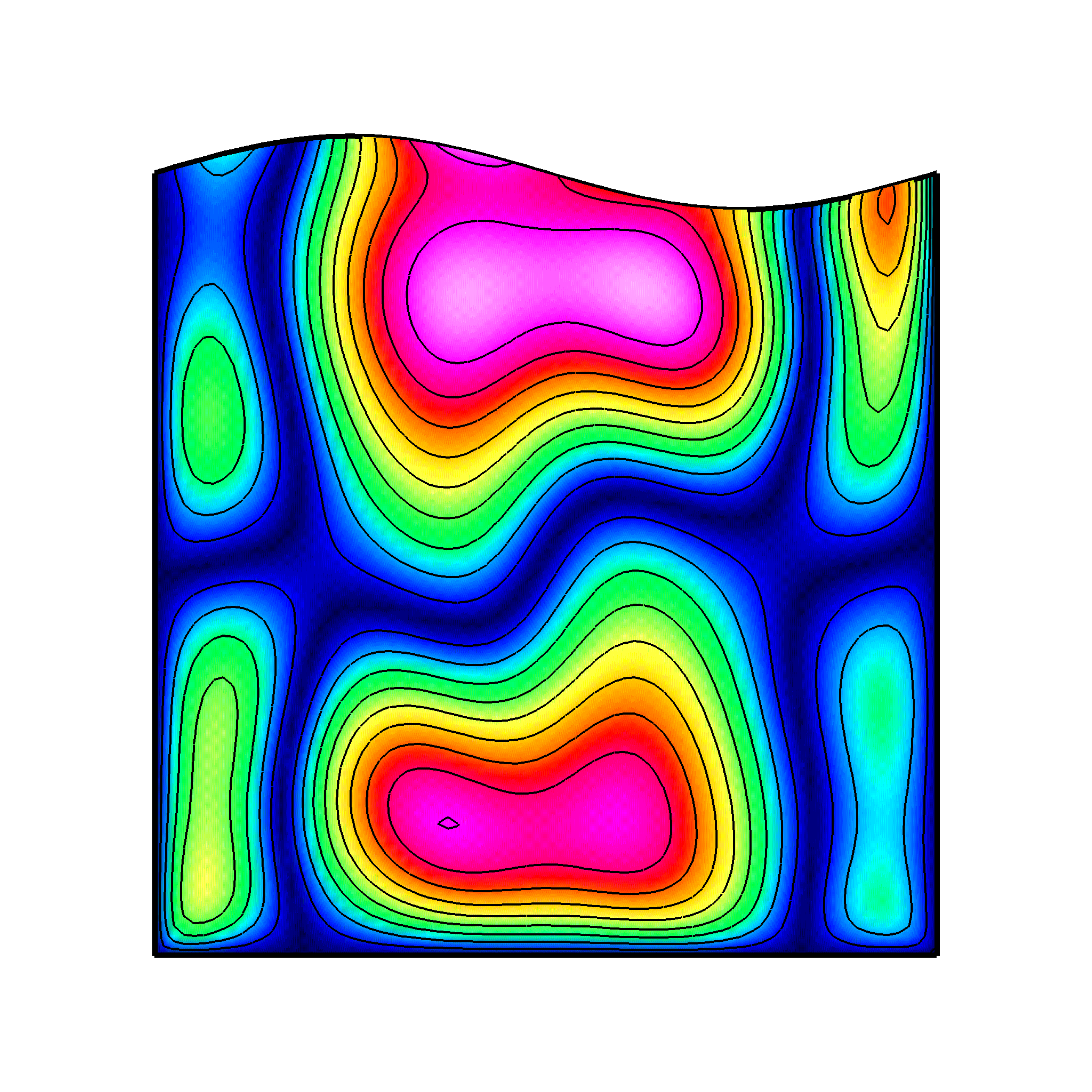}{\figWidth}};
\draw (15.,0.3) node[anchor=south west,xshift= +0pt,yshift=+0pt] {\trimfigcb{fig/colourBarLines}{\cbWidth}{\cbHeight}};
\draw (15.,0.3) node[anchor=south west,xshift= +8pt,yshift=+1pt] {\scriptsize $0$};
\draw (15.,4.5) node[anchor=south west,xshift= +8pt,yshift=-6pt] {\scriptsize $3.5\text{e}{-4}$};
\draw(10.85,6.) node[draw,fill=white,anchor=west,xshift=2pt,yshift=-16pt] {\scriptsize Error in $v_2$ at $t=0.1$};

\end{tikzpicture}
\end{center}
  \caption{Manufactured solutions: computational results at $t=0.1$ for a light beam ($\rhos\hs=10^{-3}$) 
using a grid spacing of $h=1/80$.
Left: deformed grid (coarsened version). Middle: solution component $v_2$. Right: error in $v_2$.
}
  \label{fig:lightBeamError}
\end{figure}
\let\cbWidth\undefined
\let\cbHeight\undefined
\let\trimfigcb\undefined
}

%

A convergence study is shown in  Table~\ref{tab:tableTZAmp} for light ($\rhos\hs=10^{-3}$), medium ($\rhos\hs=1$) and heavy ($\rhos\hs=10^{3}$) beams using grid spacings $\Delta r_{1,j}=\Delta r_{2,j}=h_j=1/(10j)$,
$j=1,2,\ldots$.  The maximum-norm errors, $E_j^{(q)}$ for solution component $q$ and for grid resolutions $j=1$, 2, 4 and 8, are computed for each value of $\rhos\hs$.  The convergence
rate is estimated by a least squares fit to the logarithm of the error over the grid resolutions. 
The results in the table for the components of the fluid velocity $(v_1,v_2)$ and pressure $p$, and for the beam displacement $\eta$ and velocity $\eta_t$ show the expected second-order accuracy for all three choices of the beam mass, and in particular for the case of a light beam when added-mass effects are strong.

\tableTZAmp

%
%


%
%
%
%
%
%

\newcommand{\tableParabolicBeamAmp}{%
  \begin{figure}[hbt]\tableFont %
    \begin{center}
      \begin{tabular}{|c|c|c|c|c|c|c|} \hline 
          \strutt~~$\bar{\rho}{\bar{h}}$~~& $t_{{\rm ss}}$ & $E_j^{(p)}$  & $E_j^{(v_1)}$  & $E_j^{(v_2)}$  & $E_j^{(\eta)}$    & $E_j^{(\eta_t)}$    \\ \hline
        1e-3 & 30  &  1.50e-12  & 7.06e-14  & 2.25e-13   & 2.00e-15 & 2.25e-13     \\ \hline
        1     &  50  & 4.53e-12  &  3.73e-14   & 7.15e-13   & 1.80e-13 & 7.15e-13     \\ \hline
        10 & 250  & 4.19e-12  & 3.84e-14   & 2.63e-12  & 9.08e-14 &2.63e-12   \\ \hline
      \end{tabular}
\caption{Maximum-norm errors at near steady state when $t=t_{{\rm ss}}$ for light ($\rhos\hs=10^{-3}$), medium ($\rhos\hs=1$) and heavy ($\rhos\hs=10$) beams.}
\label{tab:tableParabolicBeamAmp}
\end{center}
\end{figure}
}

\newcommand{\Gcpc}{\Gc_{\text{fp}}}
\subsubsection{Steady state beam bounding a pressurized fluid chamber} \label{sec:flexibleChannelSS}

To illustrate the effectiveness of the AMP coupling scheme in comparison to the
traditional FSI coupling, we now discuss a problem where the beam deforms as a result of a uniform 
fluid pressure applied at the bottom of the fluid domain.  For this problem, an exact solution for the
long-time steady state deflection of the beam can be found.  The steady state can be computed numerically by
integrating the time-dependent equations from an initially flat beam.
It is found that for {\em light} beams the traditional scheme is unstable, in agreement with the previous stability analysis, while the
AMP scheme remains stable and the steady state solution is computed accurately. 



For this FSI problem, the material parameters for the fluid are chosen 
as  $\rho=1$ and $\mu=0.2$, while those for the
beam are chosen as $\Ks=0$, $\Ts=5$, $\Es\Is=0$, $\Kt=\Kxxt=0$ and a finite value for $\rhos\hs$.
The boundary conditions on the left and right sides of the fluid domain are given as a no-slip walls ($\vv=0$), 
while a given pressure ($p=P_a=5$) and zero tangential velocity ($v_1=0$) are assigned on the bottom of the fluid domain.
One boundary condition is required at each end of the beam and these are taken as
$\eta=0$ at each end.  (The highest spatial derivative in the beam model is two for the chosen parameters.)
The initial conditions are given as a flat stationary beam ($\eta=\dot\eta=0$) and 
a fluid at rest. The exact steady-state solution is
\begin{align}
\vv=0, \qquad p = P_a,\qquad \eta = \frac{P_a \xx(1-\xx)}{2\Ts},\qquad \xx\in[0,1].
\end{align} 
As before, $\rhos\hs$ is varied to illustrate the stability of the time-integration of the equations using the AMP coupling scheme for a range of masses of the beam, and also to highlight the instability of the traditional FSI coupling scheme for light beams. 
%
%

Figure~\ref{fig:ParabolicBeam} shows the behavior of the displacement and velocity of the beam at
its center, $\xx=1/2$, as a function of time for the cases of a light ($\rhos\hs=10^{-3}$), medium ($\rhos\hs=1$) and heavy ($\rhos\hs=10$) beam.
For each case, there is an early transient stage consisting of oscillations of the beam, but these decay as the beam approaches a steady state.
Since the damping is supplied by the viscous dissipation in the fluid, the oscillations decay more rapidly for lighter beams.
At very long times, the solutions for beam displacement and velocity approximate the exact steady state solutions.  Since the steady
state solution is quadratic in $\xx$, we expect the agreement to be close to exact. The table in Figure~\ref{tab:tableParabolicBeamAmp} 
gives the max-norm errors in the steady state solutions for the light, medium and heavy beams.  As expected, the AMP scheme is found to be stable and accurate in all three cases.  Finally, we note that for this problem the traditional partitioned scheme (without sub-time-step iterations) is found to be unstable when $\rhos\hs$ is less than $5$ (approximately).

{
\newcommand{\figWidth}{8cm}
\newcommand{\trimfig}[2]{\trimFig{#1}{#2}{.02}{0.02}{.045}{.005}}
\begin{figure}[htb]
\begin{center}
\begin{tikzpicture}[scale=1]
  \draw(4.0,5.0) node[anchor=south west,xshift=-4pt,yshift=+0pt] {\trimfig{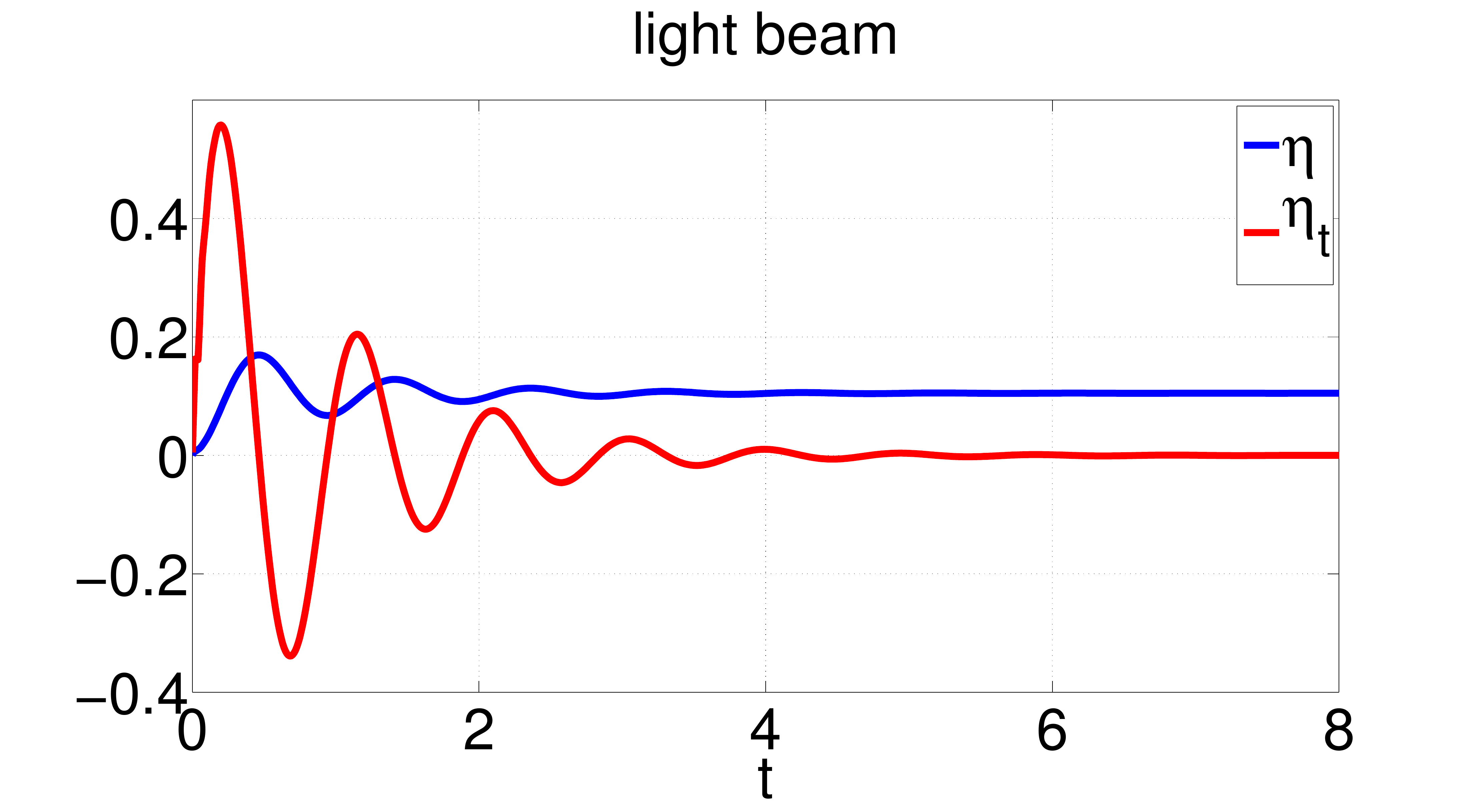}{\figWidth}};
    \draw(0.0,0.0) node[anchor=south west,xshift=-4pt,yshift=+0pt] {\trimfig{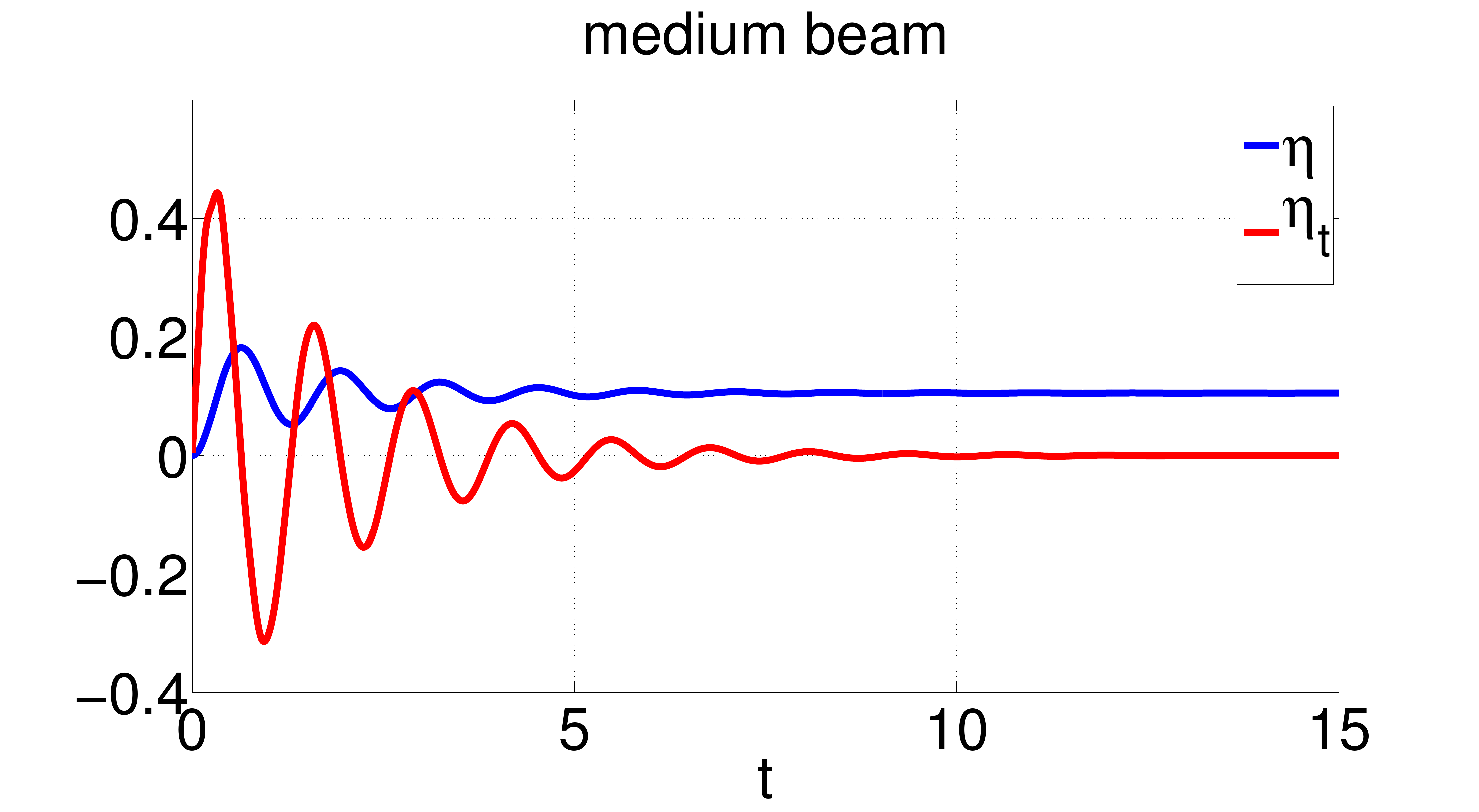}{\figWidth}};
   \draw(8.0,0.0) node[anchor=south west,xshift=-4pt,yshift=+0pt] {\trimfig{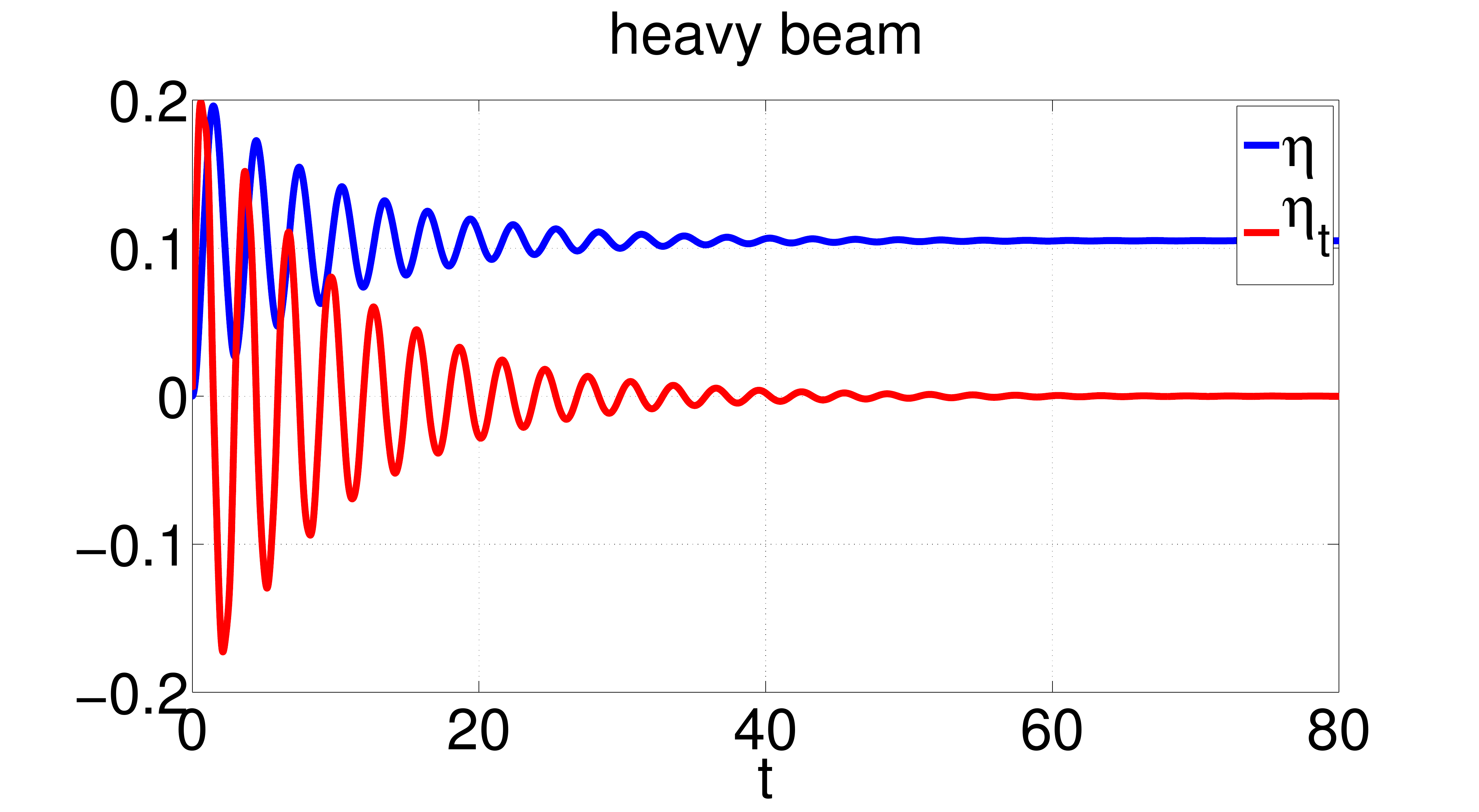}{\figWidth}};
\end{tikzpicture}
\end{center}
  \caption{Beam displacement $\eta$ and velocity $\eta_t$ at its center, $\xx=1/2$, versus time for light ($\rhos\hs=10^{-3}$), medium ($\rhos\hs=1$) and heavy ($\rhos\hs=10$) beams.}
  \label{fig:ParabolicBeam}
\end{figure}
}
\tableParabolicBeamAmp

\subsection{Verification results for a beam with finite thickness}

For the case of FSI problems involving beams with finite thickness, we first consider solutions of two further benchmark problems to verify the stability and accuracy of the AMP scheme for this more difficult configuration.  All numerical solutions described in this section, and in the subsequent sections, are computed using deforming composite grids and the beam solver based on finite elements.  As part of the verification tests, we consider both the AMP-PBA and AMP-PBF implementations of the AMP interface conditions.  For both benchmark problems, the beam has constant nonzero thickness $\hs$ and separates fluid domains on either side.

\subsubsection{Dynamic motion of a flat beam separating two fluid chambers}  \label{sec:flexibleChannelTD}

We consider a beam separating two initially rectangular fluid domains.  The beam is assumed to have finite thickness $\hs$, and
is chosen to initially lie along the $\xx$-axis and occupy the domain $\OmegaS(0)=[-1,1]\times[-\hs/2,\hs/2]$. 
At time $t=0$, the lower fluid domain is the rectangular domain $\OmegaF_-(0)  = [-1,1]\times[\Hf^{-},-\hs/2]$,
while the upper fluid domain is $\OmegaF_+(0)  = [-1,1]\times[\hs/2,\Hf^{+}]$.  For this first test problem, we assume that
the beam satisfies a {\em sliding} boundary condition $\eta_s=\Es\Is\eta_{sss}=0$
at each end, and that the fluid satisfies slip-wall conditions along the vertical boundaries.
At the top and bottom fluid boundaries, we assume that $v_1=0$ and $p=P_\pm(t)$, where $P_\pm(t)$ are specified time-dependent pressures.
For this configuration, an exact solution can be found for which the beam moves vertically, but remains flat.  In the fluid, the horizontal
component of the velocity is zero, while the vertical component of the velocity and the pressure depend only on $\yy$ and $t$.  This solution,
albeit relatively simple, provides
a good test of the stability of the AMP scheme for finite-thickness beams and for the two implementations of the AMP interface
conditions.

Working through the equations governing the fluids in each chamber, assuming no dependence on $\xx$ and that the densities in both chambers
are equal and given by $\rho$, and using the kinematic matching conditions at the top and bottom surfaces of the beam, we find that $v_1=0$ and $v_2=d\eta/dt$ in
both fluid chambers, and
\begin{equation}
p(\xv,t) = \left\{
\begin{array}{ll}
\displaystyle{
P_-(t)-\rho(\yy-\Hf^{-}){d\sp2\eta\over dt\sp2}
},\quad & \xv\in\OmegaF_-(t)=[-1,1]\times[\Hf^{-},\eta-\hs/2], \quad t>0,\medskip\\
\displaystyle{
P_+(t)-\rho(\yy-\Hf^{+}){d\sp2\eta\over dt\sp2}
}, & \xv\in\OmegaF_+(t)=[-1,1]\times[\eta+\hs/2,\Hf^{+}],\quad t>0.
\end{array}
\right.
\label{eq:chamberPressure}
\end{equation}
The equation governing an EB beam with all terms involving spatial derivatives set to zero is
\begin{equation}
   \rhos\hs{d\sp2\eta\over dt\sp2} = - \Ks_0 \eta - \Kt {d\eta\over dt} + p(\xx,\eta-\hs/2,t) - p(\xx,\eta+\hs/2,t) ,\qquad \xx\in[-1,1],\quad t>0,
\label{eq:beamOnlyTime}
\end{equation}
Using the expressions for the fluid pressures in~\eqref{eq:chamberPressure}, we find that the vertical displacement of the (flat) beam satisfies
\begin{equation}
(\rhos\hs+\rho\Hf){d\sp2\eta\over dt\sp2}+ \Kt {d\eta\over dt} + \Ks_0 \eta = \Delta P(t),\qquad t>0,
\label{eq:spring}
\end{equation}
where $\Hf=\Hf^{+}-\Hf^{-}-\hs$ is the total depth of the fluid chambers and $\Delta P(t)=P_-(t)-P_+(t)$ is the difference between
the applied pressures.  We observe that the displacement of the beam satisfies a standard spring equation with mass given by the
sum of the beam mass, $\rhos\hs$, and the fluid added mass, $\rho\Hf$, damping given by $\Kt$, spring constant given by $\Ks_0$, and a time-dependent external forcing given by $\Delta P(t)$.
Exact solutions for $\eta(t)$ satisfying~\eqref{eq:spring} for chosen initial conditions are easily found.  For example, if
$\eta=d\eta/dt=0$ at $t=0$, and if $\Kt=0$ and $\Delta P(t)$ is constant, then
\[
\eta(t) = \frac{\Delta P}{K_0} \bigl\{( 1 - \cos(\omega t) \bigr\},\qquad t>0,
\]
where the frequency of the oscillating beam is given by
\[
\omega=\sqrt{\frac{K_0}{\rhos\hs + \rho H}}.
\]

Numerical solutions of this FSI problem are computed using the full DCG methodology as discussed in Section~\ref{sec:dcg}. 
The composite grid $\Gcpc^{(j)}$ for the geometry is shown in the left plot of Figure~\ref{fig:partitionedChambersBeamGrids},
and consists of a background Cartesian grid and a body-fitted grid adjacent to the beam surface for each fluid
chamber.  The two body-fitted grids are generated independently using the hyperbolic grids generator, and each  
has $6$ grids lines in the normal direction (so that the grids becomes thinner as
the composite grid is refined).  The grid spacings in each coordinate direction are approximately equal for all component grids 
and chosen to be $h_j=1/(10 j)$.  As before, the index $j$ determines the resolution of the grid.



{
\newcommand{\figWidth}{8cm}
\newcommand{\trimfig}[2]{\trimFig{#1}{#2}{.0}{0.0}{.0}{.0}}
\begin{figure}[htb]
\begin{center}
\begin{tikzpicture}[scale=1]
  \useasboundingbox (0.0,.75) rectangle (17,8);  
  \draw(0.0,0.0) node[anchor=south west,xshift=-4pt,yshift=+0pt] {\trimfig{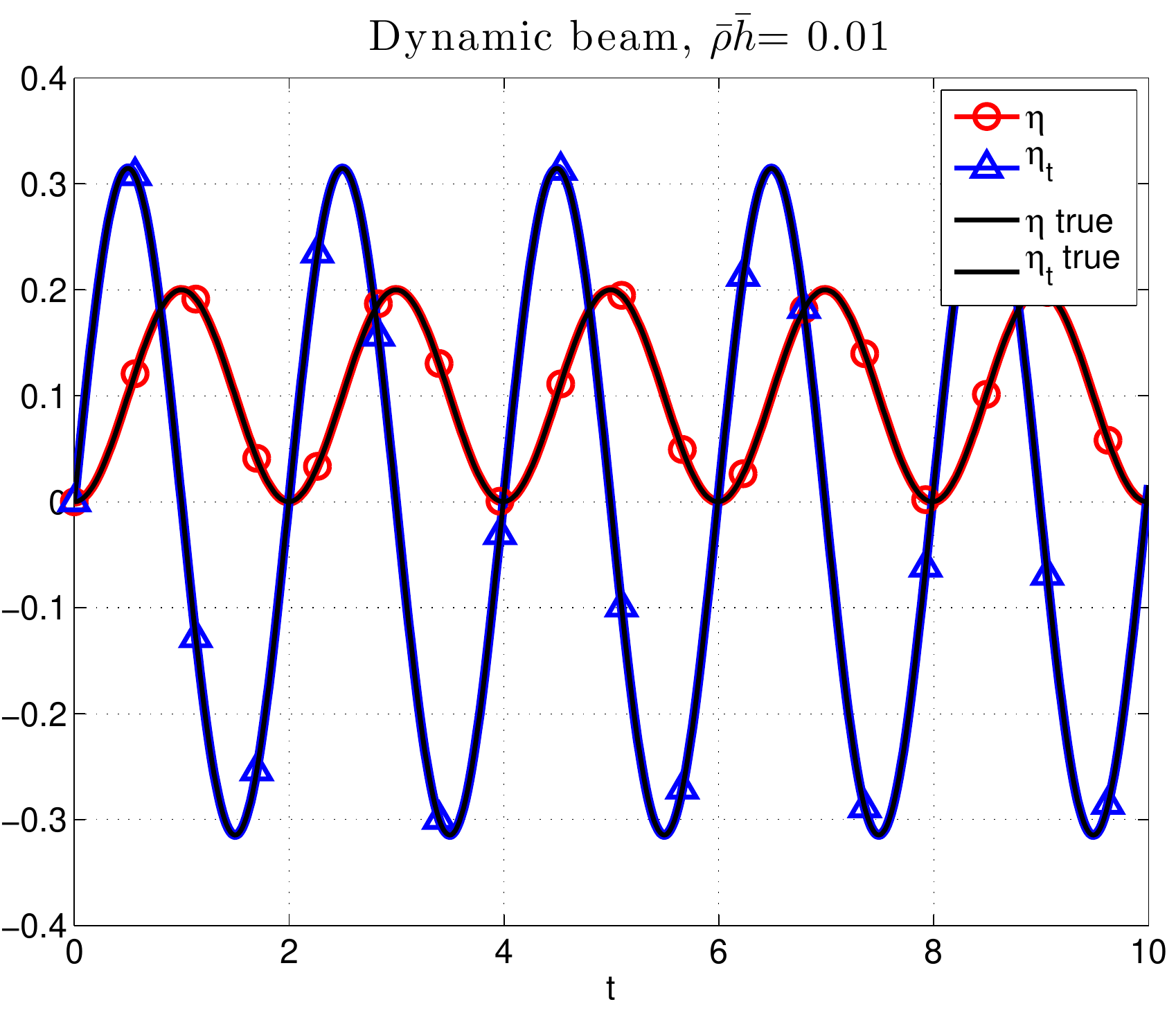}{\figWidth}};
  \draw(8.5,0.0) node[anchor=south west,xshift=-4pt,yshift=+0pt] {\trimfig{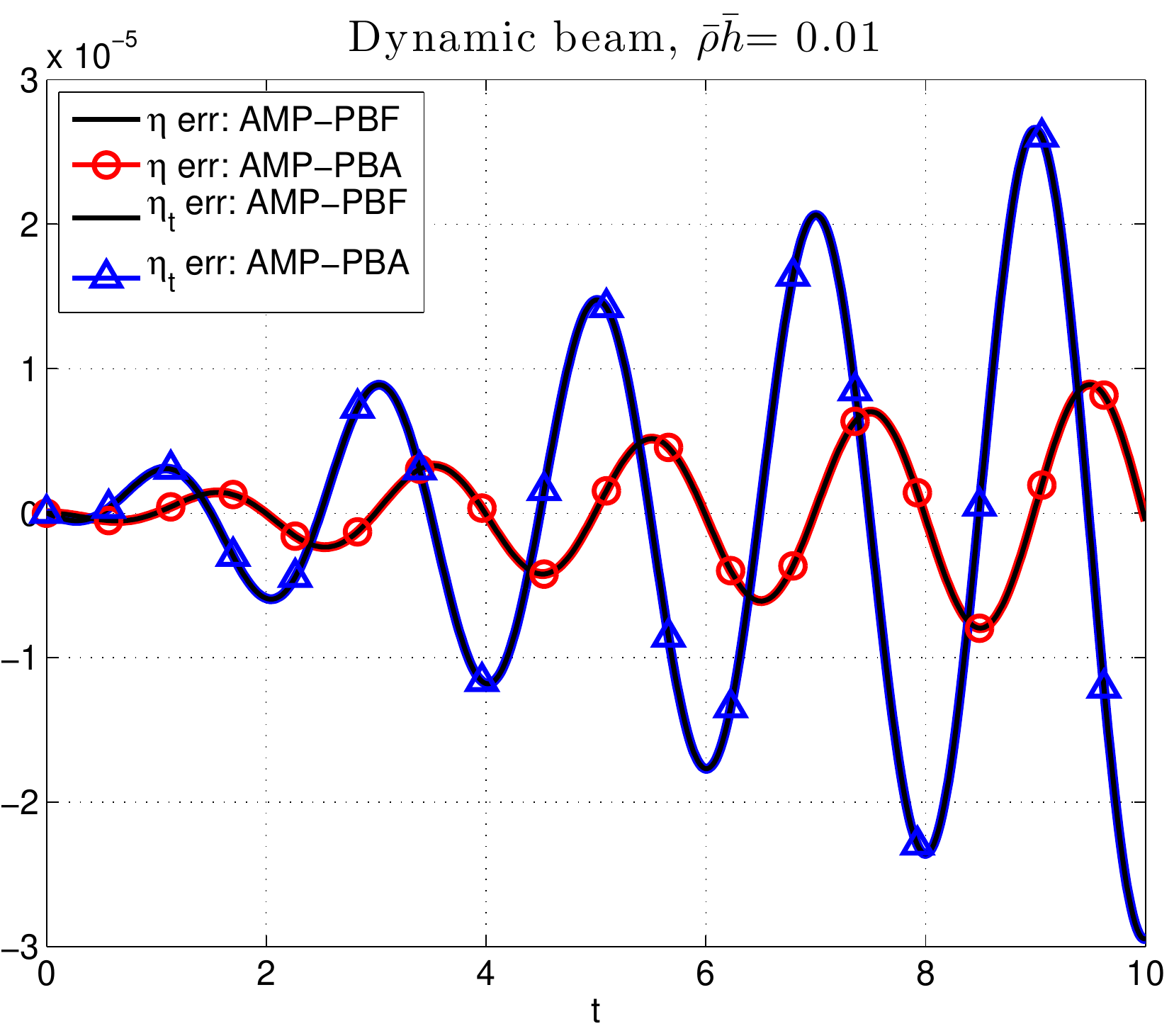}{\figWidth}};
\end{tikzpicture}
\end{center}
\caption{Oscillating beam between two fluids for $\rhos\hs=0.01$. 
Left: beam vertical displacement, $\eta$, and velocity, $\eta_t$,
compared to the exact solution. Right: errors in $\eta$ and $\eta_t$ for the AMP-PBF (predict beam internal force) and 
AMP-PBA (predict beam acceleration) schemes. Results are for the coarse grid $\Gcpc^{(2)}$.}\label{fig:oscillatingBeam}
\end{figure}
}

\bogus{
\begin{table}[hbt]\tableFont 
\begin{center}
\begin{tabular}{|l|c|c|c|c|c|c|c|c|c|c|c|} \hline 
grid &     N &    $p$   & r & $\vv$ & r & $\grad\cdot\uv$ & r &      $\eta$ & r &    $\eta_t$ & r \\ \hline 
  flexiblePartition1 &     1 & \num{5.5}{-6} &      & \num{1.2}{-5} &      & \num{1.0}{-5} &      & \num{7.9}{-6} &      & \num{9.6}{-6} &       \\ \hline
  flexiblePartition2 &     2 & \num{6.7}{-7} &  8.3 & \num{1.2}{-6} & 10.2 & \num{8.2}{-7} & 12.5 & \num{8.8}{-7} &  9.0 & \num{1.1}{-6} &  8.9  \\ \hline
  flexiblePartition4 &     4 & \num{6.2}{-8} & 10.8 & \num{1.1}{-7} & 10.5 & \num{4.9}{-8} & 16.7 & \num{8.8}{-8} & 10.0 & \num{1.1}{-7} & 10.0  \\ \hline
  flexiblePartition8 &     8 & \num{4.0}{-9} & 15.2 & \num{7.5}{-9} & 14.6 & \num{1.8}{-9} & 28.0 & \num{6.1}{-9} & 14.4 & \num{7.5}{-9} & 14.5  \\ \hline
    rate             &       &  $3.47$       &      &  $3.52$       &      &  $4.16$       &      &  $3.43$       &      &  $3.43$       &       \\ \hline
\end{tabular}
\caption{Cgins, beamPiston, ins, exact.beamPistonRhos100p0.pc.flexiblePartition.ins.order2.vector, ts=pc, $t=.7$, $\nu=.02$, dtMax=0.1, aftol=1e-10, axialAxis=2, kp=1, , cfl=0.9, -rhoBeam=100.0 -E=1.e-3 -tension=1.e-1 -thickness=.1 -K0=10 -p0=1 -addedMass=1 -fluidOnTwoSides=1  -ampProjectVelocity=1 -projectNormalComponent=1 -nis=4 -option=beamPiston -sideBC=slipWall -cfls=1., Thu Jun 18 15:26:22 2015}\label{table:exact.beamPistonRhos100p0.pc.flexiblePartition.ins.order2.vector}
\end{center}
\end{table}
}
\newcommand{\tableBeamPistonHeavy}{%
\begin{tabular}{|c|c|c|c|c|c|c|c|c|} \hline
\multicolumn{9}{|c|}{\strutt Oscillating beam, $\rhos\hs=10$} \\ \hline 
\strutt~~$h_j$~~& $E_j^{(p)}$ & $r$ & $E_j^{(\vv)}$  & $r$ & $E_j^{(\eta)}$ &  $r$ & $E_j^{(\eta_t)}$ &  $r$\\ \hline 
 1/20      &  \num{5.5}{-6} &      & \num{1.2}{-5} &      &  \num{7.9}{-6} &      & \num{9.6}{-6} &       \\ \hline
 1/40      &  \num{6.7}{-7} &  8.3 & \num{1.2}{-6} & 10.2 &  \num{8.8}{-7} &  9.0 & \num{1.1}{-6} &  8.9  \\ \hline
 1/80      &  \num{6.2}{-8} & 10.8 & \num{1.1}{-7} & 10.5 &  \num{8.8}{-8} & 10.0 & \num{1.1}{-7} & 10.0  \\ \hline
 1/160     &  \num{4.0}{-9} & 15.2 & \num{7.5}{-9} & 14.6 &  \num{6.1}{-9} & 14.4 & \num{7.5}{-9} & 14.5  \\ \hline
\rateLabel &   $3.47$       &      &  $3.52$       &      &   $3.43$       &      &  $3.43$       &       \\ \hline
\end{tabular}
}

\bogus{
\begin{table}[hbt]\tableFont 
\begin{center}
\begin{tabular}{|l|c|c|c|c|c|c|c|c|c|c|c|} \hline 
grid &     N &    $p$   & r & $\vv$ & r & $\grad\cdot\uv$ & r &      $\eta$ & r &    $\eta_t$ & r \\ \hline 
  flexiblePartition1 &     1 & \num{3.1}{-4} &      & \num{1.7}{-4} &      & \num{4.2}{-5} &      & \num{4.8}{-5} &      & \num{9.1}{-5} &       \\ \hline
  flexiblePartition2 &     2 & \num{1.3}{-5} & 22.8 & \num{5.1}{-6} & 34.1 & \num{8.4}{-7} & 49.7 & \num{2.1}{-6} & 22.8 & \num{4.8}{-6} & 19.0  \\ \hline
  flexiblePartition4 &     4 & \num{9.2}{-7} & 14.6 & \num{3.3}{-7} & 15.3 & \num{2.7}{-8} & 30.8 & \num{1.4}{-7} & 14.6 & \num{3.3}{-7} & 14.5  \\ \hline
  flexiblePartition8 &     8 & \num{5.8}{-8} & 15.9 & \num{2.1}{-8} & 16.1 & \num{1.3}{-9} & 20.6 & \num{8.9}{-9} & 16.0 & \num{2.1}{-8} & 15.9  \\ \hline
    rate             &       &  $4.10$       &      &  $4.30$       &      &  $4.98$       &      &  $4.10$       &      &  $4.02$       &       \\ \hline
\end{tabular}
\caption{Cgins, beamPiston, ins, exact.beamPistonRhosp1.pc.flexiblePartition.ins.order2.vector, ts=pc, $t=.7$, $\nu=.02$, dtMax=0.1, aftol=1e-10, axialAxis=2, kp=1, , cfl=0.9, -rhoBeam=.1 -E=1.e-3 -tension=1.e-1 -thickness=.1 -K0=10 -p0=1 -addedMass=1 -fluidOnTwoSides=1  -ampProjectVelocity=1 -projectNormalComponent=1 -nis=4 -option=beamPiston -sideBC=slipWall -cfls=2., Thu Jun 18 15:08:39 2015}\label{table:exact.beamPistonRhosp1.pc.flexiblePartition.ins.order2.vector}
\end{center}
\end{table}
}
\newcommand{\tableBeamPistonLight}{%
\begin{tabular}{|c|c|c|c|c|c|c|c|c|} \hline
\multicolumn{9}{|c|}{\strutt Oscillating beam, $\rhos\hs=0.01$} \\ \hline 
\strutt~~$h_j$~~& $E_j^{(p)}$ & $r$ & $E_j^{(\vv)}$  & $r$ & $E_j^{(\eta)}$ &  $r$ & $E_j^{(\eta_t)}$ &  $r$\\ \hline 
 1/20      & \num{3.1}{-4} &      & \num{1.7}{-4} &      &  \num{4.8}{-5} &      & \num{9.1}{-5} &       \\ \hline
 1/40      & \num{1.3}{-5} & 22.8 & \num{5.1}{-6} & 34.1 &  \num{2.1}{-6} & 22.8 & \num{4.8}{-6} & 19.0  \\ \hline
 1/80      & \num{9.2}{-7} & 14.6 & \num{3.3}{-7} & 15.3 &  \num{1.4}{-7} & 14.6 & \num{3.3}{-7} & 14.5  \\ \hline
 1/160     & \num{5.8}{-8} & 15.9 & \num{2.1}{-8} & 16.1 &  \num{8.9}{-9} & 16.0 & \num{2.1}{-8} & 15.9  \\ \hline
\rateLabel &  $4.10$       &      &  $4.30$       &      &   $4.10$       &      &  $4.02$       &       \\ \hline
\end{tabular}
}

\begin{figure}[hbt]\tableFontSize
\begin{center}
\tableBeamPistonHeavy  \\
\vskip\baselineskip
\tableBeamPistonLight
\caption{Oscillating beam between two fluids. Maximum errors and 
estimated convergence rates at $t=0.7$ computed using the AMP scheme for a heavy beam, $\rhos\hs=10$ and
a light beam $\rhos\hs=0.01$. For this problem the convergence rates are close to four since the errors
are primarily due to the beam solver which is fourth-order accurate in space. 
}
\label{tab:beamPiston}
\end{center}
\end{figure}

The parameters for the EB beam in the numerical calculations are taken 
as $\hs=0.1$, $\Ks_0=10$, $\Ts=10^{-1}$, $\Es\Is=10^{-3}$, $\Kt=\Kxxt=0$, and $\rhos$ varied to test light and heavy beams.
The parameters for the fluid are $\rho=1$, $\mu=.02$, $\Hf\sp-=-.5-\hs/2$, and $\Hf\sp+=.5+\hs/2$ so that $H=1$.
The pressure difference between the top and bottom was taken as $\Delta P=1$.
Figure~\ref{fig:oscillatingBeam} shows the computed values of $\eta$ and $d\eta/dt$ compared to
the exact solution for a light beam $\rhos\hs=10^{-2}$ on the grid $\Gcpc^{(2)}$.
The differences between the computed and exact solutions are nearly indistinguishable.
Figure~\ref{fig:oscillatingBeam} also plots the errors in $\eta$ and $d\eta/dt$ 
for the AMP-PBF (predict beam internal force) and 
AMP-PBA (predict beam acceleration) schemes, and the results are again nearly indistinguishable. 
The tables in Figure~\ref{tab:beamPiston} shows the max-norm errors and estimated convergence rates for a
light beam, $\rhos\hs=10^{-2}$, and a heavy beam, $\rhos\hs=10$ at $t=0.7$.  
The maximum errors are converging close to fourth-order accuracy. Recall that the beam solver is 
fourth-order accurate in space and second-order order accurate in time but converges as fourth-order overall
when the time-step scales as $\dt \propto h_j^2$, which is the case for the calculations presented here.



\subsubsection{Steady state beam separating two pressurized fluid chambers}

As a final benchmark test case for which an exact solution is available, we consider the long-time behavior of an EB beam
separating two initially rectangular pressurized fluid chambers.  The geometric configuration for this test case is the
same as that considered in the previous test for a beam with finite thickness $\hs$, but the boundary conditions on the
vertical sides are changed so that the beam does not remain flat as time evolves.  Here we assume no-slip boundary conditions for the fluid
on the side walls, and take $\eta=\partial\eta/\partial s=0$ as {\em clamped} boundary conditions for the EB beam.  The material properties
of the beam are taken to be $\hs=0.1$, $\Ts=0$, $\Es\Is=0.2$, $\Kt=5$ and $\Kxxt=0$, with $\rhos$ varied to consider light,
medium and heavy beams.  The properties of the fluid are $\rho=1$, $\mu=.05$, $P_+=0$ and $P_-(t)=P_0R(t)$, with $P_0=1$,
where $R(t)$ is a smooth {\em ramp} function given by
\begin{align}
    R(t) =          \begin{cases}
            (35+(-84+(70- 20 t) t) t) t^4,\quad     & \text{for $0\le t \le 1$}, \\
                1,   & \text{for $t>1$}.
               \end{cases}   \label{eq:ramp}
\end{align}
The ramp function satisfies $R=R\sp\prime=R\sp{\prime\prime}=R\sp{\prime\prime\prime}=0$ at $t=0$ and has three continuous derivatives at $t=1$. 
 Thus, the applied pressure in the lower chamber varies smoothly from $P_-(0)=0$ to $P_-(1)=P_0$, and remains equal to $P_0$ for $t>1$.  At $t=0$, the fluid in both chambers is at rest, and the displacement and velocity of the beam are both zero.

As in the FSI problem considered in Section~\ref{sec:flexibleChannelSS}, the beam oscillates in response to the nonzero applied fluid pressure in the bottom chamber, but then
approaches a steady state solution as time increases due to the viscous dissipation in the fluid
and damping in the beam.  At steady state, the
velocity in both fluid chambers is zero and the fluid pressures are uniform as before, but now for the beam parameters considered here the displacement of the beam is given by
\[
\eta  = \frac{P_0}{24 \Es\Is} \bigl(1-\xx^2\bigr)^2,\qquad \xx\in[-1,1].
\]

{
\newcommand{\figWidth}{8cm}
\newcommand{\trimfig}[2]{\trimFig{#1}{#2}{.02}{0.02}{.285}{.285}}
\begin{figure}[htb]
\begin{center}
\begin{tikzpicture}[scale=1]
  \useasboundingbox (0.0,.75) rectangle (16.,5.1);  
  \draw(0.0,0.0) node[anchor=south west,xshift=-4pt,yshift=+0pt] {\trimfig{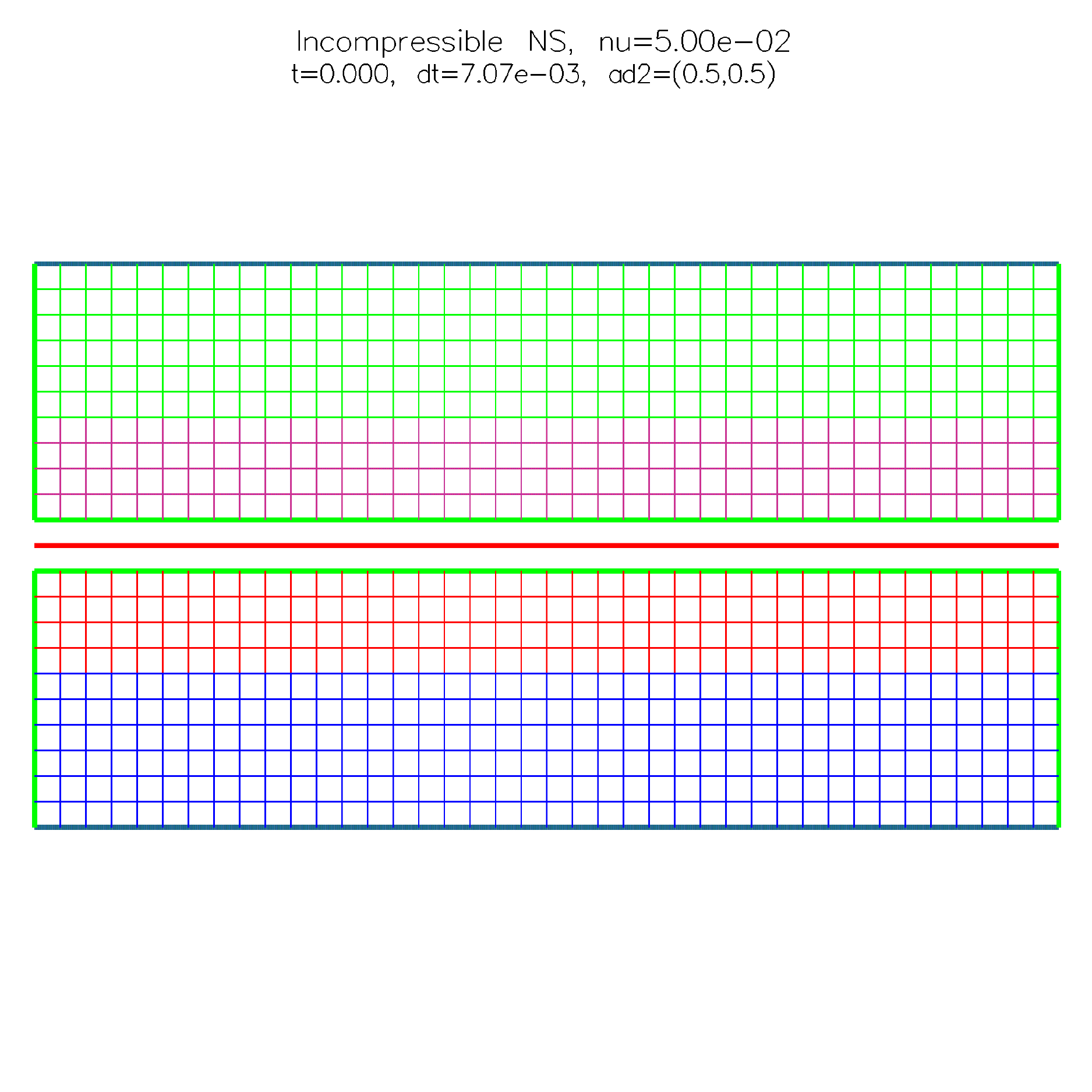}{\figWidth}};
  \draw(8.0,0.0) node[anchor=south west,xshift=-4pt,yshift=+0pt] {\trimfig{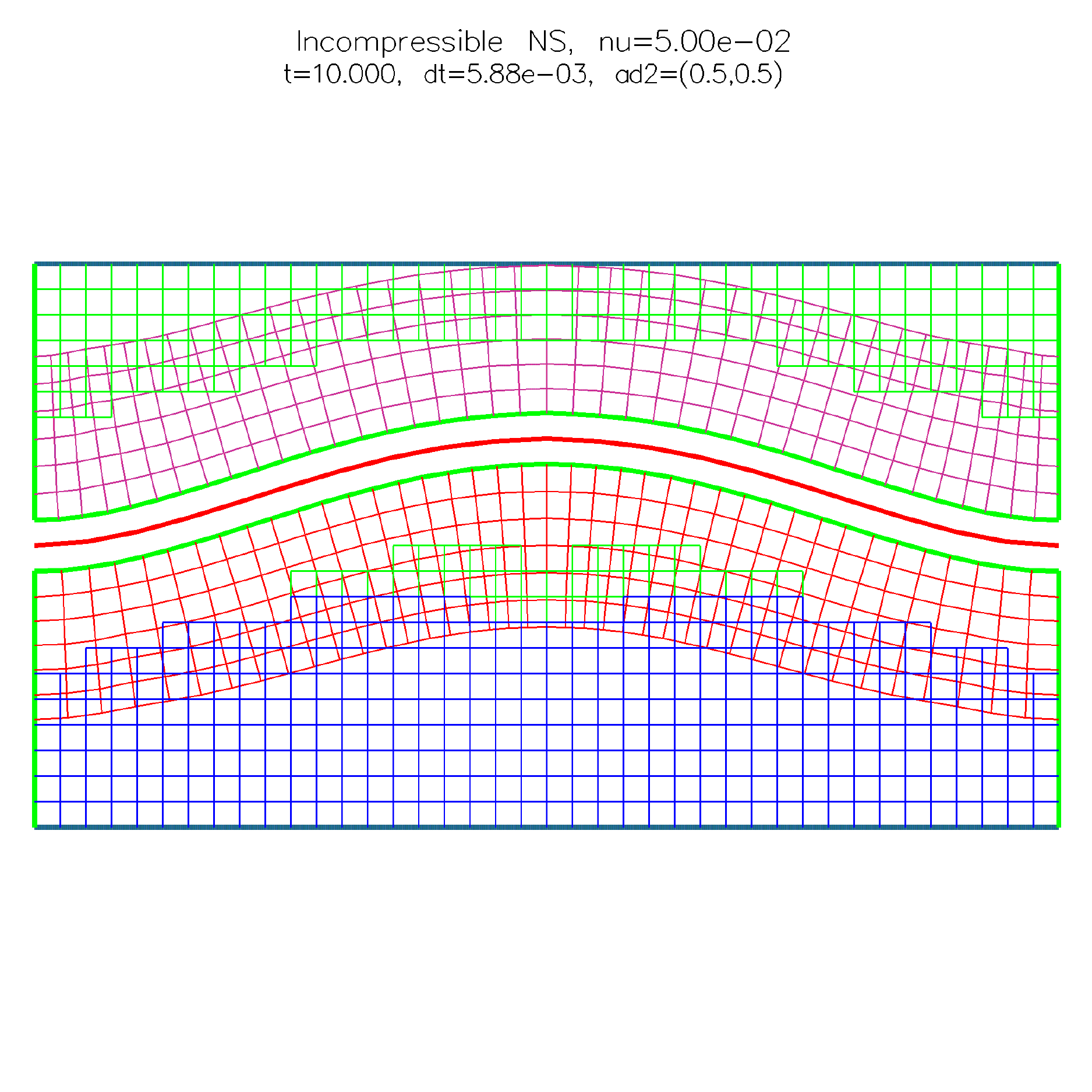}{\figWidth}};
\end{tikzpicture}
\end{center}
  \caption{Composite grids for the partitioned chambers: coarse grid $\Gcpc^{(2)}$ at $t=0$ (left) and at $t=10$ (right). 
    The solution at $t=10$ has nearly reached steady state.}
  \label{fig:partitionedChambersBeamGrids}
\end{figure}
}

Numerical solutions are computed using the composite grid $\Gcpc^{(j)}$, as before, with grid spacings approximately equal to $h_j=1/(10 j)$.  Figure~\ref{fig:partitionedChambersBeamGrids} shows a representative grid, $\Gcpc^{(2)}$, at $t=0$ (left) and at $t=10$ (right) computed for the case $\rhos\hs=0.1$.  The solution at $t=10$ has nearly reached steady state.  The maximum error in the beam displacement compared to the steady state solution is approximately $4\times 10^{-7}$  at $t=10$,
and reaches $10^{-15}$ by $t=50$.
Figure~\ref{fig:beamUnderPressureTimeHistory} shows the motion of the beam midpoint for different beam
densities and different algorithmic options with a smaller damping coefficient of $\Kt=1$.
The top-left graph in the figure shows the behavior of the
displacement, $\eta(0,t)$, and
velocity, $\eta_t(0,t)$, as functions of time. Results are shown for
a light beam, $\rhos\hs=0.01$, and a medium beam, $\rhos\hs=1$, computed on a coarse grid $\Gcpc^{(2)}$ and
a finer grid $\Gcpc^{(4)}$.  The light beam is excited more rapidly than the medium, and undergoes more rapid oscillations.
We also observe that the time histories computed using the coarse and fine grid are nearly indistinguishable for both beams.
The top-right graph in Figure~\ref{fig:beamUnderPressureTimeHistory} 
compares results from the AMP scheme and the traditional partitioned scheme
with sub-iterations (TP-SI) for the case $\rhos\hs=0.01$. Added-mass effects are relatively large for this light beam requiring a relaxation parameter
in the TP-SI iteration of $\omega=0.1$. The TP-SI scheme required an average of approximately $77$ sub-iterations
per time step to achieve a sub-iteration convergence
tolerance of $10^{-5}$. The results from the AMP and TP-SI schemes are in good agreement, but the computational cost
for the TP-SI scheme is significantly larger.

The dynamic beam problem provides a good test to compare the behavior of the two AMP variations described in
Section~\ref{sec:beamSolverNumerical} for the finite-element based beam solver. 
It is our experience for this problem, and others, 
that for {\em very light} beams the AMP-PBA (predict beam acceleration)
algorithm tends to be more robust (e.g. requires fewer or no smoothing iterations)
than the AMP-PBF (predict beam internal force). The AMP-PBF
approach requires the evaluation of the third spatial-derivative of the finite element beam displacement and this can be sensitive, especially near the ends of the beam where one-sided approximations are used (as might be expected).
The bottom graph in Figure~\ref{fig:beamUnderPressureTimeHistory} compares 
results obtained using the AMP-PBA and AMP-PBF for a moderately light beam with $\rhos\hs=0.1$. In this case
the results obtained using both AMP variations are in excellent agreement. However, for a very light beam, $\rhos\hs=0.01$, the AMP-PBF variation has some difficulty 
and shows poor behavior near the ends of the beam where the acceleration of the beam varies very rapidly during
the startup phase when the pressure forcing is turned on. 
It is entirely possible that this sub-optimal behavior of the AMP-PBF variation
could be rectified with further refinements
to the algorithm, but we have not pursued this since the alternative AMP-PBA variation behaves well and is  
simpler to implement.

{
\newcommand{\figWidth}{8cm}
\newcommand{\trimfig}[2]{\trimFig{#1}{#2}{.0}{0.0}{.0}{.0}}
\begin{figure}[htb]
\begin{center}
\begin{tikzpicture}[scale=1]
  \useasboundingbox (0.0,-6.25) rectangle (17,8);  
  \draw(0.0,0.0) node[anchor=south west,xshift=-4pt,yshift=+0pt] {\trimfig{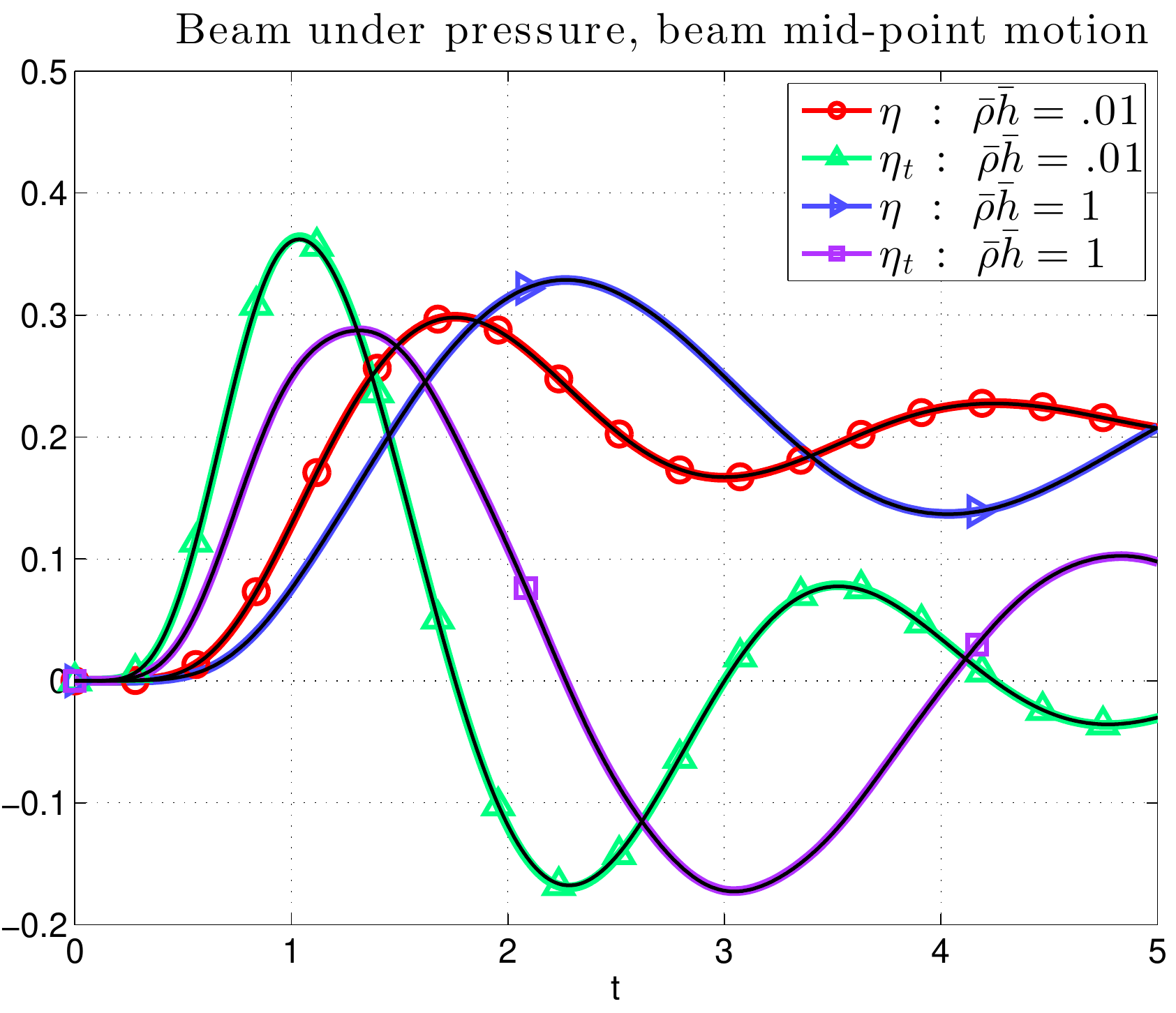}{\figWidth}};
  \draw(8.5,0.0) node[anchor=south west,xshift=-4pt,yshift=+0pt] {\trimfig{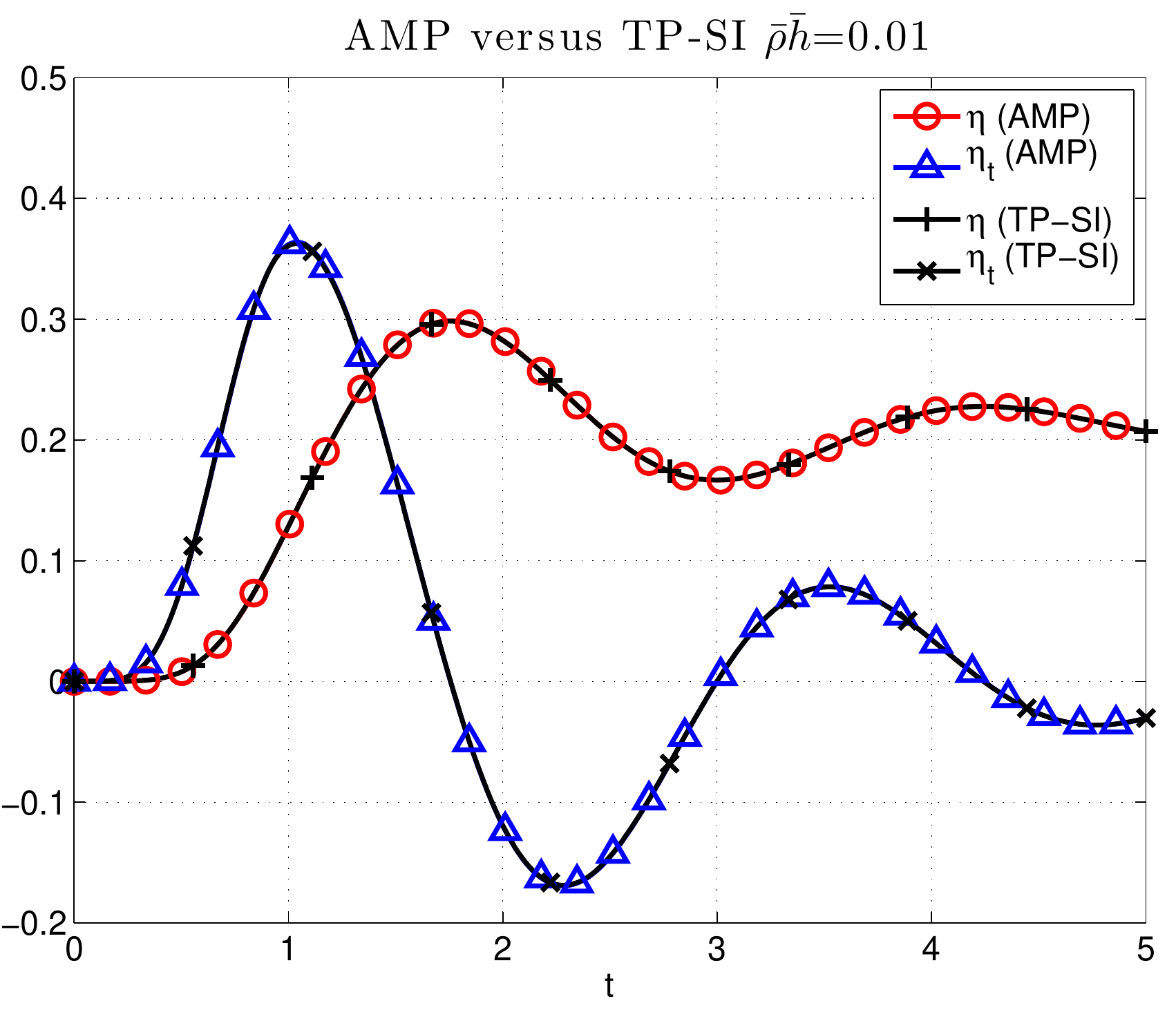}{\figWidth}};
  \draw(4.25,-7.) node[anchor=south west,xshift=-4pt,yshift=+0pt] {\trimfig{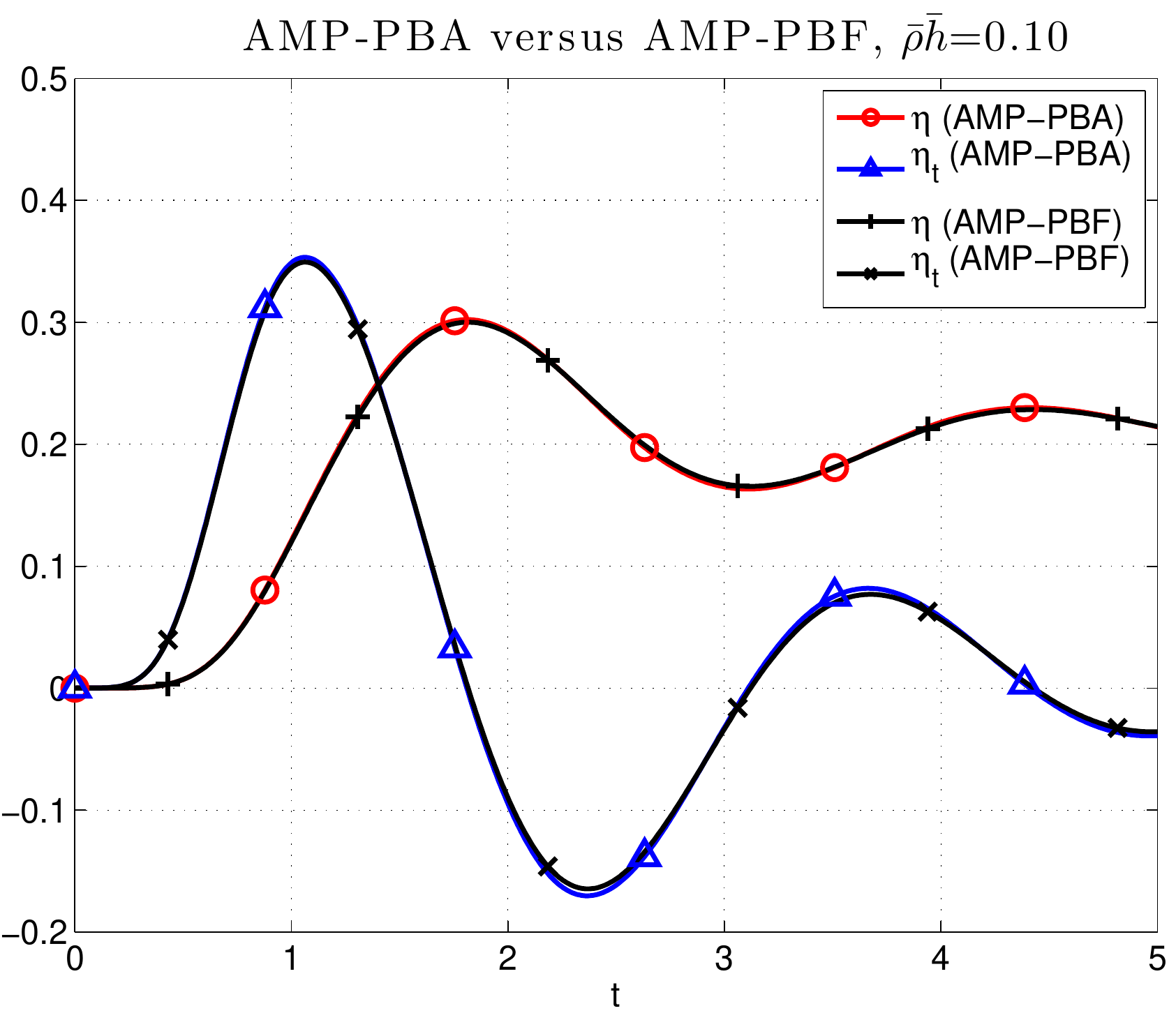}{\figWidth}};
\end{tikzpicture}
\end{center}
\caption{Beam under pressure, motion of the beam mid-point. Top left: $\eta$ and $\eta_t$ for a light ($\rhos\hs=0.01$) 
and medium ($\rhos\hs=1$) beam using a coarse grid $\Gcpc^{(2)}$ (colored curves) and a finer grid $\Gcpc^{(4)}$ (black curves).
Top right: AMP versus the TP-SI scheme for the light beam, $\rhos\hs=0.01$.  Bottom: A comparsion of the two AMP schemes for coupling
to the finite-element beam solver on grid $\Gcpc^{(4)}$ for $\rho\hs=0.1$.}
\label{fig:beamUnderPressureTimeHistory}
\end{figure}
}

\newlength{\ycbTop}%
\newlength{\ycbMid}%
\subsection{Traveling-wave pressure-pulse in an elastic tube} \label{sec:pressurePulse}

The propagation of a pressure pulse through
a two-dimensional elastic tube can be used as a model for the flow of blood
in a large artery or vein following the problem proposed in~\cite{FormaggiaGerbeauNobileQuarteroni2001},
and later considered in~\cite{GuidoboniGlowinskiCavalliniCanic2009,FernandezLandajuela2014}
to study the stability of partitioned schemes. 
Added-mass effects are important in this application
since the density of the elastic wall is close to that of blood
and also since the tube is long and thin. Traditional partitioned methods
require many sub-iterations per time step to remain stable for this problem. The AMP algorithm is shown to provide stable
results with no sub-iterations per time step. 

{
\newcommand{\figWidth}{16cm}
\newcommand{\trimfig}[2]{\trimFig{#1}{#2}{.0}{0.}{.31}{.32}}
\begin{figure}[htb]
\begin{center}
\begin{tikzpicture}[scale=1]
  \useasboundingbox (0.0,.75) rectangle (16.,7.1);  
  \draw(0.0,4.4) node[anchor=south west,xshift=-4pt,yshift=+0pt] {\trimfig{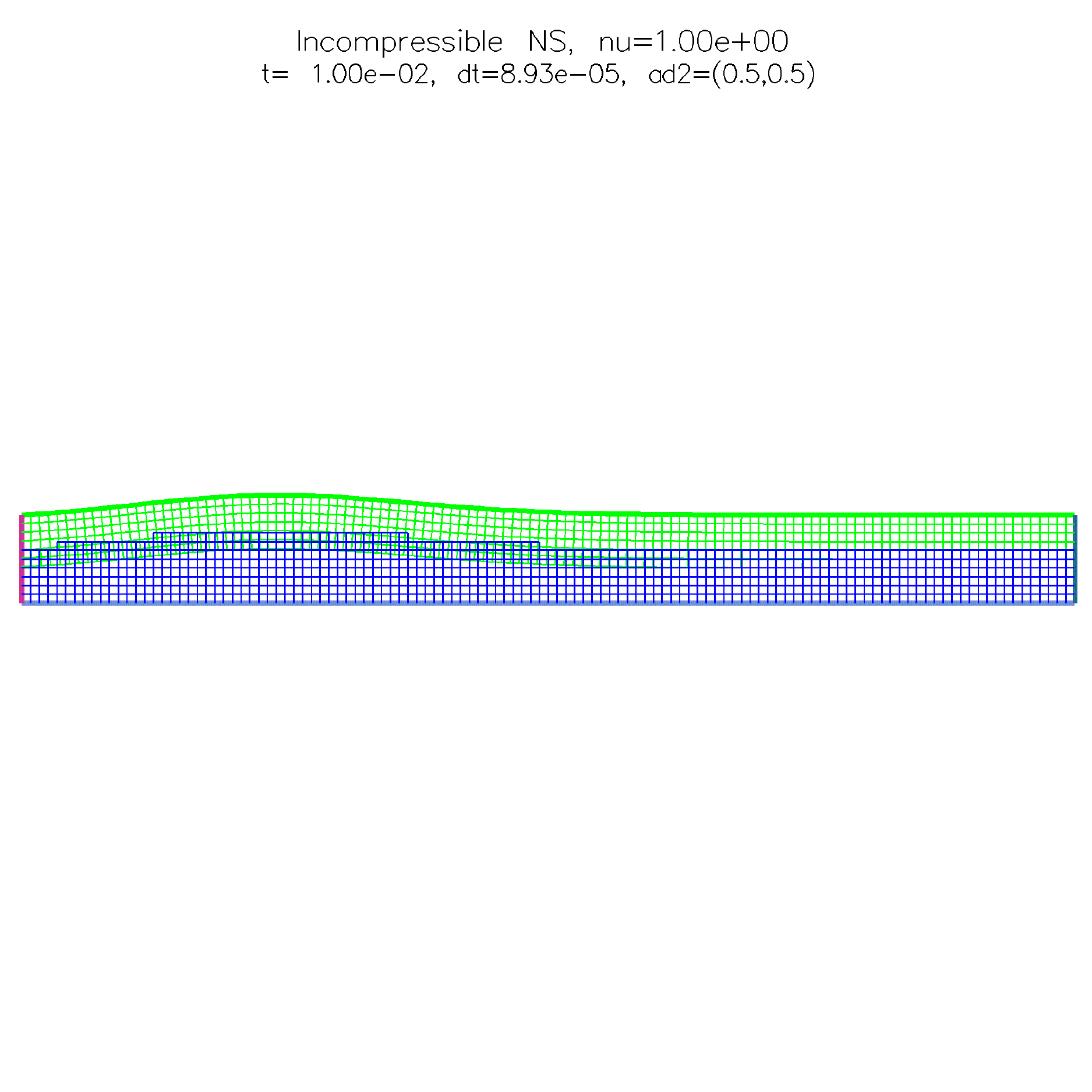}{\figWidth}};
  \draw(0.0,2.2) node[anchor=south west,xshift=-4pt,yshift=+0pt] {\trimfig{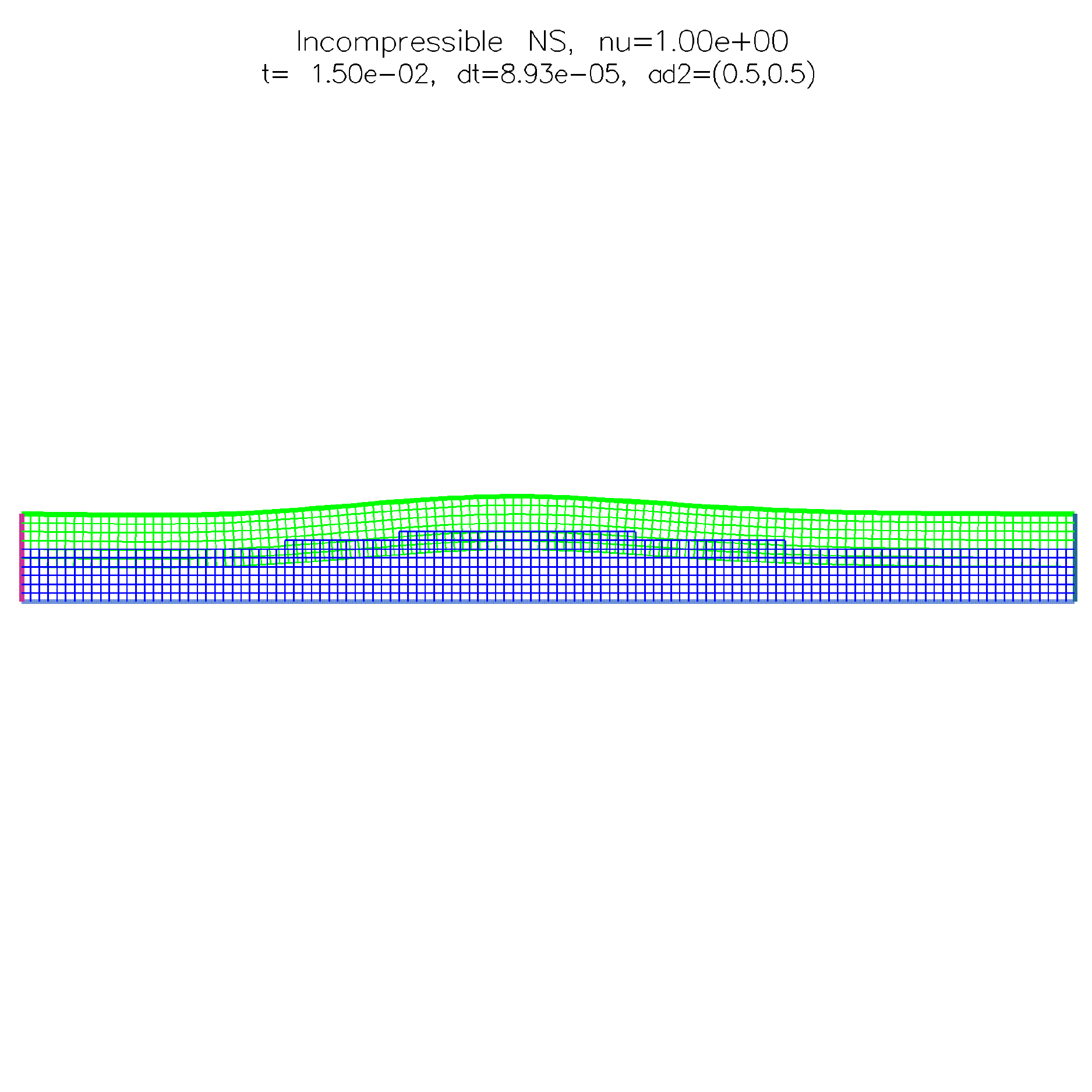}{\figWidth}};
  \draw(0.0,0.0) node[anchor=south west,xshift=-4pt,yshift=+0pt] {\trimfig{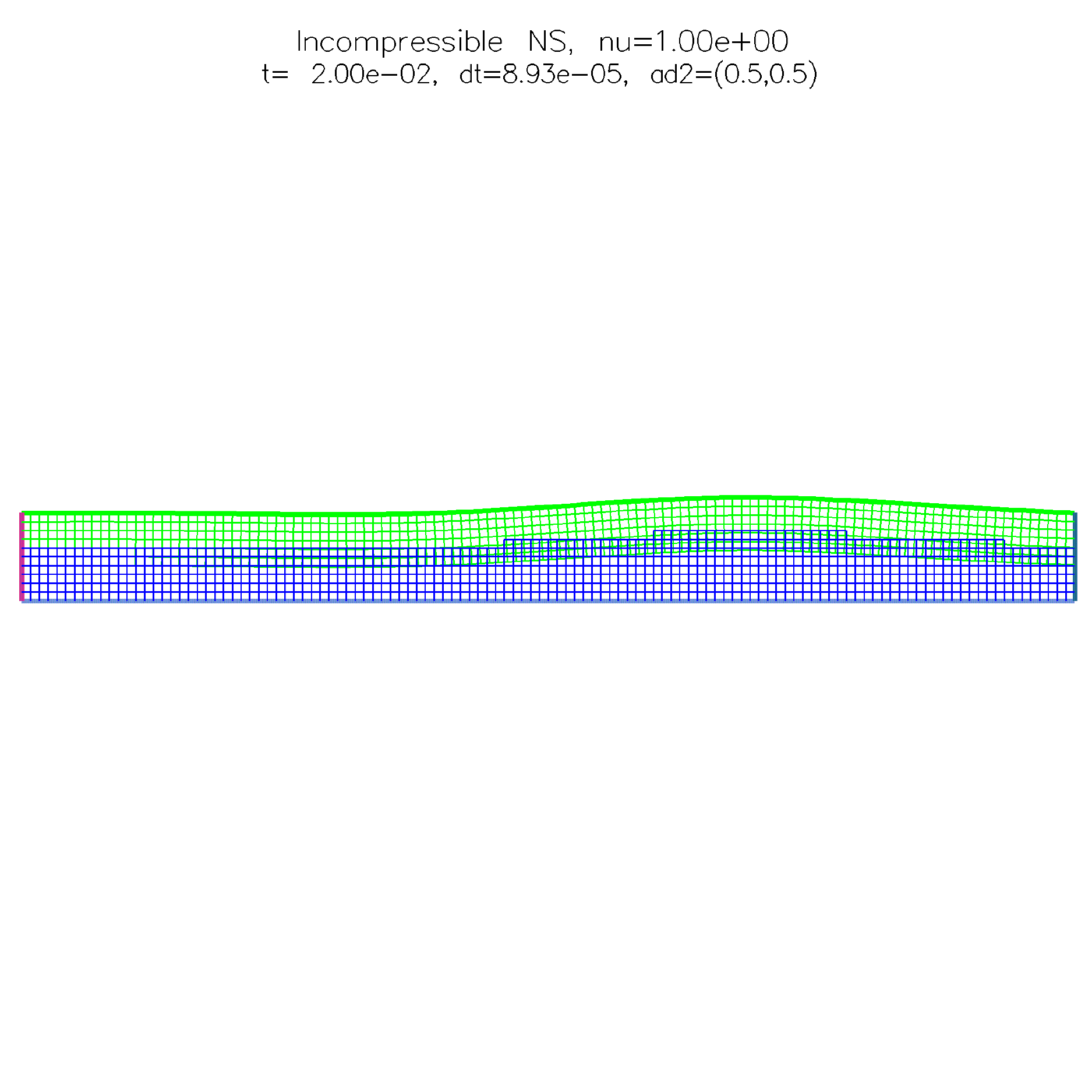}{\figWidth}};
\end{tikzpicture}
\end{center}
  \caption{Composite grids for the traveling-wave pressure-pulse in a flexible tube: coarse grid $\Gcfc^{(2)}$ at times 
   $t=1.0$, $1.5$ and $2.0$.}
  \label{fig:flexibleChannelGrids}
\end{figure}
}

The domain for the fluid is the channel $\OmegaF(t)=[0,L]\times[0,R(t)]$, where $L$ is the length
of the channel and $R(t)$ is the distance from the axis of symmetry of the channel at $\yy=0$ and the elastic wall.
We consider the case when $L=6$ and $R(0)=0.5$.  The boundary conditions for the fluid domain are as follows. 
The bottom boundary is a slip-wall where $v_2=0$ and $\partial v_1/\partial \yy=0$.
The top boundary is the flexible beam where the AMP interface
conditions are applied. On the left boundary the tangential component of the velocity is 
set to zero, $v_2=0$, and the pressure is specified to be 
the time-dependent pulse 
\begin{align}
    p(0,\yy,t) = \begin{cases}
           \pMax \, \sin( \pi t /\tMax)\quad & \text{for $0\le t \le \tMax$}, \\
           0                            & \text{for $t > \tMax$},
               \end{cases}    \label{eq:pressureBC}
\end{align}
where $\pMax$ and $\tMax$ determine the magnitude and duration of the pulse (values given below).
The right boundary is an outflow boundary where $p=0$ and the components of velocity are extrapolated.

The composite grid $\Gcfc^{(j)}$ for the fluid domain at various times is shown in Figure~\ref{fig:flexibleChannelGrids}.
The composite grid consists of a background Cartesian grid (shown in blue) together with a hyperbolic grid adjacent to the top
interface (shown in green). The hyperbolic grid is generated with $6$ grids lines in the normal direction and thus this boundary-fitted component grid becomes thinner as
the composite grid is refined. The grid spacing is chosen to be $h_j=1/(10 j)$ approximately.
%
%

We perform simulations of the elastic tube problem for two choices of the parameters of the fluid and the beam,
identified as Parameter Sets I and II.
The parameter sets are based on the parameters chosen in Fernandez et~al.~\cite{FernandezLandajuela2014}, but
with some differences as described below.  For both choices, the beam is taken to be a generalized string model
with $\Es\Is=0$ (i.e.~no fourth-order spatial derivative term in the beam equation).
The AMP-PBA (predict-beam-acceleration) variation is used for all calculations in this section
and the projected interface
velocity is smoothed with three iterations of the fourth-order filter, as described in Section~\ref{sec:velocityFilter}. 

\subsubsection{Parameter Set I}

{
\newcommand{\drawFC}[7]{%
\begin{scope}[#1]
\draw(0.0,0) node[anchor=south west,xshift=-4pt,yshift=+0pt] {\trimfig{flexibleChannel/#2}{\figWidth}};
\draw(.5,2.35) node[draw,fill=white,anchor=west,xshift=2pt,yshift=-16pt] {\scriptsize #3};
\draw(1.2,2.35) node[draw,fill=white,anchor=west,xshift=0pt,yshift=-16pt] {\scriptsize #5};
\draw (\xcb,\ycb) node[anchor=south west,xshift= +0pt,yshift=+0pt] {\trimfigcb{fig/colourBarLines}{\cbWidth}{\cbHeight}};
\draw (\xcb,\ycb) node[anchor=south west,xshift= +8pt,yshift=+1pt] {\scriptsize $#6$};
\draw (\xcb,\ycbTop) node[anchor=south west,xshift= +8pt,yshift=-6pt] {\scriptsize $#7$};
\end{scope}
}
{
\newcommand{\cbWidth}{.2cm}
\newcommand{\cbHeight}{2.8cm}
\newcommand{\xcb}{14.3cm}
\newcommand{\ycb}{-.15cm}
\setlength{\ycbTop}{\ycb+\cbHeight}
\setlength{\ycbMid}{\ycb+\cbHeight*\real{.5}}
\newcommand{\xLabel}{6.5cm}
\newcommand{\yLabel}{6.5cm}
\newcommand{\trimfigcb}[3]{\includegraphics[width=#2, height=#3, clip, trim=17cm 2.35cm 1.65cm 2.35cm]{#1}}
\newcommand{\figWidth}{14.5cm}
\newcommand{\trimfig}[2]{\trimFig{#1}{#2}{.03}{.15}{.33}{.34}}
\begin{figure}[htb]
\begin{center}
\begin{tikzpicture}[scale=1]
  \useasboundingbox (0.0,.5) rectangle (15.25,9);  
%
  \drawFC{xshift= 0.0cm,yshift= 6.0cm}{fc8b_p_t1p5}{$p$}{$p$}{$t=1.5$}{-.106}{1.40};
  \drawFC{xshift= 0.0cm,yshift= 3.0cm}{fc8b_u_t1p5}{$v_1$}{$v_1$}{$t=1.5$}{-.145}{.513};
  \drawFC{xshift= 0.0cm,yshift= 0.0cm}{fc8b_v_t1p5}{$v_2$}{$v_2$}{$t=1.5$}{-.230}{.199};
%
\end{tikzpicture}
\end{center}
  \caption{Parameter set I: Traveling-wave pressure-pulse in a flexible tube. Contours of $p$, $v_1$ and $v_2$ at $t=1.5$ computed
    on grid $\Gcfc^{(8)}$. 
    For clarity, the top surface curve is also 
    is shown in red with the displacement amplified by a factor of $5$.}
  \label{fig:flexibleChannelContours}
\end{figure}
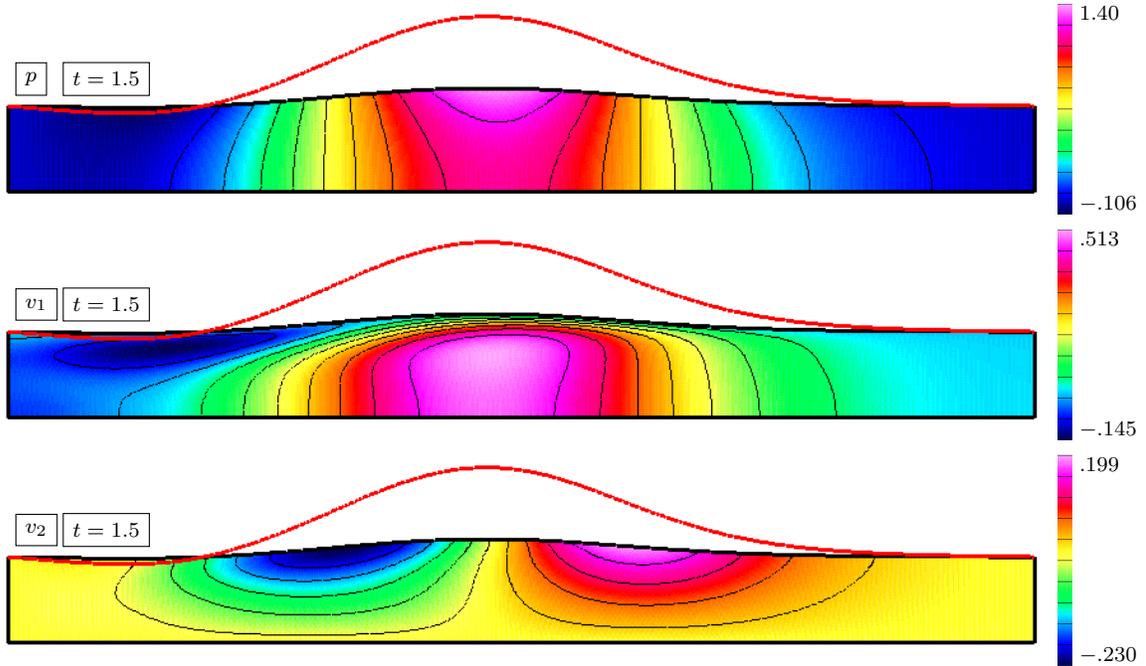
}
}

Dimensionless parameters for this set are taken to be $\rho=1$, $\mu=10^{-2}$, $\pMax=2$ and $\tMax=0.75$ for the fluid,
and 
$\rhos\hs=0.11$, 
$\Ts=5/6 \approx .83333$, $\Ks_0=40/3 \approx 13.333$,  $\Kt=\rhos\hs/100$ and $\Kxxt=5/60 \approx .08333$
for the beam.
The term with coefficient $\Ks_0$ acts as a linear restoring force while
the term with coefficient $\Kxxt$ is a {\em visco-elastic} damping term that tends to smooth high-frequency oscillations in space.
The parameters used for this first set are slightly different than the ones used in~\cite{FernandezLandajuela2014}.  Our
parameters give a larger
amplitude of the deformation of the beam and a thicker boundary layer in the fluid near the surface of the beam.  A larger
deformation of the beam is a more severe test of our DCG approach, and a thicker boundary layer enables grid-convergence tests
demonstrating second-order accuracy using coarser grids.  (The choice of
viscosity used in~\cite{FernandezLandajuela2014} requires a very fine grid to resolve the very thin
boundary layer, see Parameter Set II in Section~\ref{sec:elasticTubeFernandez}
below).  We also note that the fluid tractions and velocities are transferred directly on the
beam reference curve for the simulations described in~\cite{FernandezLandajuela2014}, whereas our approach applies
these quantities on the beam surface
as described in Section~\ref{sec:beamSolverNumerical}.
In order to approximate a beam of negligible thickness with our approach,
we choose a small value for the beam thickness of $\hs=10^{-3}$.
{
\newcommand{\figWidth}{7.5cm}
\newcommand{\figWidthb}{3.75cm}
\newcommand{\trimfig}[2]{\trimFig{#1}{#2}{.0}{0.0}{.0}{.0}}
\begin{figure}[htb]
\begin{center}
\begin{tikzpicture}[scale=1]
  \useasboundingbox (0.0,.75) rectangle (16.5,13.);  
  \draw(0.0 ,6.5) node[anchor=south west,xshift=-4pt,yshift=+0pt] {\trimfig{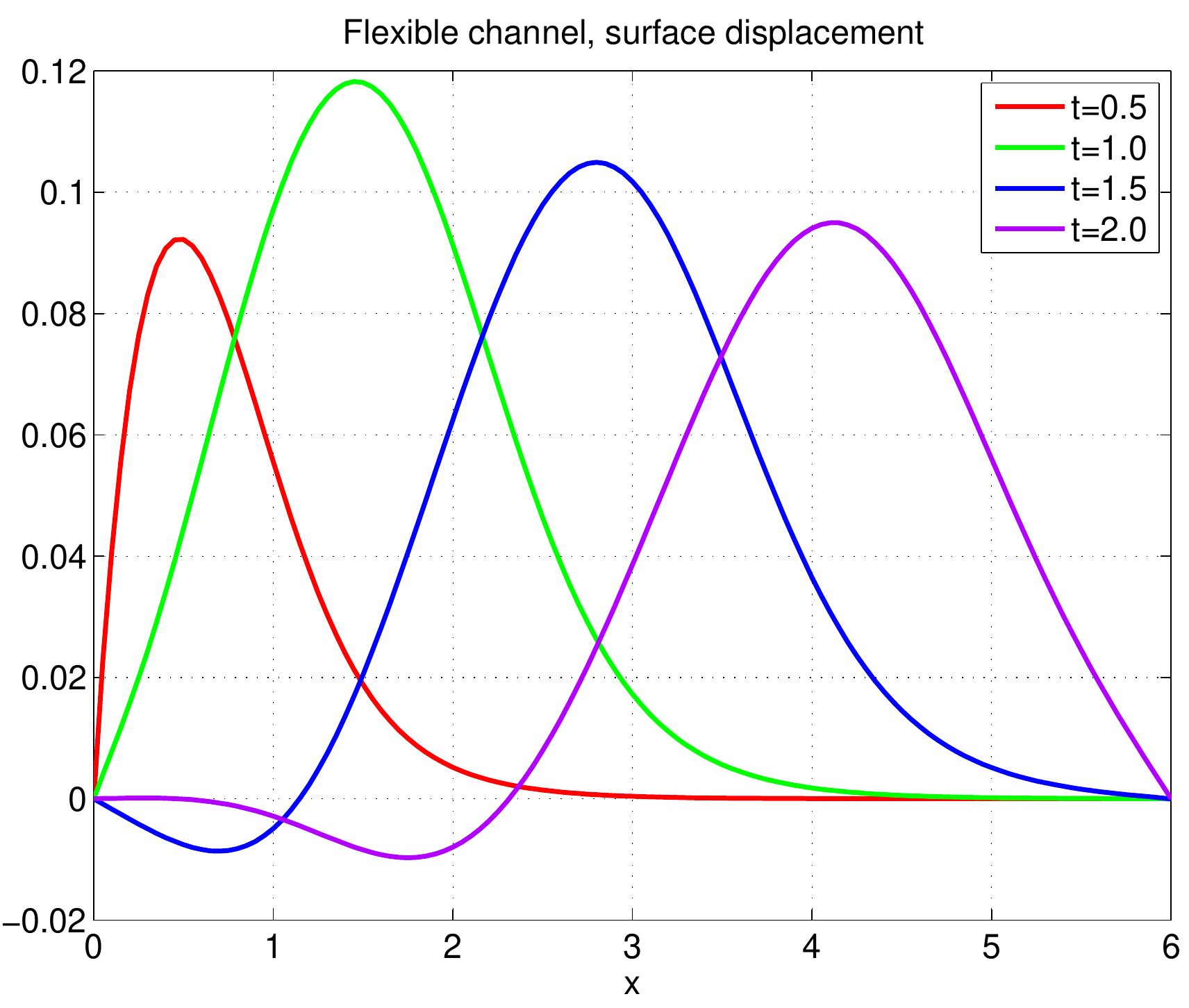}{\figWidth}};
  \draw(7.75,6.5) node[anchor=south west,xshift=-4pt,yshift=+0pt] {\trimfig{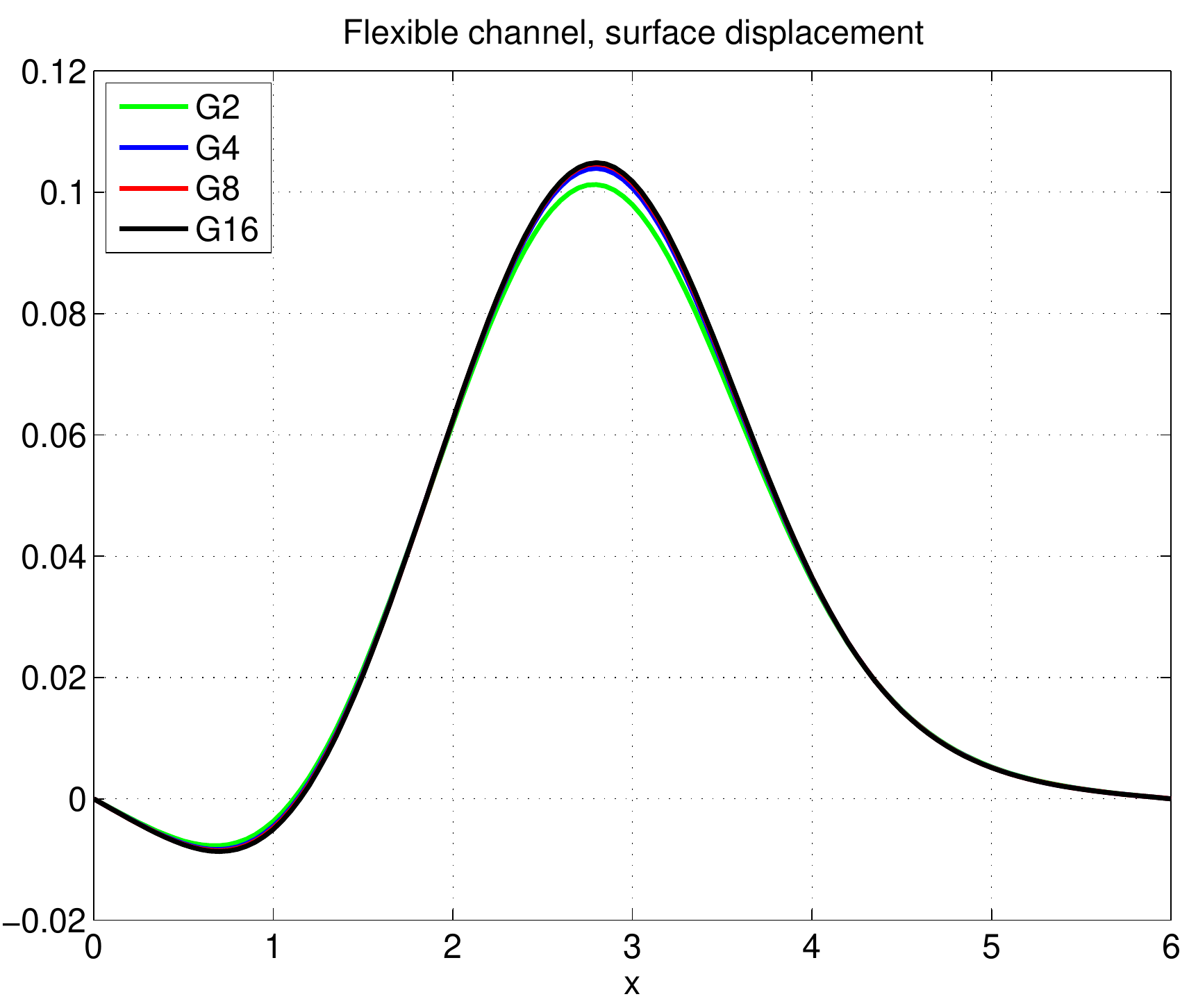}{\figWidth}};
  \draw(0.0 ,0.0) node[anchor=south west,xshift=-4pt,yshift=+0pt] {\trimfig{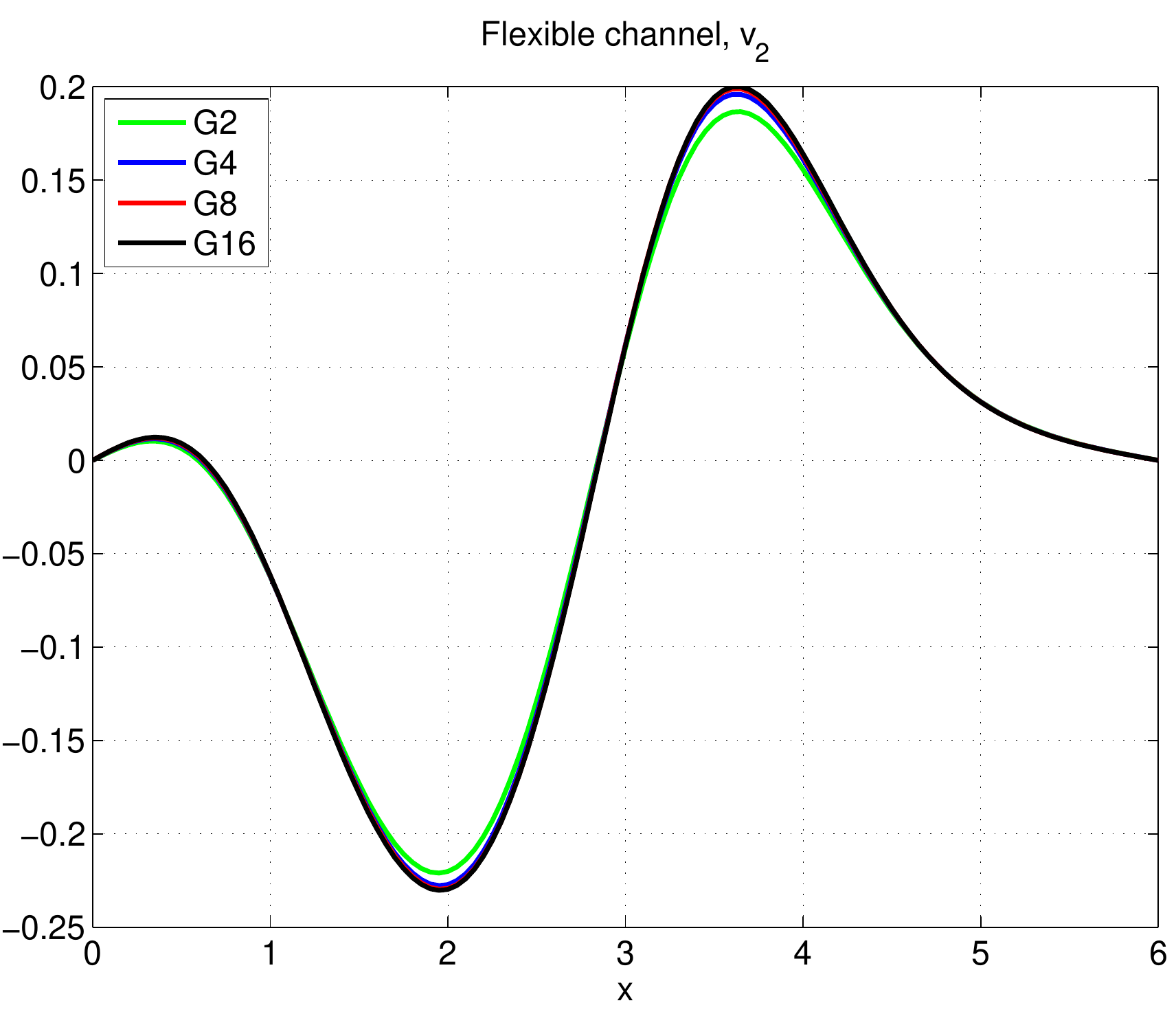}{\figWidth}};
  \draw(7.75,0.0) node[anchor=south west,xshift=-4pt,yshift=+0pt] {\trimfig{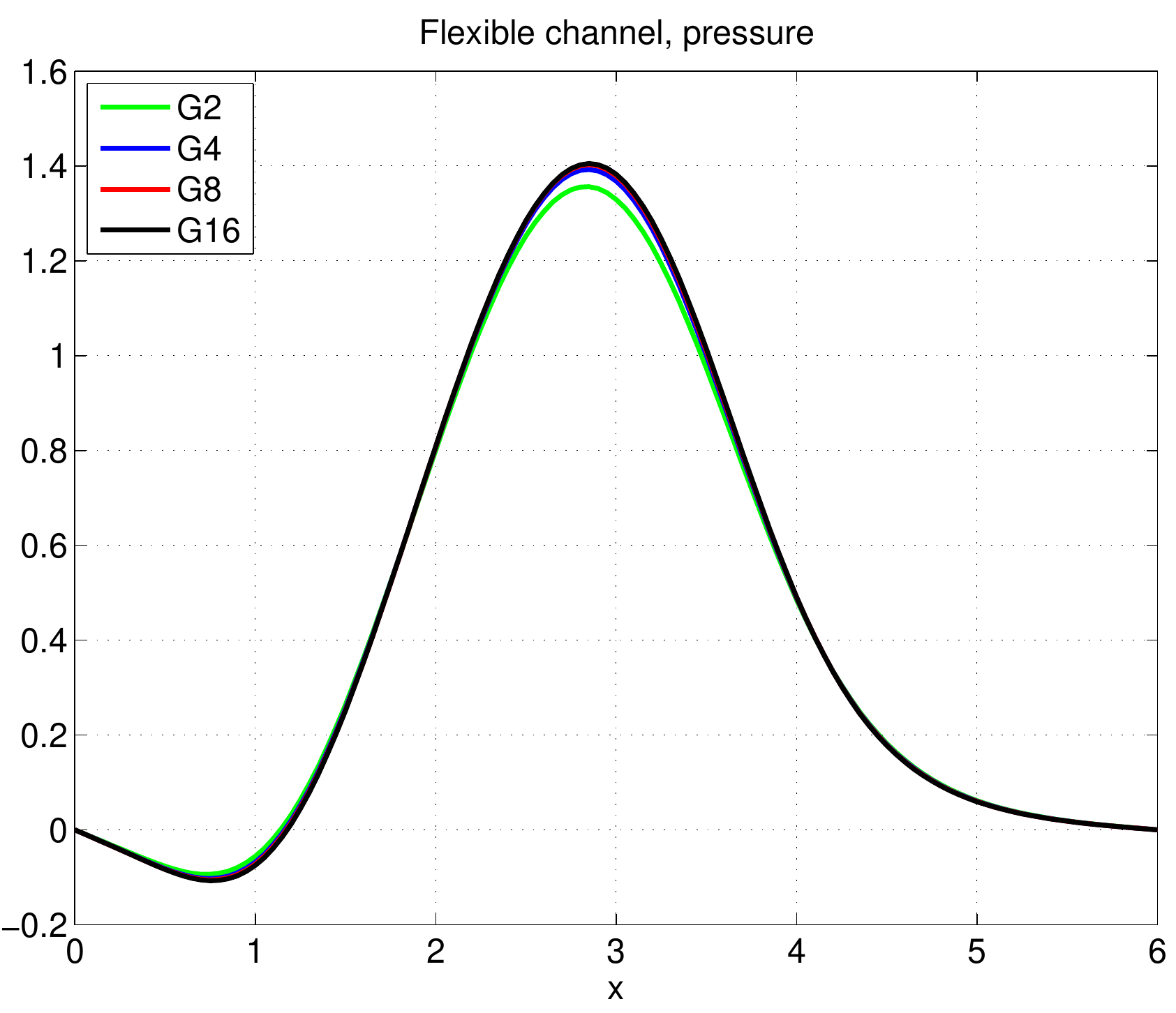}{\figWidth}};
  \draw(12.8,9.4) node[anchor=south west,xshift=-4pt,yshift=+0pt] {\trimfig{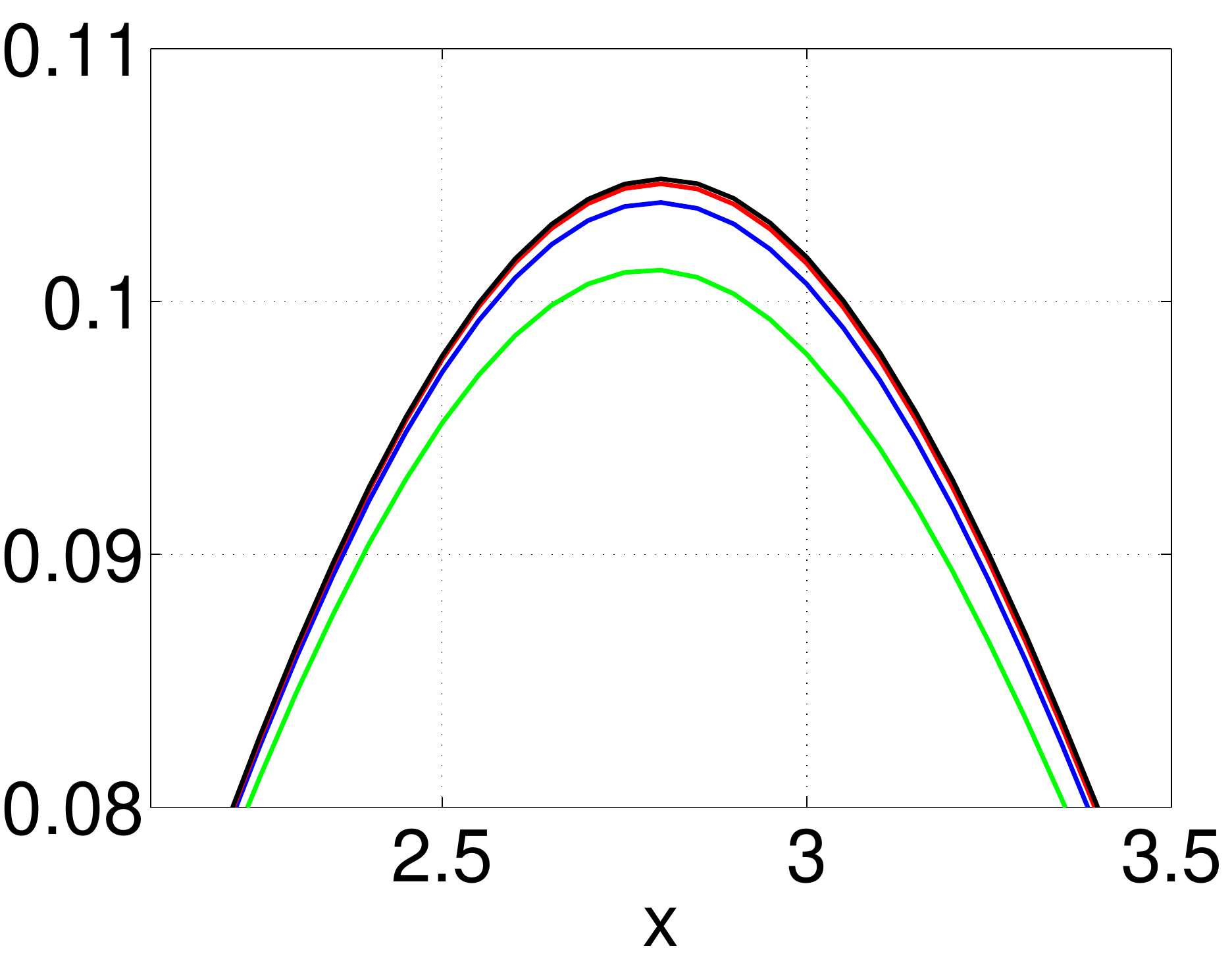}{\figWidthb}};
  \draw(3.9,.2) node[anchor=south west,xshift=-4pt,yshift=+0pt] {\trimfig{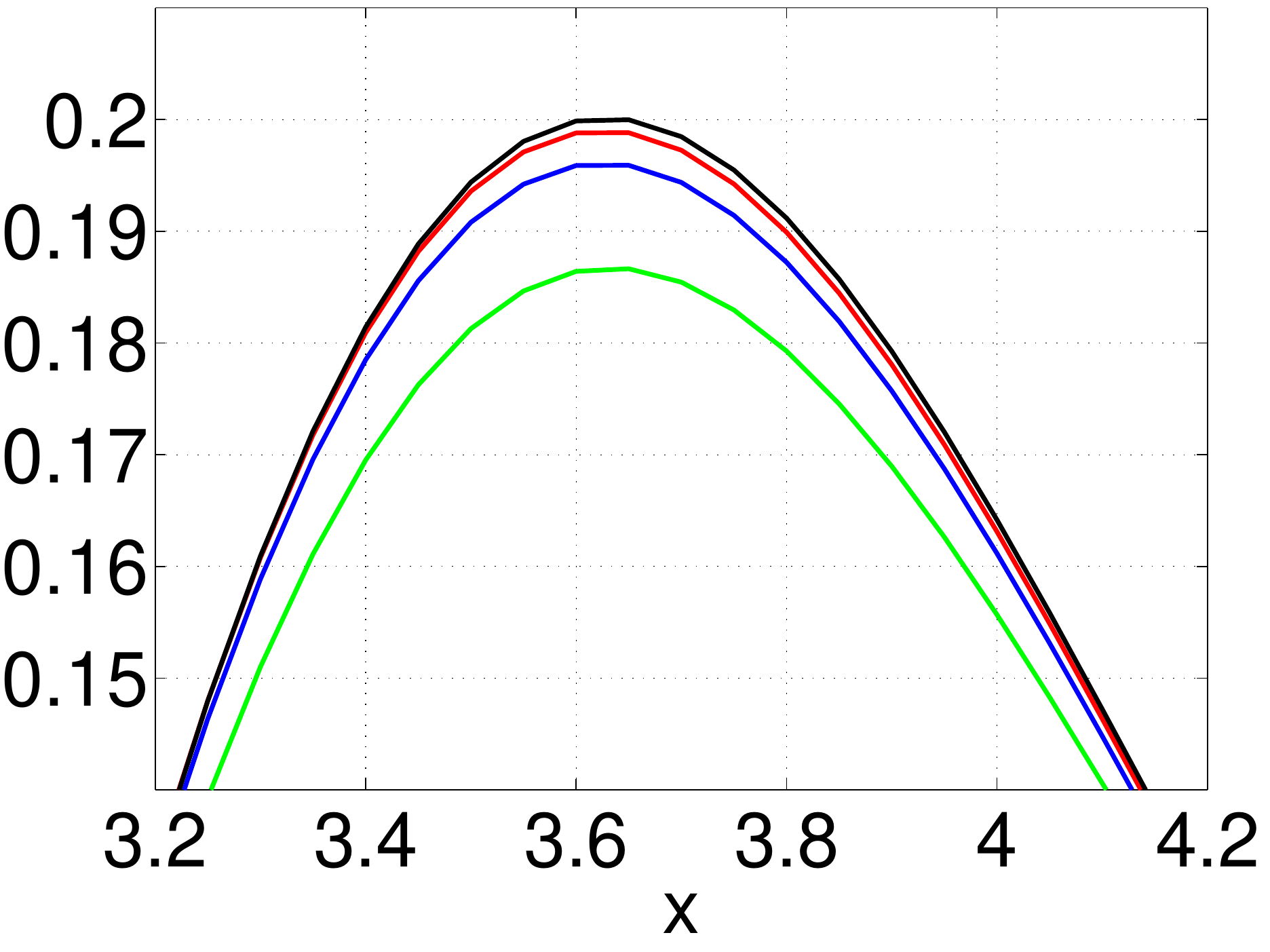}{\figWidthb}};
  \draw(12.75,2.85) node[anchor=south west,xshift=-4pt,yshift=+0pt] {\trimfig{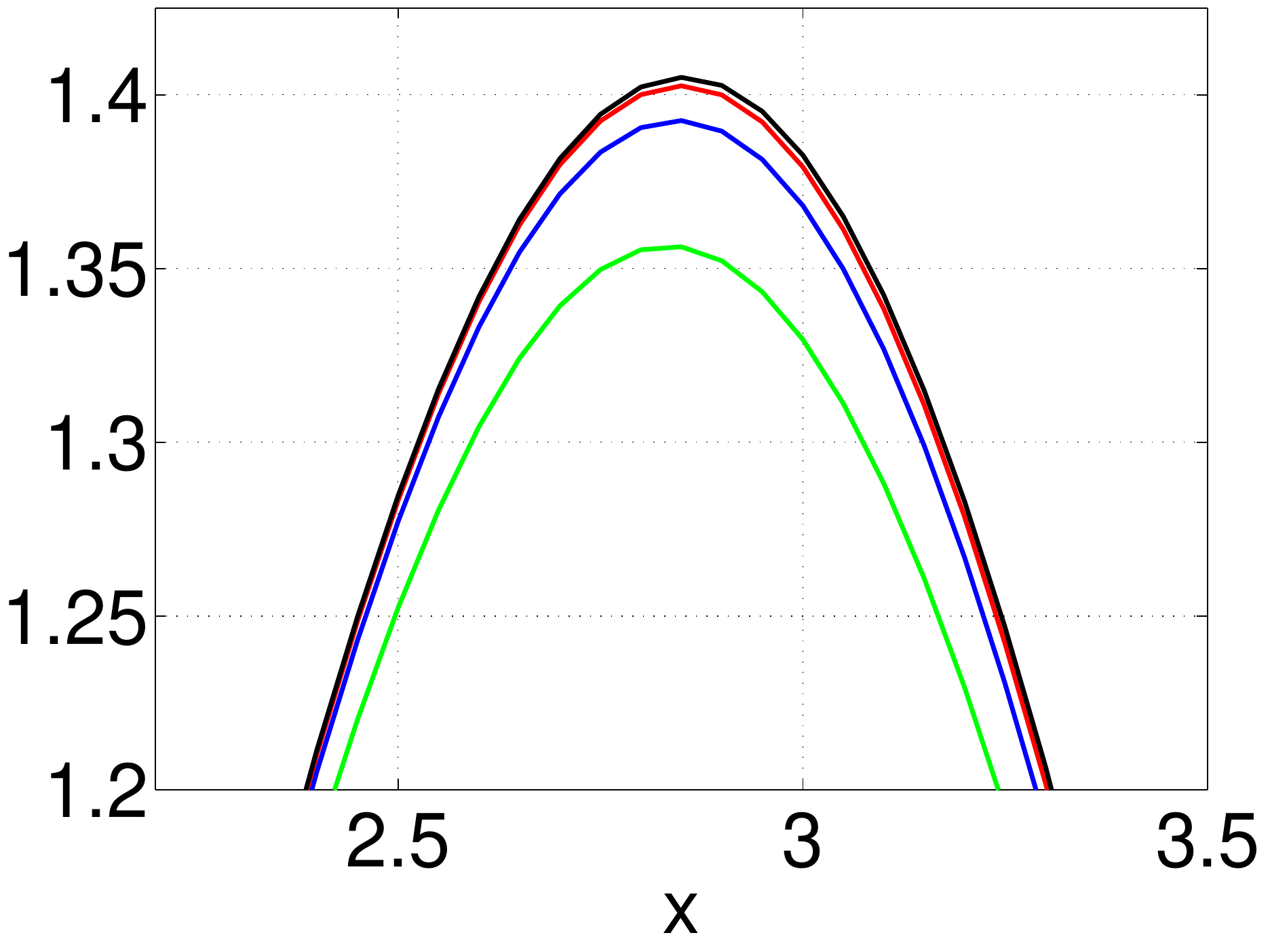}{\figWidthb}};
%
\end{tikzpicture}
\end{center}
\caption{Traveling-wave pressure-pulse in a flexible tube. 
Top left: surface displacement over time showing the formation and propogation
of the pulse. Top right: surface displacement at $t=1.5$ for grids $\Gcfc^{(k)}$, $k=2,4,8,16$ showing the
convergence of the surface as the grid is refined. Bottom left and lower right: convergence of $v_2$ and $p$, 
respectively, as the grid is refined.}
\label{fig:flexibleChannelComparisonI}
\end{figure}
}

Figure~\ref{fig:flexibleChannelContours} shows shaded contours of the fluid pressure and the components
of the fluid velocity at $t=1.5$.  The pressure pulse at this time has created a bulge in the flexible
channel near its half-length, and for clarity the centerline-curve of the beam on the top surface 
is drawn as a red curve with the displacement amplified by a factor of $5$.
The traveling pressure pulse is seen to consist of a localized high pressure region where the
beam displacement is positive. The pressure is largest at the surface near the maximum displacement of the
beam.  The horizontal component of the velocity is generally large and positive
within the pulse with a boundary layer near the top surface to match the no-slip condition there.
The vertical component of velocity is positive ahead of the pulse and negative behind it confirming
that the pulse is traveling in the positive $x_1$ direction.

The top left plot in Figure~\ref{fig:flexibleChannelComparisonI} shows the evolution of the shape of the surface at four
different times. The surface is initially displaced upwards on the left side due to the 
high pressure at inflow~\eqref{eq:pressureBC}
which acts over the time interval $[0,t_{{\rm max}}]$, where $t_{{\rm max}}=0.75$. The surface pulse then propagates to the right and
slowly decays in amplitude. The top right plot shows the beam displacement at $t=1.5$ computed using grids with increasing resolution.
The solutions are seen to converge with the curves on two finest grids being nearly indistinguishable.  The two plots
on the bottom of Figure~\ref{fig:flexibleChannelComparisonI} show the grid convergence of the vertical component of the fluid velocity
and the fluid pressure at the surface of the beam at $t=1.5$.  Here, we also observe excellent convergence of these two
fluid variables.

To provide a more quantitative estimate of the solution accuracy, a self-convergence study is
performed. Given a sequence of three or more grids of increasing resolution, a posteriori estimates
of the errors and convergence rates can be computed using the Richardson extrapolation procedure
described in~\cite{pog2008a}. These self-convergence estimates assume that the
numerical solution converges to a limiting solution, but does not assume that some chosen fine-grid
solution is the exact solution.
A posteriori estimates computed in this way are given in
the table in Figure~\ref{tab:flexibleChannelI} for four grids of increasing resolution (the convergence rate
computation uses the three finest resolutions).
The table provides estimates of the max-norm errors for various components of the solutions
and the corresponding estimates of the convergence rates.  We observe that all solution components
are converging with approximately second-order accuracy. 
Note that due to the manner in which rates are computed, 
the ratios of errors in the columns headed ``r'' are the average ratios 
and thus the same for each grid resolution. 

\bogus{
\begin{table}[hbt]\tableFont 
\begin{center}                                                                                                                                  
\begin{tabular}{|l|c|c|c|c|c|c|} \hline                                                                                                         
   grid              & \errFormat{p} &  r   & \errFormat{v_1} &  r   & \errFormat{v_2} &  r  \\ \hline                                              
 fc2bfixed.show  & \num{5.0}{-2} &      & \num{5.4}{-2} &      & \num{1.4}{-2} &      \\ \hline                                                  
 fc4bfixed.show  & \num{1.4}{-2} &  3.6 & \num{1.4}{-2} &  3.8 & \num{3.6}{-3} &  4.0 \\ \hline                                                  
 fc8bfixed.show  & \num{3.8}{-3} &  3.6 & \num{3.8}{-3} &  3.8 & \num{9.1}{-4} &  4.0 \\ \hline                                                  
 fc16bfixed.show & \num{1.0}{-3} &  3.6 & \num{1.0}{-3} &  3.8 & \num{2.3}{-4} &  4.0 \\ \hline                                                 
   rate          &     1.87      &      &     1.91      &      &     1.99      &     \\ \hline                                             
\end{tabular}                                                                                                                                   
\caption{Max-norm self convergence results, fixed-width grids -- Mon Mar 30 10:34:05 2015. } 
\label{fig:flexibleChannelI}
\end{center}                                                                                                                                    
\end{table}   
}
\newcommand{\tableChannelPSI}{%
\begin{tabular}{|c|c|c|c|c|c|c|} \hline
\multicolumn{7}{|c|}{\strutt Traveling wave pressure-pulse} \\ \hline 
\strutt~~$h_j$~~& $E_j^{(p)}$ & $r$  & $E_j^{(v_1)}$  & $r$    & $E_j^{(v_2)}$   & $r$    \\ \hline 
 1/20        & \num{5.0}{-2} &      & \num{5.4}{-2} &      & \num{1.4}{-2} &      \\ \hline   
 1/40        & \num{1.4}{-2} &  3.6 & \num{1.4}{-2} &  3.8 & \num{3.6}{-3} &  4.0 \\ \hline   
 1/80        & \num{3.8}{-3} &  3.6 & \num{3.8}{-3} &  3.8 & \num{9.1}{-4} &  4.0 \\ \hline   
 1/160       & \num{1.0}{-3} &  3.6 & \num{1.0}{-3} &  3.8 & \num{2.3}{-4} &  4.0 \\ \hline   
\rateLabel   &     1.87      &      &     1.91      &      &     1.99      &     \\ \hline    
\end{tabular}
}

\begin{figure}[hbt]\tableFontSize
\begin{center}
\tableChannelPSI
\end{center}
\caption{Traveling wave pressure-pulse (parameter set I): Maximum errors and 
estimated convergence rates at $t=1.5$ computed using the AMP scheme.
}
\label{tab:flexibleChannelI}
\end{figure}

Figure~\ref{fig:flexibleChannelComparisonTP} compares results from the AMP scheme with the
traditional-partitioned scheme with sub-time-step iterations (TP-SI) on the coarse grid $\Gcfc^{(2)}$. 
The number of sub-iterations per time step required for the TP-SI scheme is approximately $15$ on average, although near the start
of the run almost $250$ iterations are needed.  (A convergence tolerance of $10^{-4}$ on the
sub-iterations is used.)  A simple under-relaxed fixed-point iteration is used with the value for the
under-relaxation parameter taken as $\omega=0.025$ (the fact that such a small value is needed indicates that the problem has
large added-mass effects). 
The solutions from the AMP and TP-SI are in excellent agreement.

{
\newcommand{\figWidth}{7cm}
\begin{figure}
\begin{center}
\includegraphics[width=\figWidth]{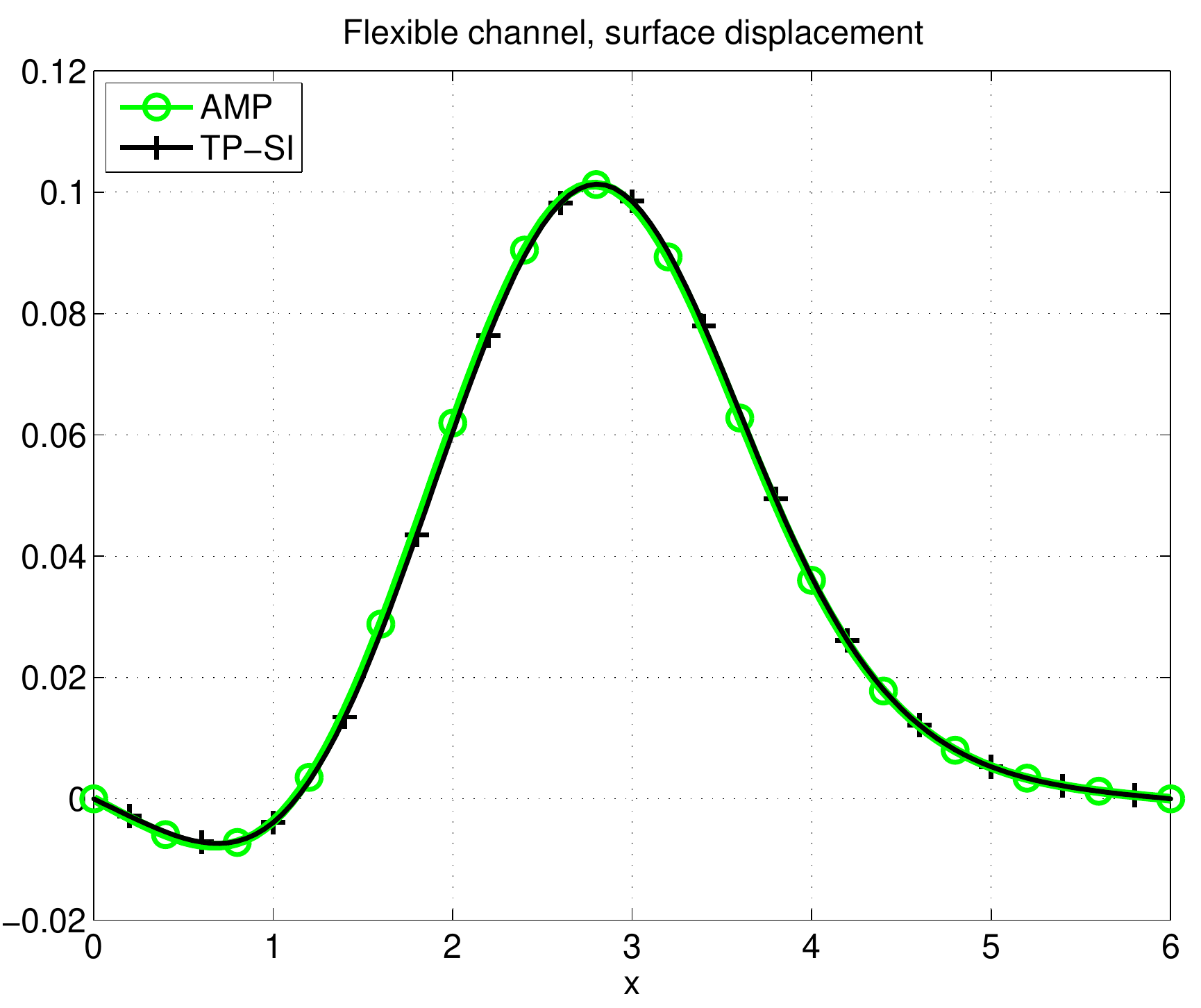}
\includegraphics[width=\figWidth]{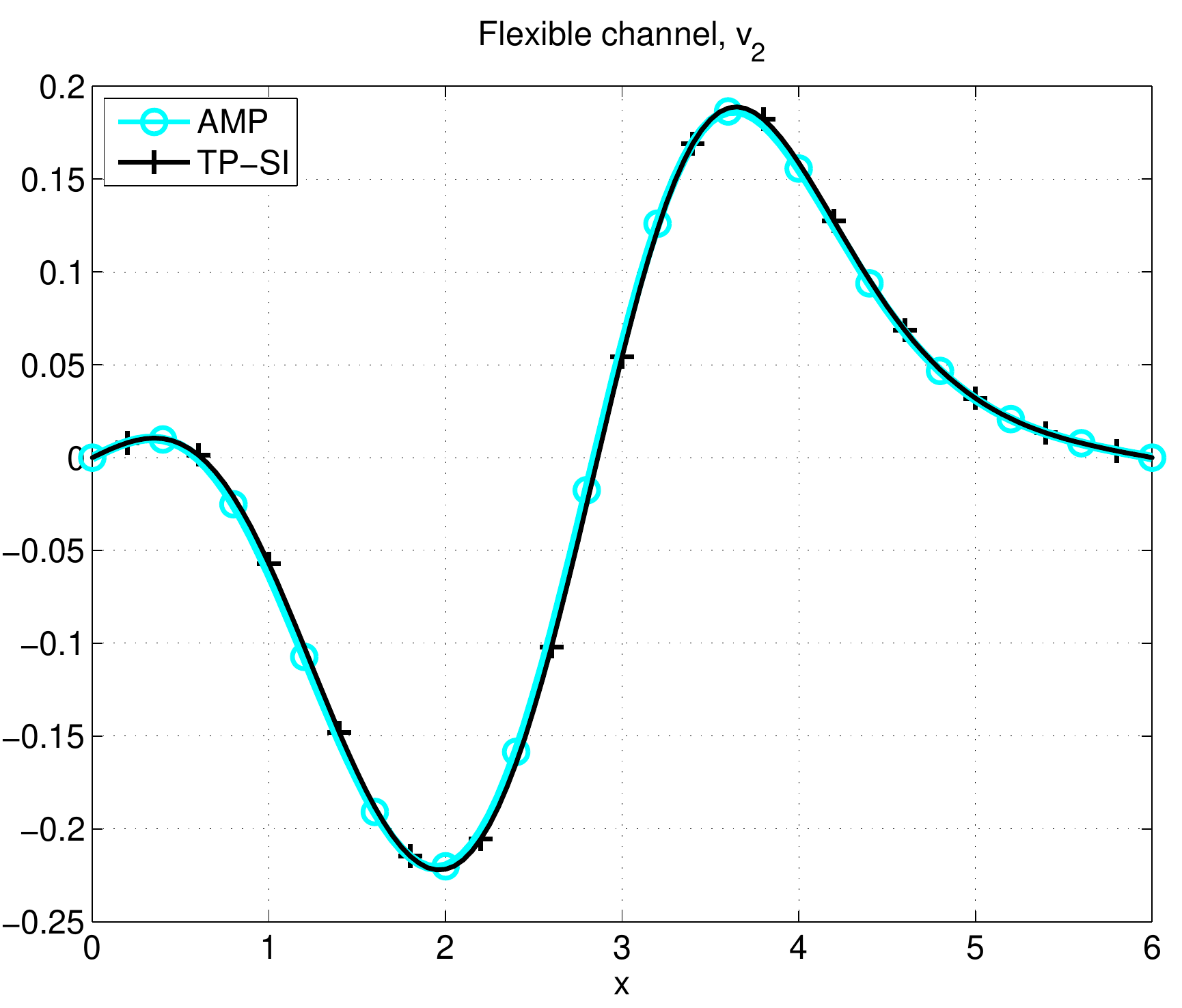}
\includegraphics[width=\figWidth]{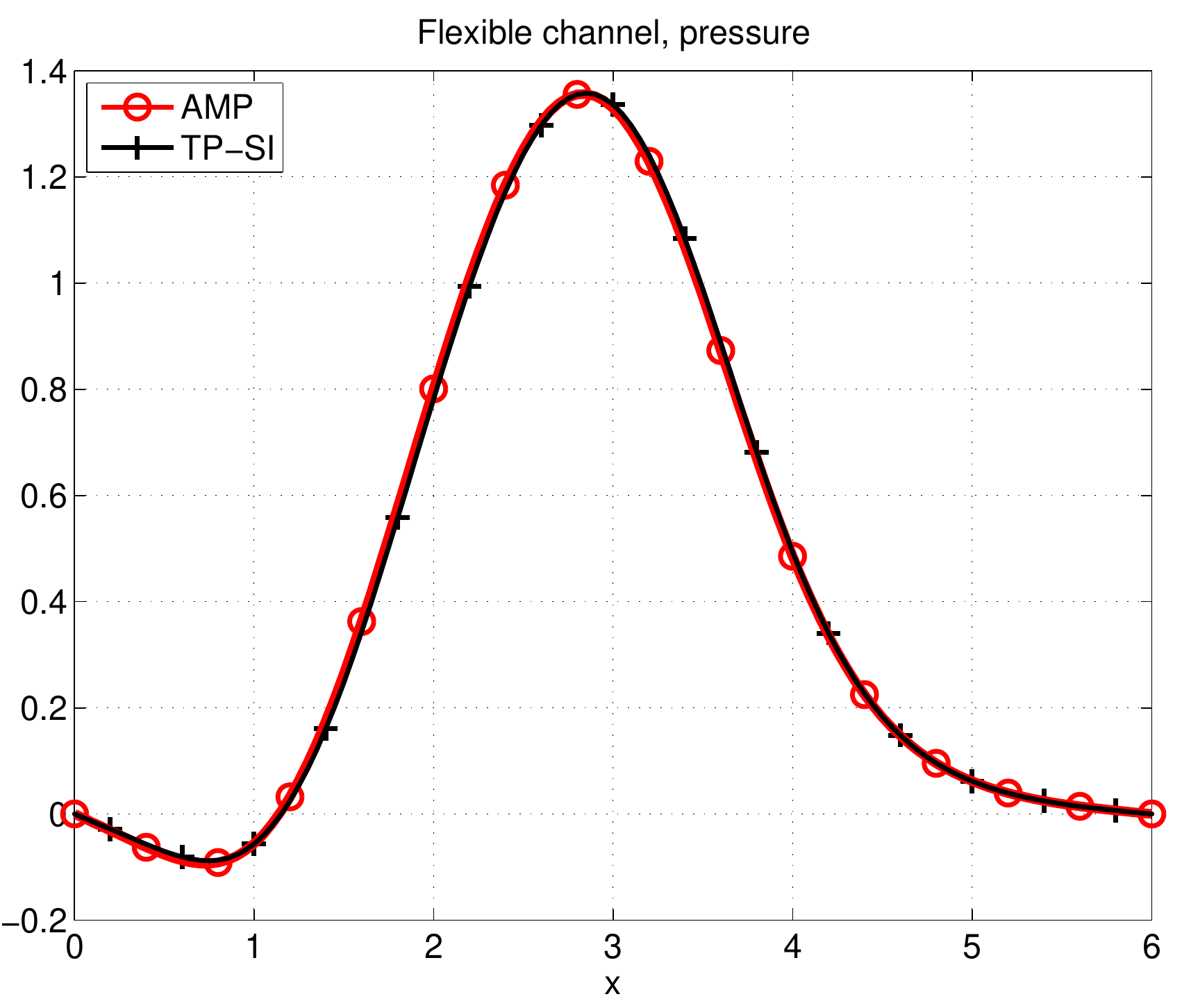}
\end{center}

\caption{Traveling-wave pressure-pulse in a flexible tube. 
A comparison between the AMP scheme and the TP scheme with sub-time-step iterations at $t=1.5$ with grid $\Gcfc^{(2)}$.
The surface displacement, velocity $v_2$ and pressure are shown (with markers placed every $8$ grid points). 
The results of the two schemes are seen to be nearly indistinguishable.
}
\label{fig:flexibleChannelComparisonTP}
\end{figure}
}

\subsubsection{Parameter Set II}  \label{sec:elasticTubeFernandez}


{
\newcommand{\figWidth}{9.75cm}
\newcommand{\figWidthb}{5.5cm}
\newcommand{\trimfig}[2]{\trimFig{#1}{#2}{.0}{0.0}{.0}{.0}}
\begin{figure}[htb]
\begin{center}
\begin{tikzpicture}[scale=1]
  \useasboundingbox (0.0,.75) rectangle (15,8.25);  
  \draw(0.0 ,0.0) node[anchor=south west,xshift=-4pt,yshift=+0pt] {\trimfig{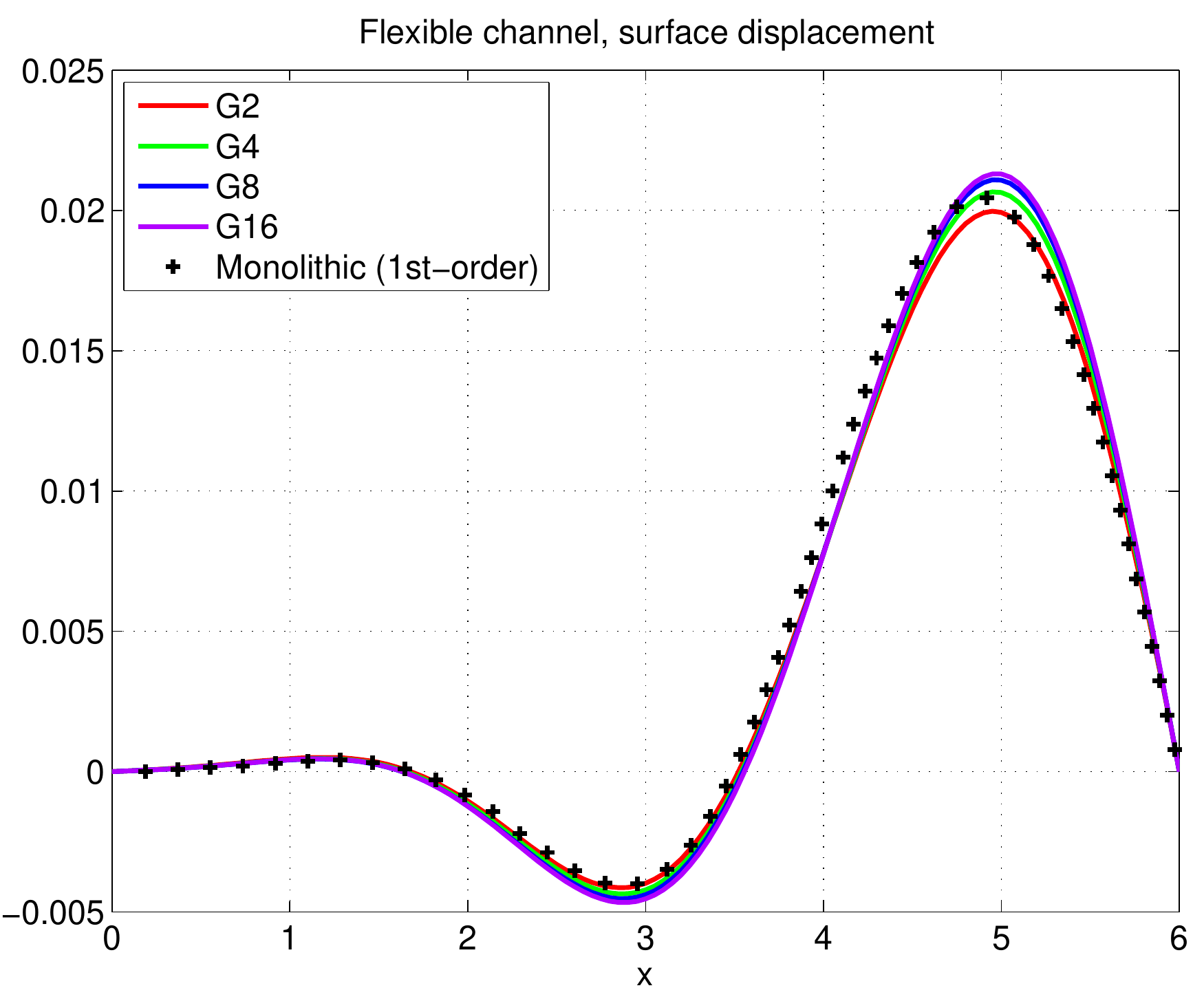}{\figWidth}};
  \draw(9.95,2.5) node[anchor=south west,xshift=-4pt,yshift=+0pt] {\trimfig{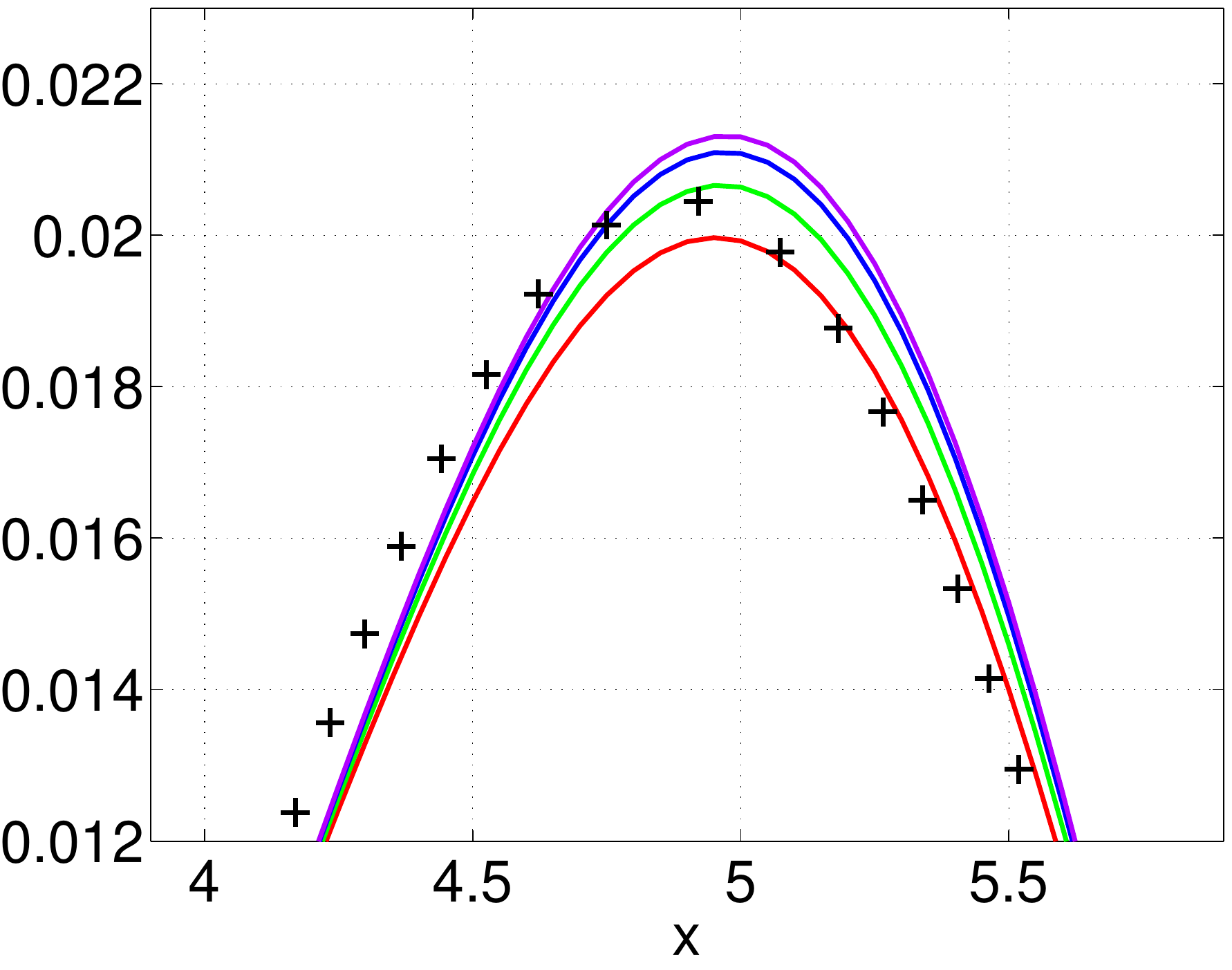}{\figWidthb}};
%
\end{tikzpicture}
\end{center}
\caption{Paramater set II: Traveling-wave pressure-pulse in a flexible tube. Surface displacement as the
grid are refined. Results for monolithic scheme are from Fernandez et.al.~\cite{FernandezLandajuela2014}.}
\label{fig:flexibleChannelComparison}
\end{figure}
}

The parameter values for this set match those of Fernandez et~al.~\cite{FernandezLandajuela2014} so that a
direct comparison of the results can be made.  We prefer to work with dimensionless parameters obtained by scaling with the reference
length $L_0=1\;{\rm cm}$, velocity $U_0=10^2\;{\rm cm}/{\rm s}$, density $\rho_0=1\;{\rm gm}/{\rm cm}^3$, and time $T_0=L_0/U_0=10^{-2}\;{\rm s}$.
These scalings give the corresponding dimensionless parameters $\rho=1$, $\mu=3.5\times 10^{-4}$, $\pMax=2$ and $\tMax=0.5$ for the fluid,
and $\rhos\hs=0.11$, $\Ts=2.5$, $\Ks_0=40$,  $\Kt=1.1\times10^{-3}$ and $\Kxxt=0.25$ for the beam.
As in the previous set of parameters, the product $\rhos\hs=0.11$ is fixed
but $\hs=10^{-3}$ and $\rhos=110$ are chosen to reduce the effects of a finite-width beam.

Figure~\ref{fig:flexibleChannelComparison} compares the surface displacement of the beam 
at $t=1.5$ computed with different grids 
using the AMP-PBA algorithm. The composite grids used in these calculations are the same as those used
in the previous section.
We observed that the surface displacement appears to be converging in a manner
consistent with second-order accuracy. 
The figure also shows the results
from Fernandez et~al.~\cite{FernandezLandajuela2014} using a first-order accurate monolithic scheme
on a grid with spacing $h=1/320$ (equal to the spacing on grid $\Gcfc^{(32)}$). 
Although the monolithic results do not exactly match the fine grid AMP results,
we note that
the monolithic results do not yet appear to be grid converged since they differ by a fair amount from the monolithic results
on a grid coarsened by a factor of 2, see~\cite{FernandezLandajuela2014}.  


{
\newcommand{\drawFC}[7]{%
\begin{scope}[#1]
\draw(0.0,0) node[anchor=south west,xshift=-4pt,yshift=+0pt] {\trimfig{flexibleChannel/#2}{\figWidth}};
\draw(.5,2.35) node[draw,fill=white,anchor=west,xshift=2pt,yshift=-8pt] {\scriptsize #3};
\draw(1.2,2.35) node[draw,fill=white,anchor=west,xshift=0pt,yshift=-8pt] {\scriptsize #5};
\draw (\xcb,\ycb) node[anchor=south west,xshift= +0pt,yshift=+0pt] {\trimfigcb{fig/colourBarLines}{\cbWidth}{\cbHeight}};
\draw (\xcb,\ycb) node[anchor=south west,xshift= +8pt,yshift=+1pt] {\scriptsize $#6$};
\draw (\xcb,\ycbTop) node[anchor=south west,xshift= +8pt,yshift=-6pt] {\scriptsize $#7$};
\end{scope}
}
{
\newcommand{\cbWidth}{.2cm}
\newcommand{\cbHeight}{2.8cm}
\newcommand{\xcb}{14.3cm}
\newcommand{\ycb}{-.15cm}
\setlength{\ycbTop}{\ycb+\cbHeight}
\setlength{\ycbMid}{\ycb+\cbHeight*\real{.5}}
\newcommand{\xLabel}{6.5cm}
\newcommand{\yLabel}{6.5cm}
\newcommand{\trimfigcb}[3]{\includegraphics[width=#2, height=#3, clip, trim=17cm 2.35cm 1.65cm 2.35cm]{#1}}
\newcommand{\figWidth}{14.5cm}
\newcommand{\trimfig}[2]{\trimFig{#1}{#2}{.03}{.15}{.33}{.34}}

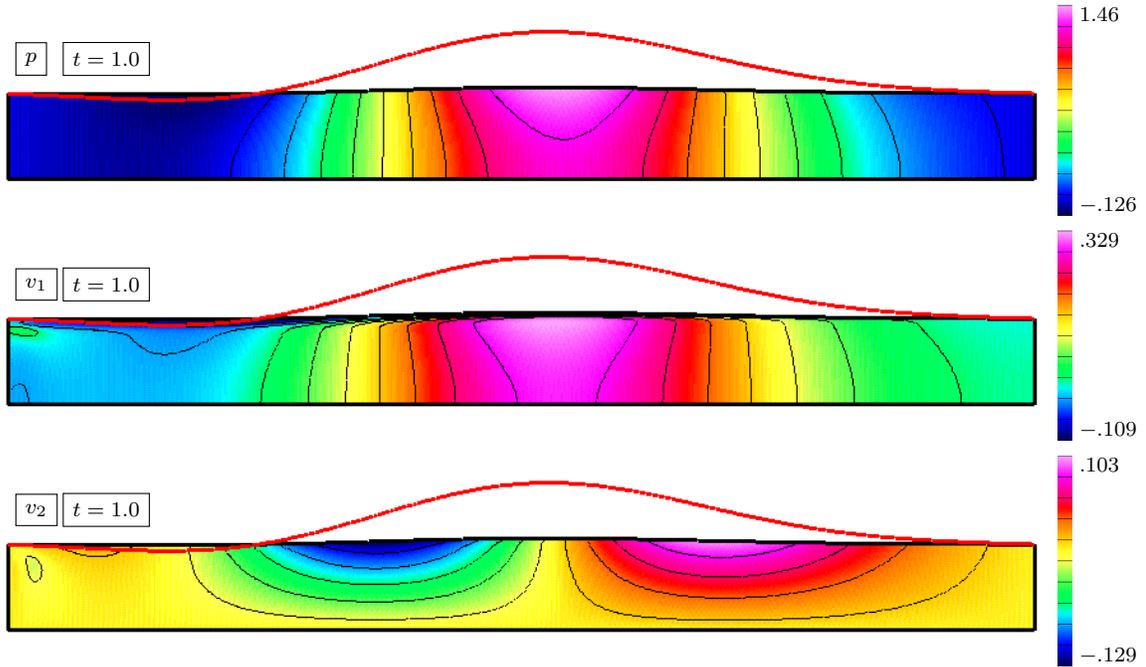
\begin{figure}[htb]
\begin{center}
\begin{tikzpicture}[scale=1]
  \useasboundingbox (0.0,.5) rectangle (15.25,9);  
  \drawFC{xshift= 0.0cm,yshift= 6.0cm}{fc8an_p_t0p01}{$p$}{$p$}{$t=1.0$}{-.126}{1.46};
  \drawFC{xshift= 0.0cm,yshift= 3.0cm}{fc8an_u_t0p01}{$v_1$}{$v_1$}{$t=1.0$}{-.109}{.329};
  \drawFC{xshift= 0.0cm,yshift= 0.0cm}{fc8an_v_t0p01}{$v_2$}{$v_2$}{$t=1.0$}{-.129}{.103};
%
\end{tikzpicture}
\end{center}
  \caption{Parameter set II: Traveling-wave pressure-pulse in a flexible tube. Contours of $p$, $v_1$ and $v_2$ at $t=1.0$ computed
    on grid $\Gcfc^{(8)}$. 
    For clarity, the top surface curve is also 
    is shown in red with the displacement amplified by a factor of $10$.}
  \label{fig:flexibleChannelContoursCaseI}
\end{figure}
}
}

Figure~\ref{fig:flexibleChannelContoursCaseI} shows shaded contours and surface displacements of the solution at $t=1.0$ computed using the AMP algorithm.
We note that the behavior of the solution is qualitatively similar to that observed previously
for Parameter Set I.  
The contour plot of the horizontal component of velocity, $v_1$, however, shows the very thin boundary layer
that forms near the interface. 

%
%

\subsection{Beam in a cross flow}  \label{sec:beamInAChannel}

As a final demonstration of our AMP algorithm, we consider flow past a flexible beam that partially blocks 
a fluid channel as shown in Figure~\ref{fig:beamInAChannelLight}.
A flow from left to right causes the beam to initially bend to the right and then oscillate. 
Eventually the flow approaches a steady-state consisting
of a long recirculation
region downstream of the deflected beam.
The computed solution from the AMP scheme for a light beam is compared with the solution obtained using the TP-SI scheme.
Results for beams of different densities are also shown and used to estimate the added mass for this configuration.


\newcommand{\xa}{x_a}
\renewcommand{\xb}{x_b}
\newcommand{\ya}{y_a}
\newcommand{\yb}{y_b}

The geometry of the problem is illustrated in Figure~\ref{fig:beamInAChannelLight} and
consists of an initially vertical beam, with length $\ls=1$, that extends into a rectangular
fluid channel of height $\yb=2$ and overall dimensions $[\xa,\xb]\times[\ya,\yb]=[-3,8]\times[0,2]$. 
The thickness of the beam at its clamped base is $\hs(0)=0.2$ (the beam reference line is located
along the centerline) and
this remains approximately constant along the beam until its rounded 
free end\footnote{ 
The reason that the beam at the tip is rounded off, rather than having sharp corners, 
is so that a high-quality grid can be constructed that resolves the rapid variation of the solution around the tip.
In addition, the solution and its derivatives will remain smooth near the rounded tip and this will give
greater confidence that the numerical solution is converging to the true solution in this region.}.
Note that while the beam reference-line (red curve in Figure~\ref{fig:beamInAChannelLight})
has an undeformed length of $\ls=1$, the actual distance from the base of the beam 
to the rounded tip at the top is slightly longer\footnote{
The actual boundary and rounded end of the beam, and hence its thickness, is defined in terms of hyperbolic trigonometric
functions through the {\em SmoothedPolygon} mapping~\cite{MAPPINGS}. The precise definition of the grid along
with scripts to run the simulations presented here
will be available with the Overture software at {\tt overtureFramework.org}.
}.
%
The base of the beam is at the origin of the computational domain, and the
upstream vertical boundary of the channel is at $\xa=-3$ while the downstream boundary is at $\xb=8$.
The grid for this geometry, denoted by $\Gbcf^{(j)}$,
was shown previously in Figure~\ref{fig:beamGrids},
and consists of a background Cartesian grid in the fluid channel together with a boundary-fitted hyperbolic grid that extends
around the beam.
The hyperbolic grid is generated with $6$ grids lines in the normal direction, and thus this grid becomes thinner as
the composite grid is refined as in previous examples.
The grid spacing is chosen to be approximately $h_j=1/(10 j)$. 

The parameters for the fluid are $\rho=1$ and $\mu=0.05$, while the parameters for the beam are $\Es\Is=5$, $\Ks=\Ts=\Kt=\Kxxt=0$,
and with a range of values for the beam density $\rhos$ as noted below.  The conditions on the boundaries of the fluid channel are
chosen as a no-slip wall on the bottom boundary at $\yy=0$, a slip-wall on the top boundary at $\yy=\yb$, inflow on the left at
$\xx=\xa$, and outflow on the right at $\xx=\xb$. The outflow boundary condition is a Robin condition on $p$ and extrapolation of
the velocity components.  The inflow velocity profile is taken to be uniform over most of the inflow boundary but parabolic near the
bottom to match the no-slip boundary condition on the lower wall. The inflow profile is given as
\begin{align}
    v_1(\xb,\yy,t) =  V_{\text{in}}(\yy) R(t),\qquad V_{\text{in}}(\yy)=
         \begin{cases}
            V_{\text{max}}      & \text{for $\yy \ge d$}, \\
           V_{\text{max}}(\yy/d)^2  & \text{for $ \yy <d$},
               \end{cases}    \label{eq:inflowProfile}
\end{align}
where $V_{\text{max}}=1$, $d=0.2$ and $R(t)$ is the ramp function given previously in~\eqref{eq:ramp} so that the horizontal component of the inflow velocity transitions smoothly in time from $v_1=0$ at $t=0$ to $v_1=V_{\text{in}}(\yy)$ for $t>1$.  The vertical component of the inflow velocity at $\xx=\xa$ is zero.

{
\newcommand{\figWidth}{14.cm}
\newcommand{\trimfig}[2]{\trimFig{#1}{#2}{.0}{.0}{.325}{.325}}
\newcommand{\figWidtha}{10.cm}
\newcommand{\trimfiga}[2]{\trimFigb{#1}{#2}{.0}{.0}{.25}{.3}}
\def\xa{.22}
\def\xb{3.925}
\def\xc{13.85}
\def\ya{.42}
\def\yb{2.9}
\def\ba{4.8}
\def\bb{5.15}
\def\bc{5.5}
\begin{figure}[htb]
\begin{center}
\begin{tikzpicture}[scale=1]
  \useasboundingbox (0.0,.0) rectangle (14.,9);  
  \draw(0.0,6.0) node[anchor=south west,xshift=-4pt,yshift=+0pt] {\trimfig{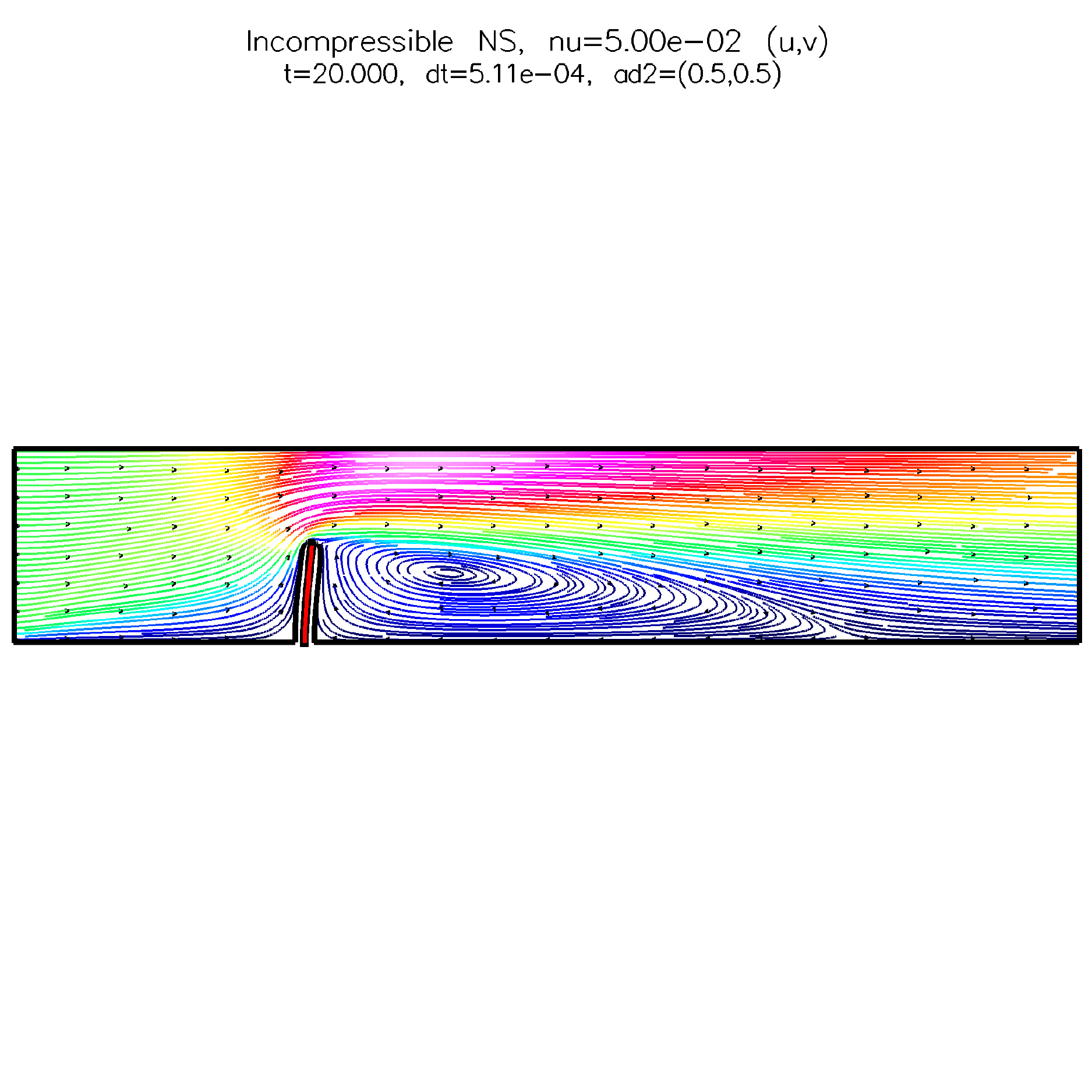}{\figWidth}};
  \draw(2.0,-.3) node[anchor=south west,xshift=-4pt,yshift=+0pt] {\trimfiga{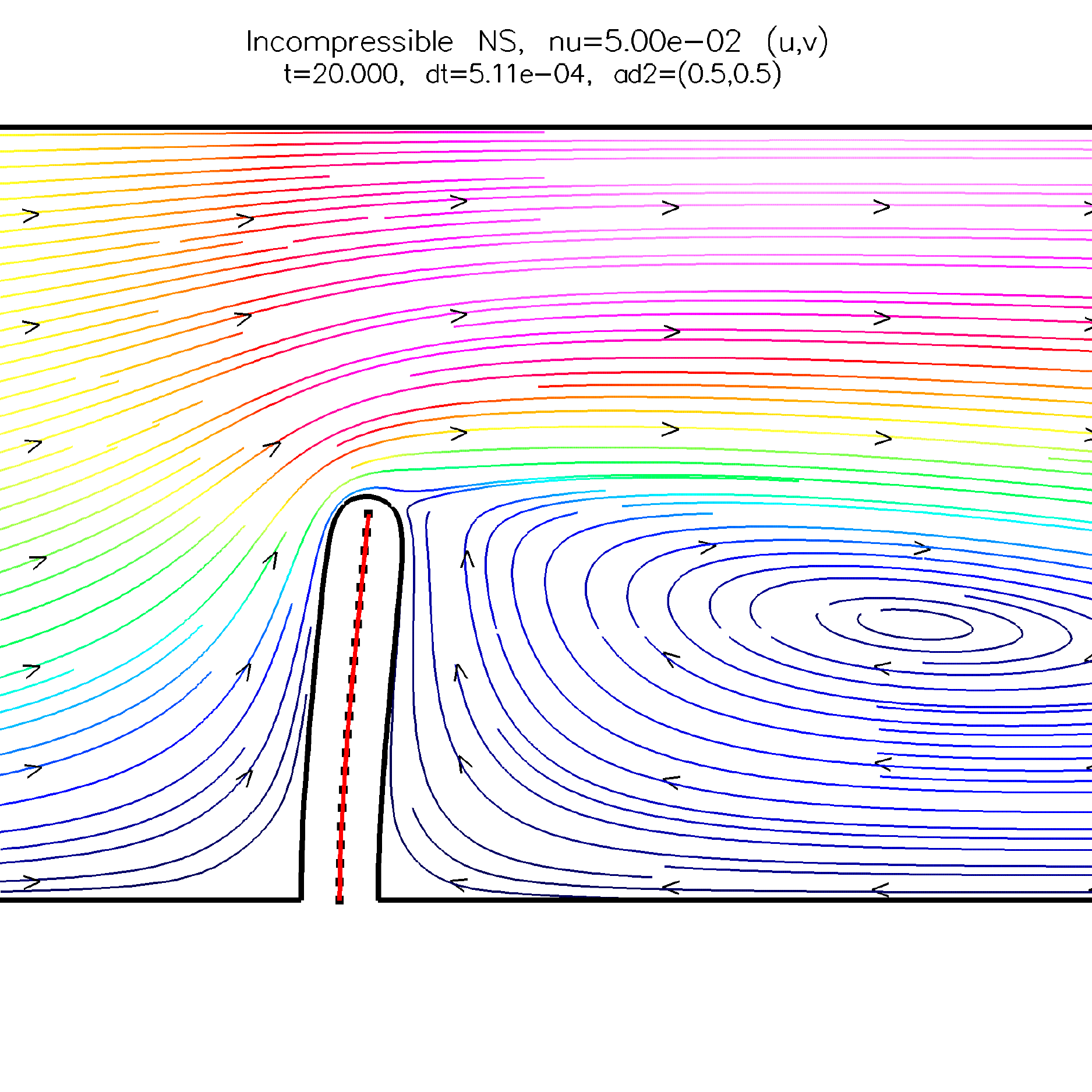}{\figWidtha}};
%
\begin{scope}[yshift=6cm]
  \begin{scope}[yshift=5pt]
    \draw[-,black,very thick] (\xa,0) node[anchor=north,yshift=0pt] {$-3$} -- 
                              (\xb,0) node[anchor=north,yshift=0pt] {$0$} -- 
                              (\xc,0) node[anchor=north,yshift=0pt] {$8$};
  \foreach \x in {\xa,\xb,\xc}
    {
     \draw[-,black,very thick] (\x,-.1) -- (\x,.1);
    }
  \end{scope}
  \draw[-,black,thick,xshift=0pt] (0,\ya) node[anchor=east,yshift=0pt] {$0$} -- 
                                  (0,\yb) node[anchor=east,yshift=0pt] {$2$};
  \foreach \y in {\ya,\yb}
    {
     \draw[-,black,very thick] (-.1,\y) -- (.1,\y);
    }
\end{scope}
%
\begin{scope}[yshift=-.3cm]
  \draw[-,black,very thick] (\ba,0) node[anchor=north,xshift=-5pt] {\scriptsize $-.1$} -- 
                            (\bb,0) node[anchor=north,yshift=0pt] {\scriptsize $0$} --
                            (\bc,0) node[anchor=north,xshift=2pt] {\scriptsize $.1$};
   \foreach \x in {\ba,\bb,\bc}
    {
     \draw[-,black,very thick] (\x,-.1) -- (\x,.1);
    } 
\end{scope}
\end{tikzpicture}
\end{center}
  \caption{Light beam with $\rhos=1$ in a cross flow: 
    Streamlines of the flow near steady state at $t=20$ computed on the composite grid $\Gbcf^{(4)}$.}
  \label{fig:beamInAChannelLight}
\end{figure}
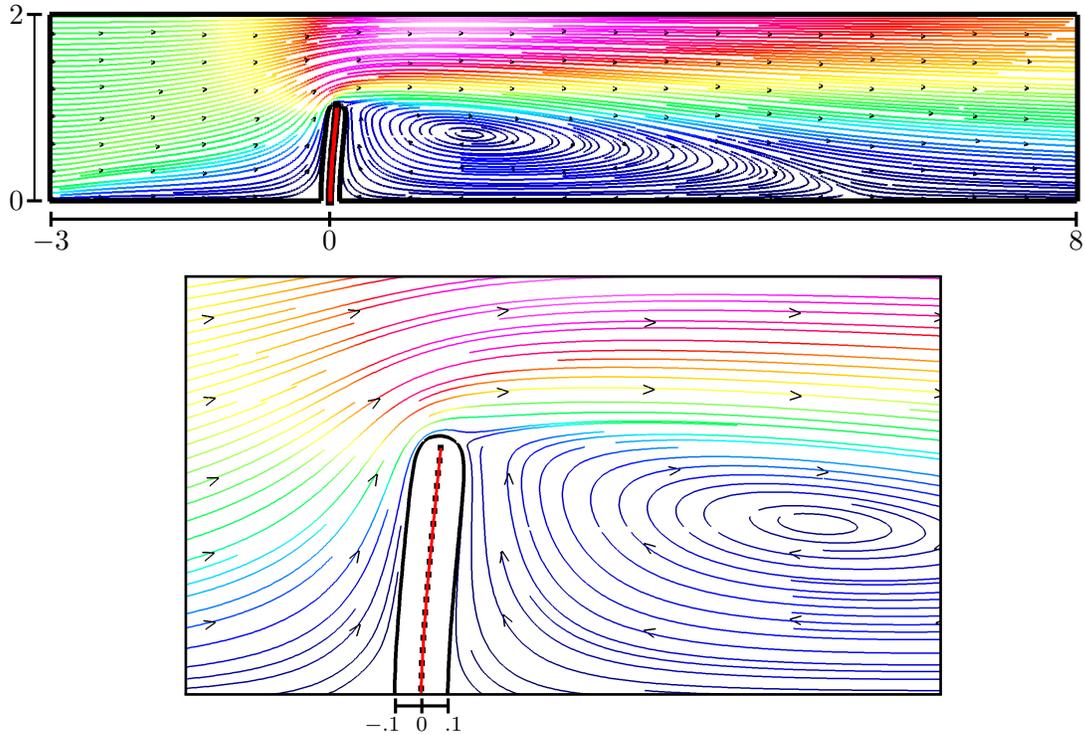
}

Figure~\ref{fig:lightBeamInAChannelTip} shows the behavior of the displacement and velocity of the tip of the beam for the case of a light beam with $\rhos=1$.  Solutions are computed using the AMP scheme with composite grids $\Gbcf^{(j)}$, $j=2,4,8$, to indicate the grid convergence.  We also show the behavior of the displacement and velocity computed using the TP-SI scheme for a fine grid $\Gbcf^{(8)}$.  The zoomed views show excellent grid convergence for the AMP scheme, and excellent agreement with the solution obtained using the TP-SI scheme. 
We note that the TP-SI scheme required, on average, approximately $19$ sub-iterations per time-step using
grid $\Gbcf^{(4)}$,  and $6$ sub-iterations using grid $\Gbcf^{(8)}$, with an iteration tolerance of $10^{-6}$ in both cases.

%
{
\newcommand{\figWidth}{7.5cm}
\newcommand{\figWidthb}{3.75cm}
\newcommand{\trimfig}[2]{\trimFig{#1}{#2}{.0}{0.0}{.0}{.0}}
\begin{figure}[htb]
\begin{center}
\begin{tikzpicture}[scale=1]
  \useasboundingbox (0.0,.75) rectangle (16.,6.75);  
  \draw(0.0 ,0.0) node[anchor=south west,xshift=-4pt,yshift=+0pt] {\trimfig{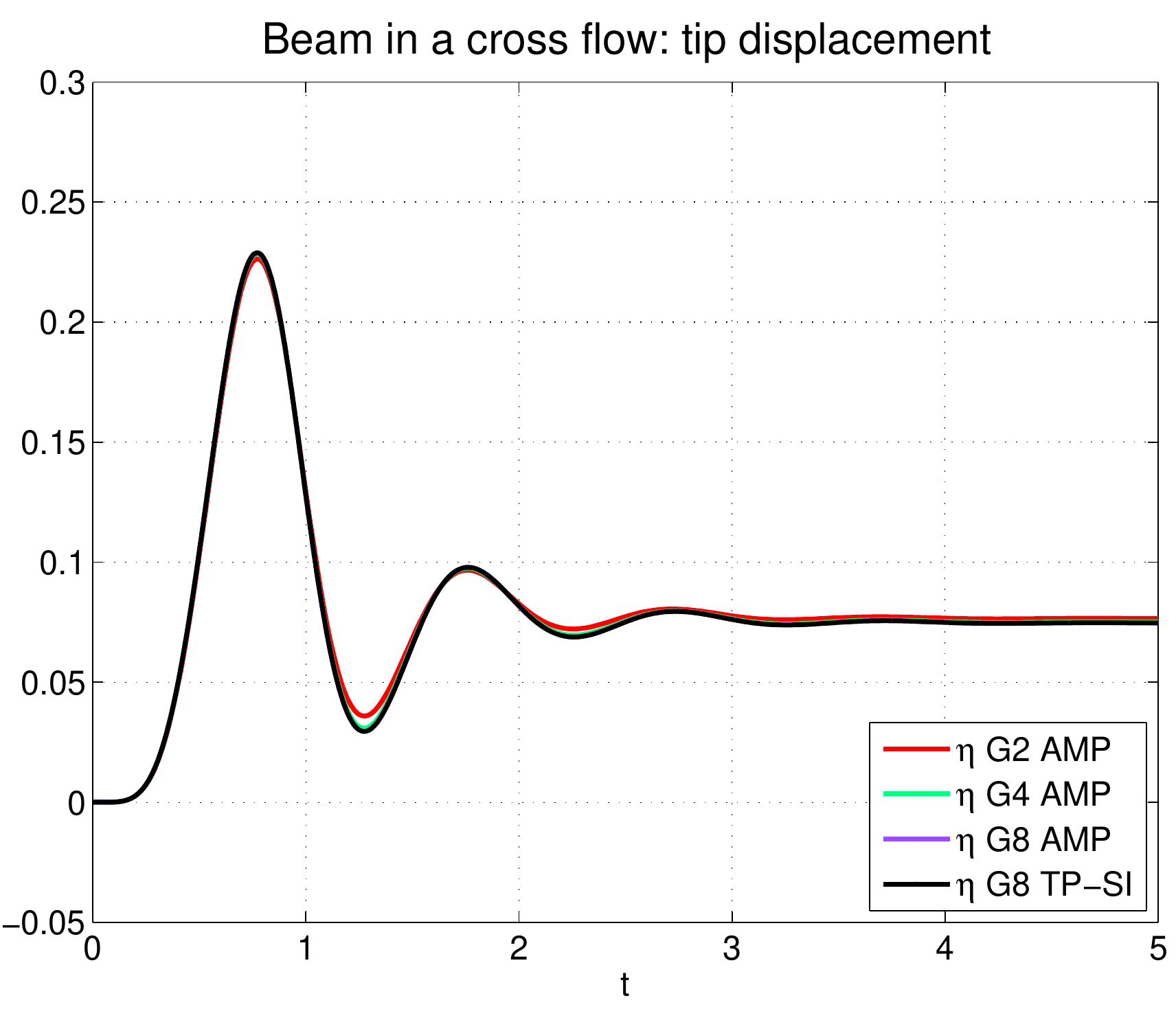}{\figWidth}};
  \draw(7.75,0.0) node[anchor=south west,xshift=-4pt,yshift=+0pt] {\trimfig{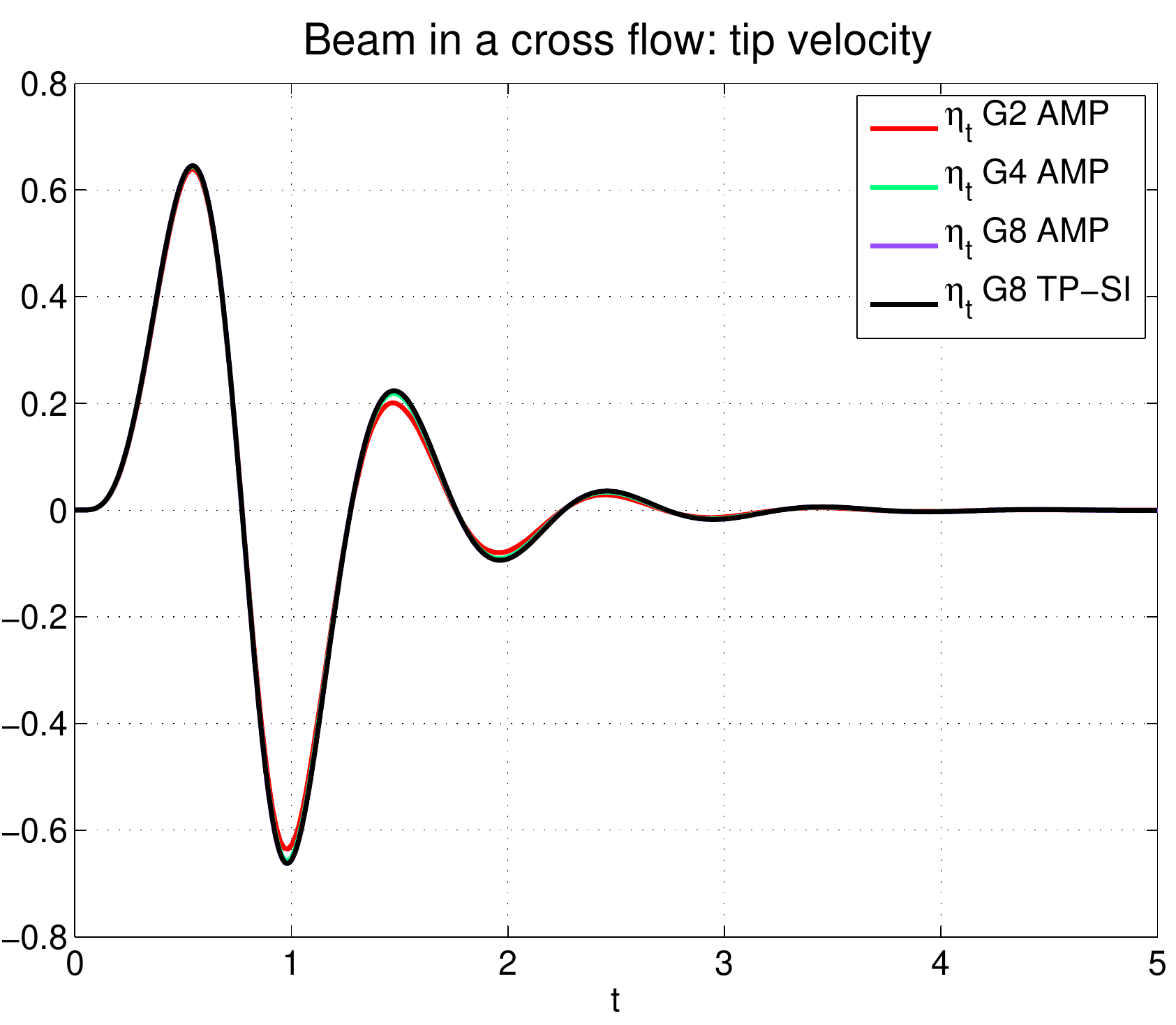}{\figWidth}};
  \draw(3.2,2.875) node[anchor=south west,xshift=-4pt,yshift=+0pt] {\trimfig{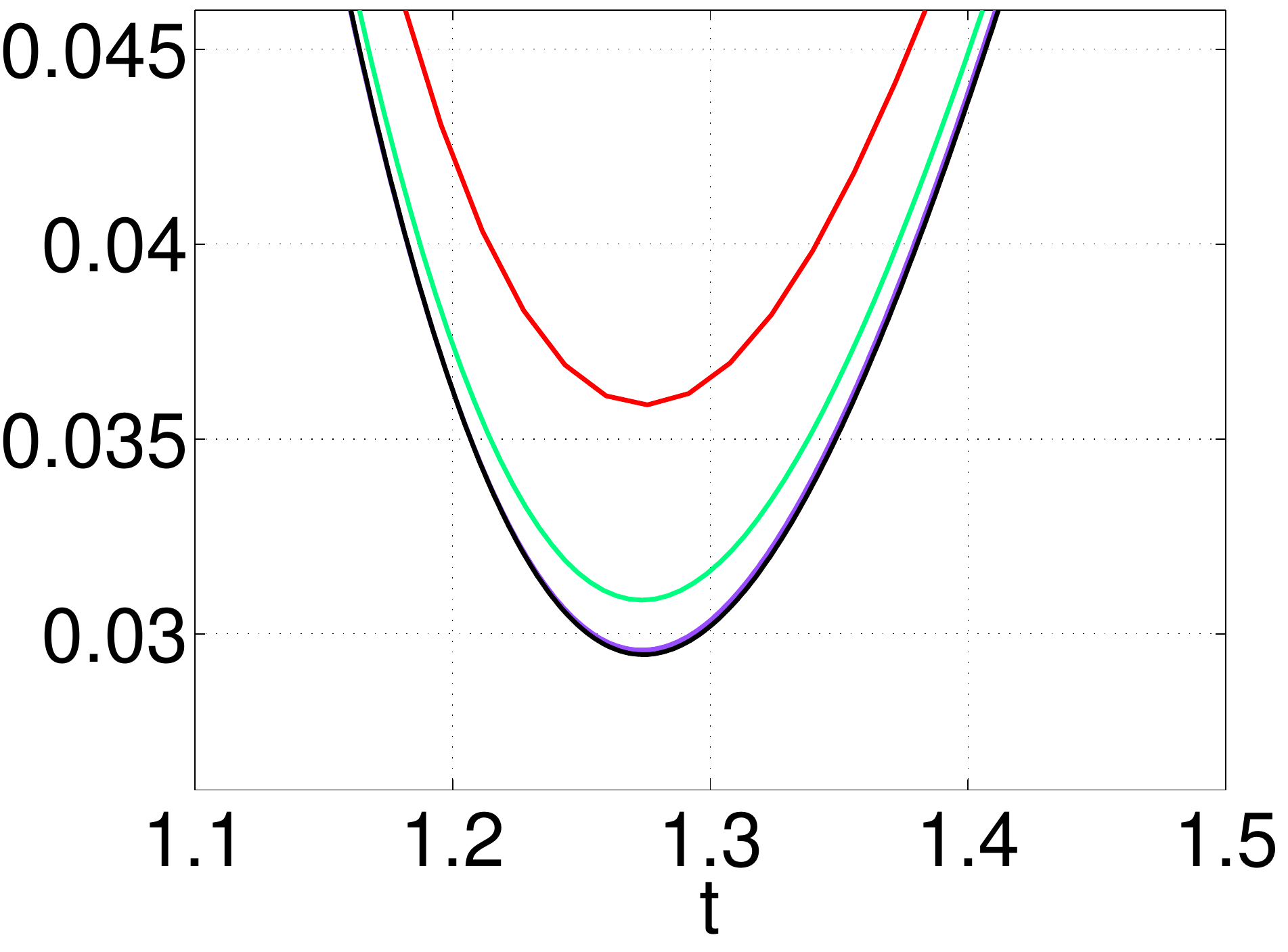}{\figWidthb}};
  \draw(11.75,.15) node[anchor=south west,xshift=-4pt,yshift=+0pt] {\trimfig{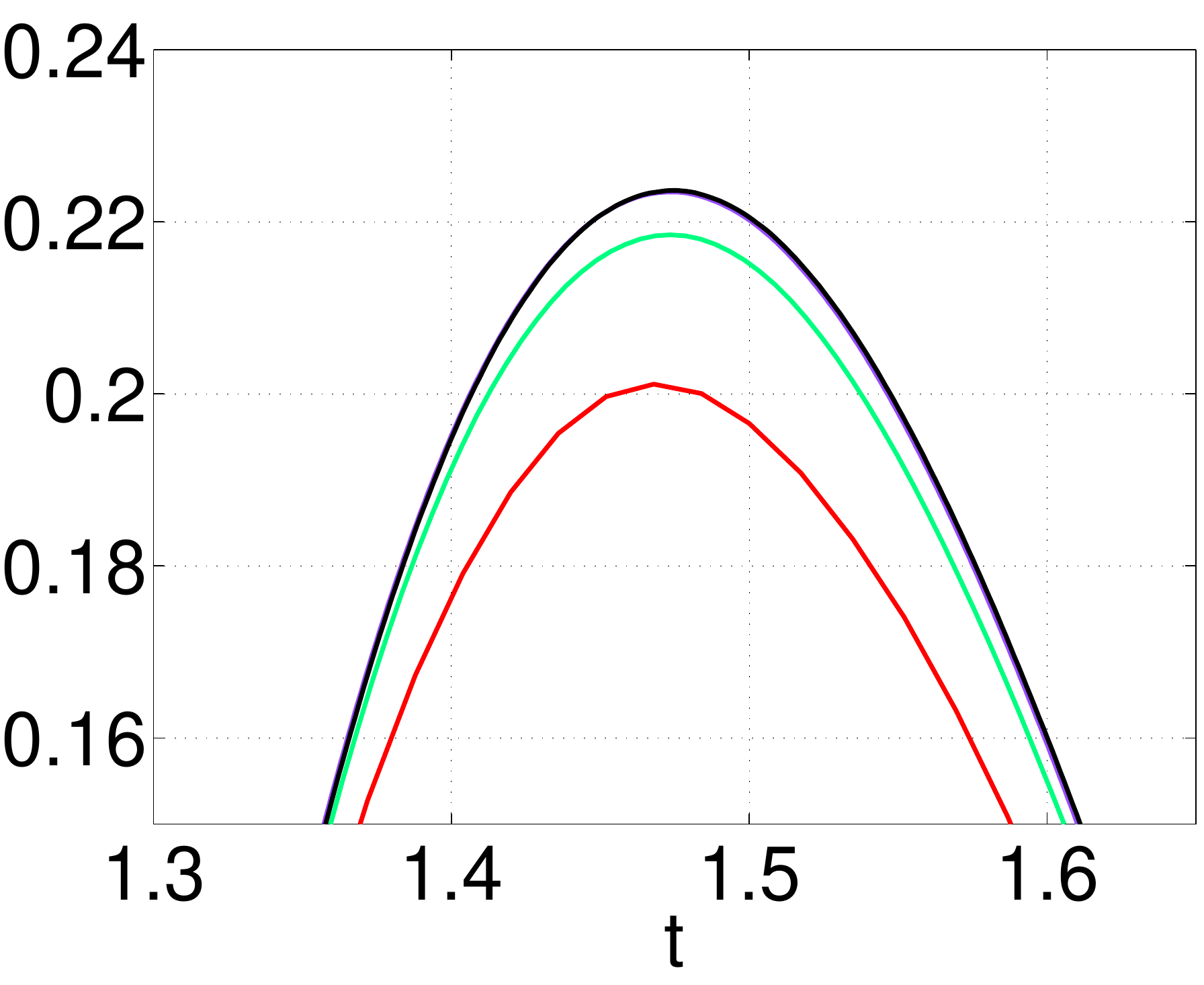}{\figWidthb}};
%
\end{tikzpicture}
\end{center}
\caption{Light beam in a cross flow: tip displacement and velocity as the grid is refined  for 
parameters $\Es\Is=5$, $\rhos=1$. }
\label{fig:lightBeamInAChannelTip}
\end{figure}
}

Figure~\ref{fig:lightBeamInAChannelTipComparison} shows the tip displacement and velocity for beams of different densities computed
using the composite grid $\Gbcf^{(4)}$.  We observe damped oscillations for all beams, with the damping for the heavy beam being the
weakest and its approximate period (time between successive peaks) being the largest.  We also note higher-frequency sub-harmonic
motion of the beam tip, especially for the heavy beam.  All solutions eventually approach a steady state where the displacement of the
beam tip is $\eta_{{\rm tip}} \approx 0.075$ independent of the beam density.

In the absence of the fluid, the frequency, $\bar\omega_0$, and period, $\bar\tau_0$,
of the lowest mode of a cantilevered beam are well known~\cite{Graff1991}
and given by
\begin{align}
 &   \bar\omega_0 \approx (1.875)^2 \sqrt{ \frac{ \Es\Is}{\rhos\hs\, \ls^4}}, \qquad 
   \bar\tau_0 = \frac{2\pi}{\bar\omega_0} .  \label{eq:cantileverBeam}
\end{align}
By comparing the computed frequency and period of the damped oscillating beam from the FSI simulation, 
the added-mass of the fluid can be estimated. In view of the
analysis discussed in Section~\ref{sec:analysis}
we expect that the frequency of the beam, when accounting for added mass effects, would be approximately given by
\begin{align}
 &   \bar\omega_{AM} \approx (1.875)^2 \sqrt{ \frac{ \Es\Is}{(\rhos\hs+M_a) \, \ls^4}}, \qquad 
   \bar\tau_{AM} = \frac{2\pi}{\bar\omega_{AM}},
     \label{eq:cantileverBeamAddedMass}
\end{align}
where $M_a$ is the added mass of the fluid (see equation~\eqref{eq:beamMotionFT}). 
 Thus, using~\eqref{eq:cantileverBeam} and~\eqref{eq:cantileverBeamAddedMass}, we have the estimate
\begin{equation}
    M_a = \rhos\hs \bigl\{ (\bar\tau_{AM}/{\bar\tau_0})^2 -1 \bigr\}.
\label{eq:AMestimate}
\end{equation}
For $\rhos=100$, $\hs=0.2$, $\ls=1$ and $\Es\Is=5$, we find $\bar\tau_0\approx 3.57$ which is close to the 
period estimated from Figure~\ref{fig:lightBeamInAChannelTipComparison} of $\bar\tau_{AM}\approx 3.6$.  Using~\eqref{eq:AMestimate}, we find $M_a/(\rhos\hs)\approx .02$, which is small as expected for this heavy beam. 

For a lighter beam with $\rhos=1$, we obtain $\bar\tau_0\approx 0.357$ and $\bar\tau_{AM}\approx .99$, which
gives $M_a\approx 1.3$, and this value is approximately $6.5$ times larger than the beam mass of $\rhos\hs=0.2$.
Based on the earlier mode analysis and the stability result in~\eqref{eq:traditionalStabLimitII}, we can estimate that the traditional scheme would become unconditionally unstable for $\rhos=M_a/\hs=1.3/0.2 = 6.5$ for $\dt$ approaching zero.
Computations of the full problem indicate that the TP scheme actually becomes unstable at a larger value of $\rhos \approx 13$. This larger value is expected since the time-step for this computation is chosen based on the stability of the domain solvers, as might be done in practice, and not reduced to account for added-mass effects.

{
\newcommand{\figWidth}{8cm}
\begin{figure}
\begin{centering}
\includegraphics[width=\figWidth]{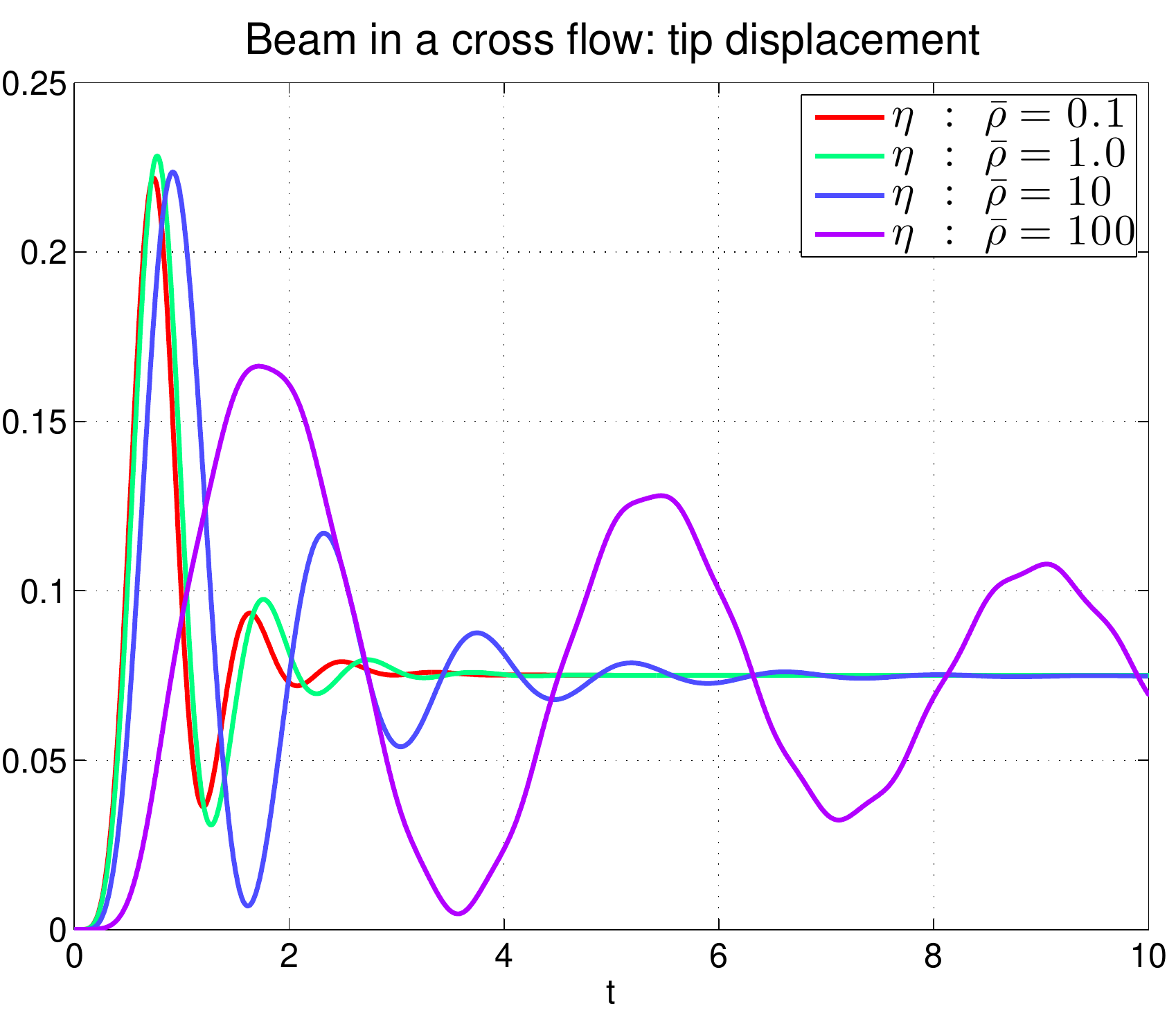} 
\includegraphics[width=\figWidth]{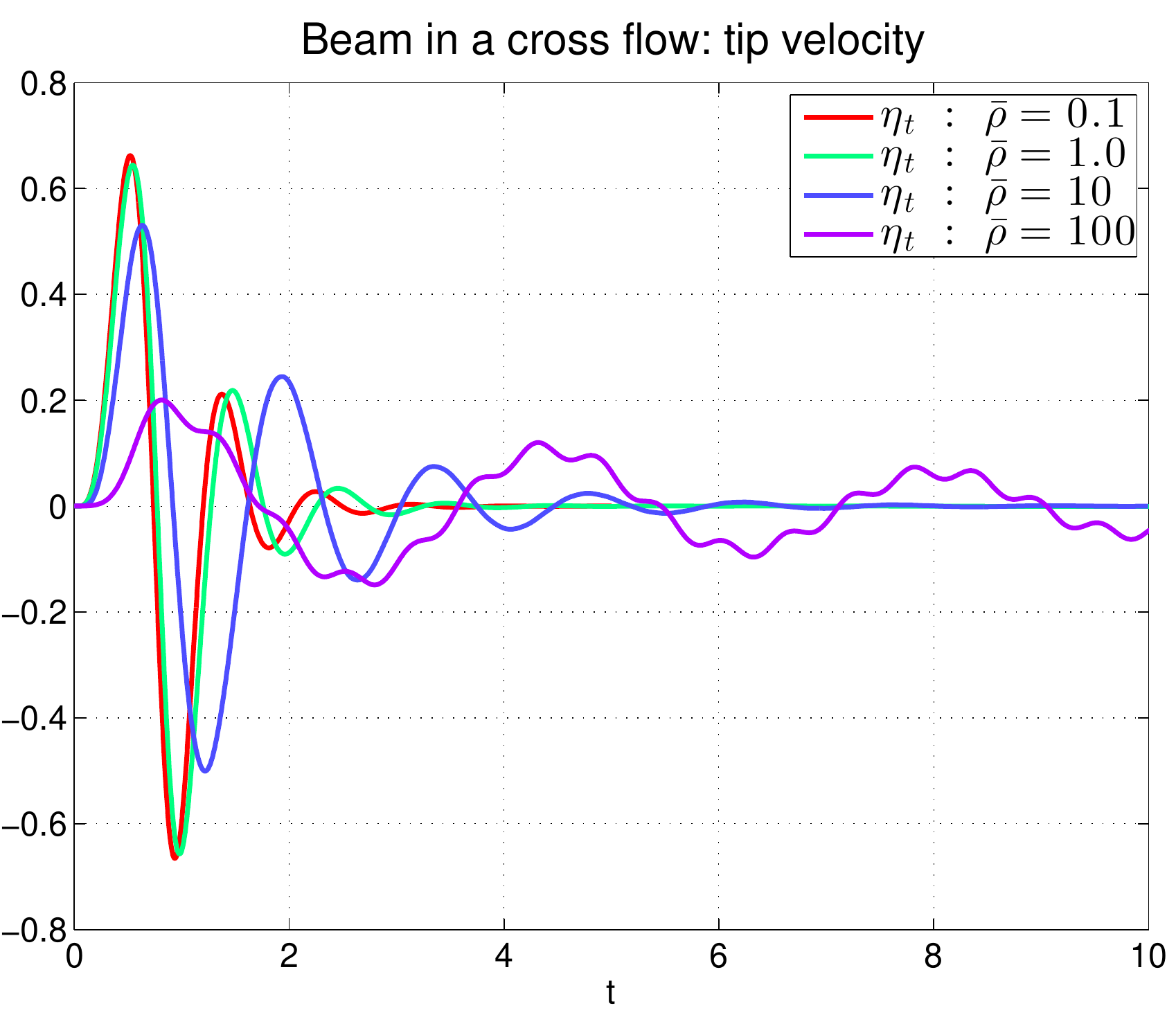}
\end{centering}

\caption{Beam in a cross flow: comparison of beam tip displacement and velocity for beams
of different densities. The solution was computed on grid $\Gbcf^{(4)}$ with $\Es\Is=5$.}
\label{fig:lightBeamInAChannelTipComparison}
\end{figure}
}

\section{Conclusions} \label{sec:conclusions}

This article makes a significant contribution to resolving an important and long-standing difficulty
for fluid-structure interaction (FSI) problems involving incompressible fluids and thin
structures. Such FSI problems arise in many applications in science and engineering, and their
accurate solution is an important and active field of research.  In a pair of recent
papers~\cite{fib2014,fis2014}, we have developed new {\em added-mass partitioned} (AMP) algorithms
for linearized FSI problems coupling incompressible (Stokes) flows and linearly elastic bulk solids
and thin structures.  An important feature of the AMP algorithms is that they are stable, without
sub-time-step iterations, independent of added-mass effects.
The present paper represents a first and important extension of the AMP algorithm in~\cite{fis2014} to fully
nonlinear FSI problems for the case of thin structures. The AMP interface conditions and 
time-stepping algorithm has been
fully described for this class of FSI problems.  The key AMP interface condition is a Robin (mixed) boundary
condition for the fluid pressure, and for a general beam this becomes a generalized Robin condition
that couples points on the circumference of the beam cross-section. 
For a two-dimensional beam with fluid on two sides, for example,
this condition involves the pressure on opposite sides of the beam.
 While the AMP algorithm has
been described for general beams, the focus has been on beams in two dimensions described by the
Euler-Bernoulli (EB) equations.  A normal-mode stability analysis has been performed for a linearized
problem involving an EB beam separating fluid domains on either side.  The analysis extends that
given in~\cite{fis2014} and shows that the AMP scheme is insensitive to
added-mass effects for this more general problem configuration, while a traditional partitioned
scheme (without sub-time-step iterations) becomes unconditionally unstable when added-mass effects
are sufficiently strong.

The general numerical framework for handling finite-amplitude deformations of the structure is based
on the use of moving and deforming composite grids, previously used in the context FSI problems
involving {\em compressible} flows coupled to rigid solids~\cite{mog2006,lrb2013}, and also to linear~\cite{fsi2012}
and nonlinear~\cite{flunsi2014r} elastic solids.  The equations governing the incompressible fluid
are solved using a fractional-step finite-difference scheme on deforming grids, while the equations
for the EB beam are solved using either finite differences or a finite-element approach.  For the
case of the finite-element beam solver, special treatment of the AMP interface conditions are
required in order to couple to the finite-difference fluid solver and we have described two
approaches for implementing this coupling.

An important feature of this paper is the use of benchmark problems of increasing complexity to
illustrate the behavior of the new AMP algorithm.  In addition to carefully verifying the stability
and second-order accuracy of the present algorithm, the results of the benchmark problems can be useful to other
researchers developing other FSI algorithms.  The AMP scheme was first verified using
manufactured solutions for the case of a deforming beam with a fluid domain on one side.  The
solution was next computed for a problem (with known steady state solution) of a deformed beam next
to one or two pressurized fluid chambers, and for a time-dependent oscillating flat beam (with known
solution).  Results were obtained using a traditional-partitioned (TP) scheme with sub-time-step iterations and
these agreed with those obtained using the AMP scheme.
 The AMP scheme was also applied to a problem from the literature for a simplified model of
an artery: solutions were computed for a traveling wave pressure-pulse propagating through a
deforming channel. Grid convergence studies were performed and a comparison was made to prior
results available in the literature.  Lastly the motion of a cantilevered beam in a fluid channel was
considered, and solutions for beams with different densities were compared and these results were
used to estimate the added mass.
In all benchmark cases considered the numerical solutions from the AMP algorithm
remained accurate and stable even for the case of light beams when added-mass effects are large.

In future work we will consider more general nonlinear beam models as well as FSI problems 
involving beams, tubes, and shells in three space dimensions.

\bibliographystyle{elsart-num}
\bibliography{journal-ISI,jwb,henshaw,henshawPapers,fsi}

\end{document}